\documentstyle[12pt,epsf,stmaryrd,amscd,amstex,amssymb]{amsart}
\newcommand{\nc}{\newcommand}
\newcommand{\rc}{\renewcommand}

\nc{\ad}{{\mbox{\bf{ad}}}}
\nc{\AJ}{{\operatorname{aj}}}
\nc{\Aut}{{\operatorname{Aut}}}
\nc{\Bls}{{{\cal B}ls}}
\nc{\Boxtimes}{{\fbox{$\times$}}}
\nc{\blt}{{\bullet}}
\nc{\bSt}{{\mbox{\bf{St}}}}
\nc{\card}{{\operatorname{card}}}
\nc{\Cch}{{\check{C}}}
\nc{\Ch}{{\operatorname{Ch}}}
\nc{\chara}{{\operatorname{char}}}
\nc{\CHom}{{\cal{H}{\mathrm{om}}}}
\nc{\Coker}{{\operatorname{Coker}}}
\nc{\codim}{{\operatorname{codim}}}
\nc{\Cone}{{\operatorname{Cone}}}
\nc{\cSgn}{{\cal{S}gn}}
\nc{\depth}{{\operatorname{depth}}}
\nc{\dirlim}{{\underset{\rightarrow}{\operatorname{lim}}}}
\nc{\dotbox}{{\overset{\bullet}{\boxtimes}}}
\nc{\dotimes}{{\overset{\bullet}{\otimes}}}
\nc{\Ed}{{\operatorname{Edge}}}
\nc{\Ext}{{\operatorname{Ext}}}
\nc{\Fac}{{\cal{F}ac}}
\nc{\Fun}{{\operatorname{F}}}
\nc{\FS}{{\cal{FS}}}
\nc{\Hom}{{\mathrm{Hom}}}
\nc{\isla}{{\mathrm{isom}}_\lambda}
\nc{\cone}{{\mathrm{cone}}}
\nc{\had}{{{\hat{\mbox{\bf{ad}}}}}}
\nc{\hgt}{{\operatorname{ht}}}
\nc{\Id}{{\operatorname{Id}}}
\nc{\id}{{\mathrm{id}}}
\nc{\Ima}{{\operatorname{Im}}}
\nc{\ind}{{\operatorname{ind}}}
\nc{\Ind}{{\mathrm{Ind}}}
\nc{\infi}{{\operatorname{inf}}}
\nc{\infh}{{\frac{\infty}{2}}}
\nc{\invlim}{{\underset{\leftarrow}{\operatorname{lim}}}}
\nc{\Jac}{{{\cal J}ac}}
\nc{\Ker}{{\operatorname{Ker}}}
\nc{\lcm}{{\operatorname{lcm}}}
\nc{\Locsys}{{{\cal L}ocsys}}
\nc{\Map}{{{\cal M}ap}}
\nc{\Mor}{{\operatorname{Mor}}}
\nc{\MS}{{\cal{MS}}}
\nc{\Ob}{{\operatorname{Ob}}}
\nc{\opp}{{\operatorname{opp}}}
\nc{\Or}{{{\cal O}r}}
\nc{\Ord}{{{\cal O}rd}}
\nc{\Part}{{{\cal P}art}}
\nc{\PGL}{{\operatorname{PGL}}}
\nc{\Pic}{{\mathrm{Pic}}}
\nc{\Rep}{{{\cal{R}}ep}}
\nc{\rk}{{\operatorname{rk}}}
\nc{\Sets}{{{\cal{S}}ets}}
\nc{\Sew}{{{\cal{S}}ew}}
\nc{\sgn}{{\operatorname{sgn}}}
\nc{\Sh}{{{\cal S}h}}
\nc{\Sign}{{{\cal S}ign}}
\nc{\Spe}{{\mbox{\bf{Sp}}}}
\nc{\supr}{{\operatorname{sup}}}
\nc{\Supp}{{\operatorname{Supp}}}
\nc{\supp}{{\operatorname{supp}}}
\nc{\Teich}{{{\cal{T}}eich}}
\nc{\tFS}{{\widetilde{\cal{FS}}}}
\nc{\Tor}{{\operatorname{Tor}}}
\nc{\totimes}{{\tilde{\otimes}}}
\nc{\tr}{{\operatorname{tr}}}
\nc{\tRep}{{\widetilde{{\cal R}ep}}}
\nc{\tTeich}{{\widetilde{{\cal T}eich}}}
\nc{\Vect}{{ {\operatorname {Vect} } }}
\nc{\Ve}{{\operatorname{Vert}}}
\nc{\wt}{{\widetilde}}

\nc{\bo}{{\mbox{\bf{0}}}}
\nc{\One}{{\mbox{\bf{1}}}}
\nc{\one}{{\mbox{\bf{1}}}}

\nc{\BA}{{\Bbb A}}
\nc{\ba}{{\mbox{\bf{a}}}}
\nc{\baB}{{\overline{B}}}
\nc{\baeta}{{\bar{\eta}}}
\nc{\baJ}{{\bar{J}}}
\nc{\BB}{{\Bbb B}}
\nc{\bA}{{\mbox{\bf{A}}}}
\nc{\bB}{{\mbox{\bf{B}}}}
\nc{\bC}{{\mbox{\bf{C}}}}
\nc{\bc}{{\mbox{\bf{c}}}}
\nc{\BC}{{\Bbb{C}}}
\nc{\bCC}{{\overline{\cal{C}}}}
\nc{\bCM}{{\overline{\cal{M}}}}
\nc{\bD}{{\bar{D}}}
\nc{\BD}{{\overline{D}}}
\nc{\bd}{{\mbox{\bf{d}}}}
\nc{\BE}{{\overline{E}}}
\nc{\BF}{{\overline{F}}}
\nc{\bF}{{\mbox{\bf{F}}}}
\nc{\bg}{{\mbox{\bf{g}}}}
\nc{\bG}{{\mbox{\bf{G}}}}
\nc{\BG}{{\Bbb G}}
\nc{\bGamma}{{\overline{\Gamma}}}
\nc{\bbH}{{\bar{\mbox{\bf{H}}}}}
\nc{\bH}{{\mbox{\bf{H}}}}
\nc{\bI}{{\mbox{\bf{I}}}}
\nc{\bII}{{ \bbb{\II}   }}
\nc{\bL}{{\mbox{\bf{L}}}}
\nc{\BL}{{\Bbb{L}}}
\nc{\blambda}{{\bar{\lambda}}}
\nc{\bM}{{\mbox{\bf{M}}}}
\nc{\bmu}{{\vec{\mu}}}
\nc{\bN}{{\mbox{\bf{N}}}}
\nc{\BN}{{\Bbb{N}}}
\nc{\bnu}{{\mbox{\boldmath{${\nu}$}}}}
\nc{\bof}{{\mbox{\bf{f}}}}
\nc{\BP}{{\Bbb P}}
\nc{\bP}{{\mbox{\bf{P}}}}
\nc{\BPO}{{\overset{\circ}{\BP}}}
\nc{\BQ}{{\Bbb Q}}
\nc{\bq}{{\mbox{\bf{q}}}}
\nc{\BR}{{\Bbb{R}}}
\nc{\bR}{{\mbox{\bf{R}}}}
\nc{\br}{{\mbox{\bf{r}}}}
\nc{\breta}{{\bar{\eta}}}
\nc{\bs}{{\mbox{\bf{s}}}}
\nc{\bS}{{\mbox{\bf{S}}}}
\nc{\bt}{{\mbox{\bf{t}}}}
\nc{\bU}{{\mbox{\bf{U}}}}
\nc{\bV}{{\mbox{\bf{V}}}}
\nc{\bu}{{\mbox{\bf{u}}}}
\nc{\BUpsilon}{{\bar{\Upsilon}}}
\nc{\bw}{{\mbox{\bf{w}}}}
\nc{\bx}{{\mbox{\bf{x}}}}
\nc{\bX}{{\mbox{\bf{X}}}}
\nc{\BZ}{{\Bbb{Z}}}
\nc{\bz}{{\mbox{\bf{z}}}}
\nc{\bZ}{{\mbox{\bf{Z}}}}
\nc{\bzero}{\mbox{\boldmath{$0$}}}

\nc{\CA}{{\cal A}}
\nc{\CAD}{{\overset{\bullet}{\cal{A}}}}
\nc{\CAO}{{\overset{\circ}{\cal{A}}}}
\nc{\CB}{{\cal B}}
\nc{\CalD}{{\cal D}}
\nc{\CE}{{\cal E}}
\nc{\CF}{{\cal F}}
\nc{\CG}{{\cal G}}
\nc{\CH}{{\cal H}}
\nc{\CI}{{\cal I}}
\nc{\CID}{{\overset{\bullet}{\cal{I}}}}
\nc{\CJ}{{\cal J}}
\nc{\CK}{{\cal K}}
\nc{\CL}{{\cal L}}
\nc{\CM}{{\cal M}}
\nc{\CN}{{\cal N}}
\nc{\CO}{{\cal O}}
\nc{\CP}{{\cal P}}
\nc{\CPO}{{\overset{\circ}{\cal{P}}}}
\nc{\CQ}{{\cal Q}}
\nc{\CR}{{\cal R}}
\nc{\CS}{{\cal S}}
\nc{\CT}{{\cal T}}
\nc{\CTD}{{\overset{\bullet}{\cal{T}}}}
\nc{\CTPO}{{\overset{\circ}{\cal{T}\cal{P}}}}
\nc{\CU}{{\cal{U}}}
\nc{\CV}{{\cal V}}
\nc{\CW}{{\cal W}}
\nc{\CX}{{\cal X}}
\nc{\CY}{{\cal Y}}
\nc{\CZ}{{\cal Z}}

\nc{\dCL}{{\overset{\bullet}{\cal{L}}}}
\nc{\dd}{{\operatorname{d}}}
\nc{\ddelta}{{\overset{\bullet}{\delta}}}
\nc{\dfu}{{\overset{\bullet}{\frak{u}}}}
\nc{\dlambda}{{\overset{\bullet}{\lambda}}}
\nc{\DO}{{\overset{\circ}{D}}}
\nc{\dpar}{{\partial}}
\nc{\dS}{{\overset{\bullet}{S}}}
\nc{\dT}{{\overset{\bullet}{T}}}

\nc{\hCH}{{\hat{\cal{H}}}}
\nc{\hCI}{{\hat{\cal{I}}}}
\nc{\hfC}{{\hat{\frak{C}}}}
\nc{\hfg}{{\hat{\frak{g}}}}
\nc{\hL}{{\hat{L}}}
\nc{\OH}{{\overset{\circ}{H}}}
\nc{\hpsi}{{\hat{\psi}}}
\nc{\hx}{{\hat{x}}}

\nc{\jo}{{\overset{\circ}{j}}}

\nc{\phid}{{\overset{\bullet}{\phi}}}

\nc{\tA}{{\tilde{A}}}
\nc{\ta}{{\tilde{a}}}
\nc{\tB}{{\tilde{B}}}
\nc{\tb}{{\tilde{b}}}
\nc{\tBP}{{\tilde{\BP}}}
\nc{\tC}{{\tilde{C}}}
\nc{\tc}{{\tilde{c}}}
\nc{\tCA}{{\tilde{\cal{A}}}}
\nc{\tCC}{{\tilde{\cal{C}}}}
\nc{\tCH}{{\tilde{\cal{H}}}}
\nc{\tCI}{{\tilde{\cal{I}}}}
\nc{\tCO}{{\tilde{\cal{O}}}}
\nc{\tCP}{{\tilde{\cal{P}}}}
\nc{\tCT}{{\tilde{\cal{T}}}}
\nc{\tD}{{\tilde{D}}}
\nc{\tDelta}{{\tilde{\Delta}}}
\nc{\tE}{{\tilde E}}
\nc{\tF}{{\tilde F}}
\nc{\tfD}{{\tilde{\frak{D}}}}
\nc{\tfF}{{\tilde{\frak{F}}}}
\nc{\tff}{{\tilde{\frak{f}}}}
\nc{\tfu}{{\tilde{\frak{u}}}}
\nc{\tJ}{{\tilde{J}}}
\nc{\tj}{{\tilde{j}}}
\nc{\tK}{{\tilde K}}
\nc{\tL}{{\tilde{L}}}
\nc{\tM}{{\tilde{M}}}
\nc{\tP}{{\tilde{P}}}
\nc{\tPhi}{{\tilde{\Phi}}}
\nc{\tpi}{\tilde{\pi}}
\nc{\TPO}{{\overset{\circ}{T\BP}}}
\nc{\tR}{{\tilde{R}}}
\nc{\tS}{{\tilde S}}
\nc{\tT}{{\tilde{T}}}
\nc{\ttau}{{\tilde{\tau}}}
\nc{\ttheta}{{\tilde{\theta}}}
\nc{\tU}{{\tilde{U}}}
\nc{\tUpsilon}{{\tilde{\Upsilon}}}
\nc{\tW}{{\tilde W}}
\nc{\ty}{{\tilde y}}
\nc{\tY}{{\tilde Y}}
\nc{\txi}{{\tilde{\xi}}}

\nc{\UD}{{\overset{\bullet}{U}}}
\nc{\UO}{{\overset{\circ}{U}}}

\nc{\vA}{{\vec{A}}}
\nc{\valpha}{{\vec{\alpha}}}
\nc{\vbeta}{{\vec{\beta}}}
\nc{\vc}{{\vec{c}}}
\nc{\vD}{{\vec{D}}}
\nc{\vd}{{\vec{d}}}
\nc{\vgamma}{{\vec{\gamma}}}
\nc{\vK}{{\vec{K}}}
\nc{\vlambda}{{\vec{\lambda}}}
\nc{\vmu}{{\vec{\mu}}}
\nc{\vnu}{{\vec{\nu}}}
\nc{\vo}{{\vec{0}}}
\nc{\vu}{{\vec{u}}}
\nc{\vx}{{\vec{x}}}
\nc{\vy}{\vec{y}}
\nc{\vzero}{\vec{0}}

\nc{\XO}{{\overset{\circ}{X}}}

\nc{\ya}{{\operatorname{aj}}}

\nc{\nen}{\newenvironment}
\nc{\ol}{\overline}
\nc{\ul}{\underline}

\nc{\Lra}{\Longrightarrow}

\nc{\Llra}{\Longleftrightarrow}
\nc{\hra}{\hookrightarrow}
\nc{\iso}{\overset{\sim}{\lra}}
\nc{\rlh}{\rightleftharpoons}

\nc{\IC}{{\cal{IC}}}
\nc{\PS}{{\cal{PS}}}
\nc{\oCG}{{\overline{\cal G}}}
\nc{\oCQ}{{\overline{\cal Q}}}
\nc{\oCZ}{{\overline{\cal Z}}}
\nc{\dZ}{{\overset{\bullet}{\cal Z}}{}}
\nc{\ddZ}{{\ddot{\cal Z}}{}}
\nc{\oZ}{{\overset{\circ}{\cal Z}}{}}

\nc{\dP}{{\overset{\bullet}{\cal P}}{}}
\nc{\oP}{{\overset{\circ}{\cal P}}{}}
\nc{\oQ}{{\overset{\circ}{\cal Q}}{}}
\nc{\obp}{{\overset{\circ}{{\bf p}}}{}}
\nc{\tbj}{{\tilde{\bf j}}{}}
\nc{\tbp}{{\tilde{\bf p}}{}}
\nc{\tfC}{{\widetilde{\frak C}}{}}
\nc{\tfE}{{\widetilde{\frak E}}{}}
\nc{\tfj}{{\widetilde{\frak j}}{}}
\nc{\tfQ}{{\widetilde{\frak Q}}{}}
\nc{\tfp}{{\widetilde{\frak p}}{}}
\nc{\ofQ}{{\overset{\circ}{{\frak Q}}}{}}
\nc{\tGQ}{{\widetilde{\cal{GQ}}}{}}
\nc{\oGQ}{{\overset{\circ}{\cal{GQ}}}{}}
\nc{\ooGQ}{{\overset{\circ\circ}{\cal{GQ}}}{}}
\nc{\oGZ}{{\overset{\circ}{\cal{GZ}}}{}}
\nc{\tGZ}{{\widetilde{\cal{GZ}}}{}}

\nc{\Ue}{{U_\varepsilon}}
\nc{\Upe}{{\Upsilon_\varepsilon}}
\nc{\crho}{{\check{\rho}}}
\nc{\ctheta}{{\check{\theta}}}

\nc{\pr}{\protect}

\nc{\nn}{{\newline}}
\nc{\np}{{\newpage}}    %\nc{\nnn}{{\newpage}}

\nc{\lab}{  \label}
\nc{\npp}{{ \newpage\setcounter{page}{0}    }}
\nc{\setpart}{{     \setcounter{part}   }}
\nc{\setpage}{{     \setcounter{page}   }}
\nc{\setsection}{{  \setcounter{section}    }}
\nc{\nd}{ $$\text{ This version is preliminary and approximate,
it is not for distribution. }$$ }

\nc{\noi}{{\noindent}}
\nc{\nop}{{\noindent {\em Proof.}} }
\nc{\cont}{\tableofcontents}

\nc{\sbr}{{\smallpagebreak}}
\nc{\mbr}{{\medpagebreak}}
\nc{\bbr}{{\bigpagebreak}}
\nc{\bbb}{ \boldsymbol         }
\nc{\bul}{ \bullet         }    %Bullet

\nc{\bem}{{ \begin{em}  }}
\nc{\eem}{{ \end{em}    }}

\nc{\bbox}{{    \blackbox   }}  %BLACKbox
\nc{\ra}{{  \rightarrow }}
\nc{\laa}{{ \leftarrow  }}  %Forget it and use @<<<
\nc{\lra}{{\longrightarrow}}
\nc{\lr}{{\leftrightarrow}}         % "corresponds"
\nc{\lrs}{{\rightleftarrows}}       % One map in each direction.

\nc{\imp}{{\Rightarrow}}            %Implication.
\nc{\eq}{{\Leftrightarrow}}         %Equivalence.
    \nc{\Ra}{{\Rightarrow}}             %Implication.
    \nc{\LRa}{{\Leftrightarrow}}            %Equivalence.
\nc{\inj}{{\pr  \hookrightarrow }}          %Injective map-right.
\nc{\injj}{{\pr \hookleftarrow  }}          %Injective map.

\nc{\sur}{{ \twoheadrightarrow  }}  %Surjective map-right.
\nc{\surr}{{    \twoheadleftarrow   }}  %Surjective map-left

\nc{\mm}{{\mapsto}}                     %Map on elements.
\nc{\mmm}{{\mapsfrom}}                  %Map on elements: to the left.

\nc{\va}{{\uparrow}}                    %Vertical Up-arrow
\nc{\bb}{\pr\underset}           %Write Bellow
\rc{\aa}{\pr\overset}            %and Above

\nc{\indd}{{ ${} \ \ \ \ \  \ \        {} $ }}  %Indent
\nc{\inddd}{{   \indd\indd          }}  %Indent 3times
\nc{\nnd}{{     \nn  \indd          }}  %Indent in a new line.
\nc{\nndb}{{    \nn  \indd $\bul$       }}  %Indent+Bullet in new.
\nc{\bss}{{\backslash}}                 %%% The other slash
\nc{\barr}{     \overline   }           %Long bar accent
\nc{\ud}{   \underline  }       %Underline

\nc{\ti}{\tilde}              %Tilde
\nc{\tii}{\widetilde}         %Tilde: wide

\nc{\hatt}{\widehat}                %Hat: wide
\nc{\hata}{{    \bbb{ \hat{} }      }}  %Hat: after the letter

\nc{\ch}{\check}                        %Check above the letter
\nc{\cha}{{     \bbb{ \check{} }    }}      %Check after the letter

\nc{\sub}{{ \subseteq   }}         %to the right
\nc{\subb}{{    \supseteq   }}         %       legt
\nc{\nsub}{{    \nsubseteq  }}         %Negated
\nc{\nsubb}{{   \nsupseteq  }}         %

\nc{\nin}{{ \notin  }}

\nc{\lb}{\langle}                           %Square brackets
\nc{\rb}{\rangle}

\nc{\lB}{   \left(  }                       %Oval   brackets
\nc{\rB}{   \right) }
\nc{\BBl}{{ \bbb{ \left( \right.}   }}              %Bold Oval brackets
\nc{\BBr}{{ \bbb{ \left. \right)}   }}

\nc{\mat} {     \left(      \matrix }
\nc{\emat}{     \endmatrix  \right) }
\nc{\sm} {      \left(      \smallmatrix    }
\nc{\esm}{      \endsmallmatrix \right) }
\nc{\smat} {        \left(      \smallmatrix    }
\nc{\esmat}{        \endsmallmatrix \right) }
\nc{\imat} {        \left.      \matrix }
\nc{\eimat}{        \endmatrix  \right. }
\nc{\ism} {     \left.      \smallmatrix    }
\nc{\eism}{     \endsmallmatrix \right. }

\nc{\ca}{       \left\{     \smallmatrix    }   %SMALL
\nc{\eca}{      \endsmallmatrix \right\}    }
\nc{\Ca}{       \left\{     \matrix     }
\nc{\eCa}{      \endmatrix  \right\}    }
\nc{\Eca}{      \endmatrix  \right.     }
\nc{\com}{  \begin{diagram} }
\nc{\ecom}{   \end{diagram} }

\nc{\tab}{  \begin{tabular}     }
\nc{\etab}{ \end{tabular}       }   %does not work?
\nc{\hl}{{  \hline          }}

\nc{\Eq}{   \begin{equation}    }
\nc{\Eeq}{  \end{equation}  }
\nc{\aln}{  \begin{align}   }
\nc{\ealn}{ \end{align} }

\nc{\se}{   \section        }
\nc{\sus}{  \subsection     }
\nc{\sss}{  \subsubsection      }

\nc{\lem}{  \subsubsection{\bf Lemma}       }
\nc{\slem}{     \subsubsection*{\bf Lemma}      }
\nc{\Lemm}{     \subsection{Lemma}      }
\nc{\lemm}{     \subsubsection{\bf Lemma}       }
\nc{\slemm}{    \subsubsection*{\bf Lemma}      }

\nc{\que}{  \subsubsection{\bf Question}        }
\nc{\sque}{     \subsubsection*{\bf  Question}      }

\nc{\sublem}{   \subsubsection{\bf Sublemma}        }
\nc{\ssublem}{  \subsubsection*{\bf Sublemma}       }

\nc{\Pro}{  \subsection{Proposition}    }
\nc{\pro}{  \subsubsection{\bf Proposition} }
\nc{\spro}{     \subsubsection*{\bf Proposition}    }

\nc{\Corr}{     \subsection{Corollary}      }
\nc{\corr}{     \subsubsection{\bf Corollary}   }
\nc{\scorr}{    \subsubsection*{\bf Corollary}  }

\nc{\Theo}{     \subsection{Theorem}        }
\nc{\theo}{     \subsubsection{\bf Theorem}     }
\nc{\stheo}{    \subsubsection*{\bf Theorem}    }

\nc{\rem}{  \subsubsection{Remark}      }
\nc{\srem}{     \subsubsection*{Remark} }
\nc{\remm}{     \subsubsection{Remarks}     }
\nc{\sremm}{    \subsubsection*{Remarks}    }
\nc{\rema}{     \subsubsection{Remarks}     }

\nc{\conj}{     \subsubsection{Conjecture}  }
\nc{\sconj}{    \subsubsection*{Conjecture} }

\nc{\ex}{   \subsubsection{Example}     }
\nc{\sex}{  \subsubsection*{Example}    }
\nc{\exs}{  \subsubsection{Examples}    }
\nc{\sexs}{     \subsubsection*{Examples}   }

\nc{\h}{{   \hslash }}  %Planck
\nc{\All}{{ \forall }}
\nc{\yy}{\infty}
\nc{\ys}{{  \frac{\infty}{2}  }}

\nc{\pl}{{\oplus}}                      %Direct sum
\nc{\tim}{{\times}}
\nc{\btim}{{\boxtimes}}
\nc{\ltim}{\ltimes}                     %
\nc{\rtim}{\rtimes}         %
\nc{\ltr}{\triangleleft}        %
\nc{\rtr}{\triangleright}       %

\nc{\ten}{{ \otimes     }}
\nc{\Lten}{{    \aa{L}\otimes   }}            %Derived tensoring.
\nc{\tenA}{ \bb{A}\ten  }
\nc{\tenB}{ \bb{B}\ten  }
\nc{\tenZ}{ \bb{\Z}\ten }
\nc{\tenR}{ \bb{\R}\ten }
\nc{\tenC}{ \bb{\C}\ten }
\nc{\tenk}{ \bb{\k}\ten }
\nc{\bten}{{\boxtimes}}                 %Tensoring: outer

\nc{\con}{{ @>{\protect\cong}>> }}      %Isomorphism with a right arrow
\nc{\conl}{{    @>{\cong}>> }}      %Lower right isomorphism
\nc{\conn}{{    @<{\cong}<<     }}      %Isomormphism with a left arrow
\nc{\Con}{{ \equiv      }}  %Congruence
\nc{\appr}{{    \sim        }}  %Approximately
\nc{\eqr}{{ \sim        }}  %Equivalence relation

\nc{\fra}{  \frac   }       %%%xxxxxxxx  Another name for fractions

\nc{\ha}{{ \frac{1}{2} }}           %%%Half
    \nc{\half}{{ \frac{1}{2} }}
\nc{\ci}{{\circ}}               %Circle dot
\nc{\cd }{{\cdot}}              %Multiplication dot
\nc{\cdd}{{\cdot}}              %Multiplication dot
        \nc{\cdx}{{\cdot}}            %old
\nc{\cddd}{{\cdot\cdot\cdot}}   %Three dots
\nc{\ox}{{  \OO_X       }}               %Functions on $X$

\nc{\omx}{{ \om_X       }}               %Canonical line bundle
\nc{\Omx}{{ \Om_X^1     }}               %1-forms

\nc{\cupp}{\bigcup}             %cup,cap
\nc{\capp}{\bigcap}
\nc{\pll}{\bigoplus}

\nc{\pii}{\prod}                %Product
\nc{\ppii}{\bigprod}            %Big product
\nc{\cci}{\sqcup}              %Coprodduct
\nc{\ccii}{\bigsqcup}

\nc{\wwe}{\bigwedge}            %wedge
\nc{\cce}{\bigcoprod}           %cowedge

\nc{\aaa}{  \stackrel   }

\nc{\edd}{{ \end{document}  }}
\nc{\tx}{   \text       }       %Text

\nc{\df}{{  \overset{\text{def}} = }}   %definition
\nc{\inv}{{     {}^{-1}             }}  %inverse
\nc{\thh}{  ^{\text{th}}        }       %nth

\nc{\emp}{{   \emptyset}}               %Emptyset

\nc{\we}{{\wedge}}              %Wedge
\nc{\wee}{{ \aa{\bul}\wedge }}      %All wedge powers
\nc{\wetwo}{{     \overset{2}\wedge       }}    %2nd wedge power

\nc{\limp}{{    \underset {\leftarrow}\lim  }}
\nc{\limi}{{    \underset {\rightarrow}\lim }}
\nc{\plimp}{{\pro\underset {\leftarrow}\lim }}
\nc{\plimi}{{\pro\underset {\rightarrow}\lim    }}

\nc{\ppp}{{ {\Bbb P}^1 }}                   %P1
\nc{\ppn}{{ {\Bbb P}^n }}                   %Pn
\nc{\pt}{   { \mathrm{pt} }   }       %point

\nc{\qlb}{{ \barr{{\Bbb Q}_l} }}            %Q-el-bar (l-adic)
\nc{\ffq}{{  {\Bbb F}_q  }}                 %Ef-q (finite field)
\nc{\ffp}{{  {\Bbb F}_p  }}                 %Ef-q (finite field)
\nc{\tw}{   {}^{(1)}    }       %Frobenius twist

\nc{\Spec}{{ \mathrm{Spec}            }}
\nc{\Specf}{{ \text{Specf}              }}
\nc{\aand}{{\ \ \ \text{and}\ \ \   }}
\nc{\oor}{{\ \  \text{or}\ \    }}
\nc{\hk}{{       \mathrm{ hyperk\ddot{a}hler }    }}

\rc{\Im}{{  \text{Im}   }}
\nc{\rank}{{    \ \text{rank}\  }}
\nc{\Res}{{ \  \mathrm{Res}   }}
\nc{\End}{{ \mathrm{End}  }}
\nc{\RHom}{{    \mathrm{RHom} }}
\nc{\HHom}{{    \HH{\mathrm{om}}  }}
\nc{\EEnd}{{    \EE{\mathrm{nd}} }}
\nc{\AAut}{{    \text{$\AA ut$} }}
\nc{\RHHom}{{   \mathrm{R}\HH\mathrm{om} }}
\nc{\Der}{{ \text{Der}  }}

\nc{\ord    }{{ \text{ord} }}           %order of zero
\nc{\divv   }{{ \text{div} }}           %divisor
\nc{\Lie    }{{ \mathrm{Lie} }}
\nc{\Tens}{{ \mathrm{Tens}  }}

\nc{\timA} {{   \underset{A}\tim             }}
\nc{\timB} {{   \underset{B}\tim             }}
\nc{\timC} {{   \underset{C}\tim             }}
\nc{\timG} {{   \underset{G}\tim             }}
\nc{\timH} {{   \underset{H}\tim             }}
\nc{\timN} {{   \underset{N}\tim             }}
\nc{\timP}{{    \underset{P}\tim             }}
\nc{\timQ}{{    \underset{Q}\tim             }}
\nc{\timS} {{   \underset{S}\tim             }}
\nc{\timT} {{   \underset{T}\tim             }}
\nc{\timU} {{   \underset{U}\tim             }}
\nc{\timV} {{   \underset{V}\tim             }}
\nc{\timX} {{   \underset{X}\tim             }}
\nc{\timY} {{   \underset{Y}\tim             }}
\nc{\timZ} {{   \underset{Z}\tim             }}

\nc{\ab}{{       ^{\text{ab}}           }}
\nc{\af}{{       ^{\mathrm{aff}}          }}
\nc{\aff}{{      _{\mathrm{aff}}          }}
\nc{\dom}{{      _{\mathrm{dom}}          }}
\nc{\AZ}{{      _{\mathrm{AZ}}          }}

\rc{\mod}{{\mathrm{mod}}}

\nc{\cod}{\text{codim}} %Codimension

\rc{\AA}{{\cal A}}
\rc{\BB}{{\cal B}}
\nc{\CC}{{\cal C}}
\nc{\DD}{{\cal D}}
\nc{\EE}{{\cal E}}
\nc{\FF}{{\cal F}}
\nc{\GG}{{\cal G}}
\nc{\HH}{{\cal H}}
\nc{\II}{{\cal I}}
\nc{\JJ}{{\cal J}}
\nc{\KK}{{\cal K}}
\nc{\LL}{{\cal L}}
\nc{\MM}{{\cal M}}
\nc{\NN}{{\cal N}}
\nc{\OO}{{\cal O}}
\nc{\PP}{{\cal P}}
\nc{\QQ}{{\cal Q}}
\nc{\RR}{{\cal R}}
\rc{\SS}{{\cal S}}
\nc{\TT}{{\cal T}}
\nc{\UU}{{\cal U}}
\nc{\VV}{{\cal V}}
\nc{\WW}{{\cal W}}
\nc{\ZZ}{{\cal Z}}
\nc{\XX}{{\cal X}}
\nc{\YY}{{\cal Y}}

\nc{\A}{{\Bbb A }}
\nc{\B}{{\Bbb B}}
\nc{\C}{{\Bbb C}}
        \nc{\cc}{{\Bbb C}}
\nc{\Cs}{{\Bbb C^*}}
        \nc{\cs}{{\Bbb C^*}}
        \nc{\ccs}{{\Bbb C^*}}
\nc{\D}{{\Bbb D}}
\nc{\E}{{\Bbb E}}
\nc{\F}{{\Bbb F}}
\nc{\G}{{\Bbb G}}
    \nc{\hH}{{\Bbb H}}
\nc{\I}{{\Bbb I}}
\nc{\J}{{\Bbb J}}
\nc{\K}{{\Bbb K}}
    \nc{\lL}{{\Bbb L}}
\nc{\M}{{\Bbb M}}
\nc{\N}{{\Bbb N}}
    \nc{\oO}{{\Bbb O}}
    \nc{\pP}{{\Bbb P}}
\nc{\Q}{{\Bbb Q}}
\nc{\R}{{\Bbb R}}
    \nc{\sS}{{\Bbb S}}
\nc{\T}{{\Bbb T}}
\nc{\U}{{\Bbb U}}
\nc{\V}{{\Bbb V}}
\nc{\W}{{\Bbb W}}
\nc{\Z}{{\Bbb Z}}
\nc{\X}{{\Bbb X}}
\nc{\Y}{{\Bbb Y}}

\nc{\k}{{\Bbbk}}
\let\O\oO

\nc{\fA}{{\frak A}}
\nc{\fB}{{\frak B}}
\nc{\fC}{{\frak C}}
\nc{\fD}{{\frak D}}
\nc{\fE}{{\frak E}}
\nc{\fF}{{\frak F}}
\nc{\fG}{{\frak G}}
\nc{\fH}{{\frak H}}
\nc{\fI}{{\frak I}}
\nc{\fJ}{{\frak J}}
\nc{\fK}{{\frak K}}
\nc{\fL}{{\frak L}}
\nc{\fM}{{\frak M}}
\nc{\fN}{{\frak N}}
\nc{\fO}{{\frak O}}
\nc{\fP}{{\frak P}}
\nc{\fQ}{{\frak Q}}
\nc{\fR}{{\frak R}}
\nc{\fS}{{\frak S}}
\nc{\fT}{{\frak T}}
\nc{\fU}{{\frak U}}
\nc{\fV}{{\frak V}}
\nc{\fW}{{\frak W}}
\nc{\fZ}{{\frak Z}}
\nc{\fX}{{\frak X}}
\nc{\fY}{{\frak Y}}
\nc{\fa}{{\frak a}}
\nc{\fb}{{\frak b}}
\nc{\fc}{{\frak c}}
\nc{\fd}{{\frak d}}
\nc{\fe}{{\frak e}}
\nc{\ff}{{\frak f}}
\nc{\fg}{{\frak g}}
\nc{\fh}{{\frak h}}
\nc{\fiI}{{\frak i}}  %!!!
    \nc{\ffi}{{\frak i}}  %!!!
\nc{\fj}{{\frak j}}
\nc{\fk}{{\frak k}}
\nc{\fl}{{\frak{l}}}
\nc{\fm}{{\frak m}}
\nc{\fn}{{\frak n}}
\nc{\fo}{{\frak o}}
\nc{\fp}{{\frak p}}
\nc{\fq}{{\frak q}}
\nc{\fr}{{\frak r}}
\nc{\fs}{{\frak s}}
\nc{\ft}{{\frak t}}
\nc{\fu}{{\frak u}}
\nc{\fv}{{\frak v}}
\nc{\fw}{{\frak w}}
\nc{\fz}{{\frak z}}
\nc{\fx}{{\frak x}}
\nc{\fy}{{\frak y}}

\nc{\al}{{\alpha }}
\nc{\be}{{\beta }}
\nc{\ga}{{\gamma }}
\nc{\de}{{\delta }}
\nc{\del}{{\partial }}
\nc{\ep}{{\varepsilon }}
\nc{\vap}{{\epsilon }}

\nc{\ze}{{\zeta }}
\nc{\et}{{\eta }}
\rc{\th}{{\theta }}

\nc{\io}{{\iota }}
\nc{\ka}{{\kappa }}
\nc{\la}{{\lambda }}
\nc{\vrho}{{\varrho}}
\nc{\si}{{\sigma }}
\nc{\ups}{{\upsilon }}
\nc{\vphi}{{\varphi }}
\nc{\om}{{\omega }}

\nc{\Ga}{{\Gamma }}
\nc{\De}{{\Delta }}
\nc{\nab}{{\nabla}}
\nc{\Th}{{\Theta }}
\nc{\La}{{\Lambda }}
\nc{\Si}{{\Sigma }}
\nc{\Ups}{{\Upsilon }}
\nc{\Om}{{  \Omega      }}
\nc{\toc}{{\tableofcontents}}
\nc{\addl}{ \addcontentsline{toc}{subsection}   }

\def\square{\hbox{\vrule\vbox{\hrule\phantom{o}\hrule}\vrule}}

\nc{\ii}{{  i\in I      }}
\nc{\tww}{{ {}^{*,1}    }}
\nc{\zhc}{{ \fZ_{{\mathrm{HC}}} }}
\nc{\zfr}{{ \fZ_{\mathrm{Fr}} }}
\nc{\fzp}{{ \fz_{p} }}
\nc{\pp}{ ^{[p]} }

\nc{\utx}{{     \DD_X       }}
\nc{\utb}{{     \DD_\BB     }}
\nc{\dx}{   \DD_X       }
\nc{\db}{   \DD_\BB     }
\nc{\z}{{   ^{\bbb 0}   }}
\nc{\LLoc}{{    \LL_{\hat 0}        }}
\let\Phi\LLoc
\nc{\hzc}{{ \hatt{(\chi,0)} }}
\nc{\hrc}{{ \hatt{\chi,-\rho}   }}
\nc{\mr}{{  {-\rho} }}

\nc{\ob}{{  \OO_\BB     }}               %Functions on $\BB$
\nc{\uc}{{  u_\chi      }}
\nc{\Fp}{{  {\Bbb F}_p  }}
\nc{\frx}{{ {\mathrm{Fr}}_X        }}

\nc{\Fnew}{ \mbox{Fr}_X }
\nc{\Frx}{  \Fnew       }

\theoremstyle{remark}

\nc{\proof}{{\it Proof. }}

\newcommand{\imbed}{\hookrightarrow}
\renewcommand{\iso}{{\tii \longrightarrow}}

\newcommand{\Lotimes}{\overset{\rm L}{\otimes}}

\newcommand{\oplusl}{\bigoplus\limits}
\newcommand{\cupl}{\bigcup\limits}

\def\square{\hbox{\vrule\vbox{\hrule\phantom{o}\hrule}\vrule}}

\renewcommand{\E}{{\mathcal E}}
\renewcommand{\V}{{\mathcal V}}
\renewcommand{\N}{{\mathcal N}}
\renewcommand{\O}{{\mathcal O}}

\renewcommand{\F}{{\mathcal F}}
\renewcommand{\G}{{\mathcal G}}

\renewcommand{\B}{{\mathcal B}}

\newcommand{\Db}{   {\mathrm{D}}^{b}    }
\newcommand{\Dmin}{ {\mathrm{D}}^-  }
\let\Dm\Dmin

\nc{\til}{  \tilde  }

\nc{\Mod}{{ \mathrm{mod} }}
\newcommand{\Wap}{{     W'\aff    }}

\newcommand{\g}{{\mathfrak g}}
\renewcommand{\h}{{\mathfrak h}}
\renewcommand{\b}{{\mathfrak b}}

\renewcommand{\bu}{{        \bullet     }}

\newcommand{\epf}{\square}

\newcommand{\Zet}{{\mathbb Z}}
\newcommand{\Pn}{{{\mathbb P}^n}}

\renewcommand{\sur}{\twoheadrightarrow}

\renewcommand{\M}{{\mathcal M}}
\renewcommand{\N}{{\mathcal N}}
\renewcommand{\B}{{\mathcal B}}
\renewcommand{\O}{{\mathcal O}}
\renewcommand{\F}{{\mathcal F}}

\newcommand{\gtil}{{\tii{\mathfrak g}}}

\nc{\diff}{{    \bbb{\bss}  }}

\nc{\hla}{{\hatt\la}}

\renewcommand{\Fp}{{{\mathbb F}_p}}

\nc{\FN}{{   \mathrm{FN}  }}
\rc{\U}{{   \mathrm{U}  }}
\nc{\uU}{{  \U      }}
\rc{\i}{\item}

\newcommand{\gt}{{\tilde \g }}

\nc{\Kbar}{{\bar{\KK}}}

\nc{\prr}{{\mathrm{pr}}}

\rc{\aa}{\pr\overset}            %and Above

\nc{\cor}{     \subsubsection{\bf Corollary}   }
\nc{\scor}{     \subsubsection*{\bf Corollary}   }
\nc{\rems}{     \subsubsection{Remarks}         }
\nc{\srems}{    \subsubsection*{Remarks}        }

\nc{\dff}{{\    \overset{\text{def}}=\  }}       %definition
\rc{\Im}{{      \text{Im}       }}

\rc{\AA}{{\cal A}}
\rc{\BB}{{\cal B}}
\rc{\SS}{{\cal S}}
\let\O\oO

\nc{\Dd}{{  \text{D}    }}
\nc{\Hh}{{  \text{H}    }}

\nc{\Ll}{{  \text{L}    }}

\nc{\Rr}{{  \mathrm{R}    }}
\rc{\th}{{\theta }}
\def\square{\hbox{\vrule\vbox{\hrule\phantom{o}\hrule}\vrule}}

\renewcommand{\iso}{{\tii \longrightarrow}}

\def\square{\hbox{\vrule\vbox{\hrule\phantom{o}\hrule}\vrule}}

\renewcommand{\N}{{\mathcal N}}

\renewcommand{\O}{{\mathcal O}}

\renewcommand{\F}{{\mathcal F}}
\renewcommand{\G}{{\mathcal G}}

\renewcommand{\B}{{\mathcal B}}
\nc{\hattt}{{               }}
\renewcommand{\h}{{\mathfrak h}}

\renewcommand{\bu}{{            \bullet         }}

\newcommand{\Ce}{{\mathbb C}}

\newcommand{\bbA}{{\mathbb A}}

\renewcommand{\M}{{\mathcal M}}
\renewcommand{\N}{{\mathcal N}}
\renewcommand{\B}{{\mathcal B}}
\renewcommand{\O}{{\mathcal O}}
\renewcommand{\F}{{\mathcal F}}

\renewcommand{\b}{{\mathfrak b}}
\nc{\Coh}{{ \CC{\mathrm{oh}}  }}
\nc{\bi}{   \begin{itemize}\item        }
\rc{\i}{    \item           }
\nc{\ei}{ \end{itemize} }
\nc{\ben}{  \begin{enumerate}\item      }
\nc{\een}{  \end{enumerate}         }

 \nc{\Iaff}{{    I\aff }}
\nc{\Waff}{{    W\aff          }} 
\nc{\Baff}{{    B\aff }}
\nc{\RGa}{{     \Rr\Gamma           }}

\nc{\uHom}{{ \underline{\rm Hom}   }}
\nc{\uuHom}{{ \bf Hom  }}

\nc{\PR}{{\bf \mathrm{pr}}}
\nc{\AS}{{\mathrm{AS}}}
\nc{\f}[1]{ \fbox{$\blacklozenge$\footnote{ \fbox{!}#1 }$\blacklozenge$ }   }
\nc{\bDD}{{ \bbb{\DD}   }}

\nc{\ie}{{,\ \     \text{i.e.,}\ \  }}
\nc{\iif}{{\ \     \text{if}\ \     }}
\nc{\hence}{{\ \ \ \text{hence}\ \ \    }}
\nc{\while}{{\ \ \ \text{while}\ \ \    }}
\nc{\with}{{\ \ \  \text{with}\ \ \     }}
\nc{\foor}{{\ \     \text{for}\ \   }}
\nc{\suchthat}{{\ \     \text{such that}\ \     }}

\nc{\ftt}[1]{{\footnote{#1}}}
\nc{\fttt}[1]{{$^($\footnote{#1}$^)$}}
\nc{\bftt}[1]{\footnote{#1}}

\nc{\st}{{  \tii s  }}

\theoremstyle{remark}

\renewcommand{\A}{{\mathcal A}}
\renewcommand{\O}{{\mathcal O}}
\newcommand{\Nt}{{\widetilde{\N}}}
\renewcommand{\b}{{\mathfrak b}}

\nc{\pmo}{{ \pm 1   }}

\nc{\unr}{  _{\mathrm{unr}}  }
\nc{\PPunr}{    _{\PP-{\mathrm{unr}}}  }
\nc{\bl}{{ \fbox{$\blacklozenge\to$} }}
\nc{\wl}{{ \fbox{$@<<<\lozenge$}     }}

\let\iso\con
\nc{\IIII}{I}   \let\Sigma\IIII

\nc{\mud}{{ \aa{\bu}\mu     }}
\nc{\tI}{{  \text{\bf I}        }}
\nc{\mn}{   _{\mu|\nu}      }
\nc{\cen}[1]{   \begin{center}  {\em  #1}   \end  {center}  }
\nc{\heart}{{\tiny \cen{\tiny $\heartsuit $ }   }} %Heart

\nc{\bOO}{{ \bbb{\OO}   }}
\nc{\bio}{{    \bbb{\io}   }}
\nc{\vpi}{{    \varpi	   }}

\nc{\WT}{{    \operatorname{wt} }} 

\newcommand{\RE}{{\mathbb R}}

\newcommand{\LW}{{^\Theta W}}
\newcommand{\LB}{{^\Theta B}}

\newcommand{\LI}{{^\Theta I}}

\newcommand{\LC}{{^\Theta C}}

\newcommand{\bbI}{{\mathbb I}}

\newcommand{\bbD}{{\mathbb D}}

\newcommand{\ot}{\otimes}

\nc{\Fr}{{\operatorname{Fr} }}

\title[Singular localization in prime characteristic]{Singular localization and intertwining functors for  reductive Lie algebras in prime characteristic}

\author{
Roman Bezrukavnikov
}
\address{\small
Department of Mathematics, Massachusetts Institute of Technology,
77 Massachusetts ave.,
Cambridge, MA 02139, USA
}
\email{ bezrukav@@math.mit.edu }
\author{    Ivan Mirkovi\'c     }
\address{\small
Department of Mathematics and Statistics,
University of Massachusetts,
\   Amherst, MA
01003, USA
}
\email{                mirkovic@@math.umass.edu        }

\author{
Dmitriy Rumynin
}
\address{\small
Mathematics Department, University of Warwick, Coventry,
\
CV4 7AL, England
}
\email{        rumynin@@maths.warwick.ac.uk         }

\thanks{
R.B. was partially supported
by NSF grant DMS-0505466 and  Sloan Foundation,
D.R.  by EPSRC
and
I.M.  by NSF grants.
}

\begin{document}

\begin{abstract}
In \cite {BMR}
 %the paper {\em Localization of modules for a semisimple Lie
%algebra in prime characteristic},
 we observed that, on the level of
derived categories, representations of the Lie algebra of a
semisimple algebraic group over a field of finite characteristic
with a given (generalized) {\em regular} central character can be identified
with coherent sheaves on the formal neighborhood of the
corresponding (generalized) Springer fiber. In the present paper  we
treat  singular  central characters.

The basic step is the Beilinson-Bernstein localization of modules
with a fixed (generalized) central character $\la$ as sheaves on the 
partial flag
variety corresponding to the singularity of $\la$. These sheaves are
modules over a sheaf of algebras which is a version of twisted
crystalline differential operators. %, but is actually larger.
 We
discuss {\em translation functors} and {\em intertwining functors}.
The latter generate an action of the affine braid group on the
derived category of modules with a regular (generalized) central
character, which intertwines different localization functors.
 We also
describe the standard duality on Lie algebra modules in terms of
$\DD$-modules and coherent sheaves.

\end{abstract}

\maketitle

\begin{flushright}
                %\em \bf Dedication: \hspace{0.5in}
To George Lusztig with admiration.
\end{flushright}

%\begin{flushright}{\em Didier Wampas declared he was critical
%of the
% attitude\\ of certain artists (Manu Chao, Noir D\' esir)
%who,\\
% while showing themselves critical of the system, \\
% make a very good living
%out of it.
%\\
%At length he argued the advantages of a swift
%\\
%uncritical death.}
%\end{flushright}

\toc

\setcounter{section}{-1}

\se{\bf Introduction }

This  is a sequel to \cite{BMR}. In the first chapter we extend the
localization construction for modular representations from
\cite{BMR} to arbitrary infinitesimal characters. This is used  in
the second chapter to study the translation functors. We use
translation functors to construct {\em intertwining functors}, which
generate an action of the affine braid group on the derived
categories of regular blocks in modular representation categories.
An application is given in chapter three, where we describe the
duality operation on representations in a regular block in terms of
localization.

More precisely, we consider representations of a
Lie algebra $\fg$ of a reductive  algebraic group $G$ over a field
$\k$ of positive characteristic $p$. (Assumptions on $G$ and $p$ are
in section \ref{Restrictions on G and p}.)
The main result of \cite{BMR} was the
description of the derived category of representations of $\fg$
with a given (generalized) {\em regular} central character, in terms
of coherent sheaves.

The center $\fZ$ of the enveloping algebra $\U=\U(\fg)$ contains the
``Harish-Chandra part'' $\zhc\df\ \U(\fg)^G$ isomorphic to the
Weyl group invariants $S(\fh)^W$ in the symmetric algebra of  the
Cartan algebra $\fh$ of $\fg$. The center  also has the ``Frobenius
part'' $\zfr\cong\OO(\fg^*\tw)$, the functions  on the Frobenius
twist of the dual of the Lie algebra. So, a character of the center
is given by a compatible pair $(\la,\chi)$ of $\la\in\fh^*$ and
$\chi\in \fg^*$. For a regular $\la$, the derived representation
category corresponding to $(\la,\chi)$ is identified in \cite{BMR}
with the derived category of coherent sheaves on the formal
neighborhood of the (generalized) Springer fiber corresponding to
$\chi$. Here we extend this description of the derived category of
representations to singular
 infinitesimal characters $\la$.
For each character of the Harish-Chandra center $\zhc$ there is
$\la$ in the corresponding $W$-orbit in $ \fh^*$ and a partial flag
variety $\PP=G/P$  such that the singularity of $\la$ is of
$\PP$-type, i.e. the stabilizer ${\mathrm{Stab}}_W(\la)$ is the corresponding
parabolic
Weyl group. % (see section \ref{Localization theorem} below for more
%details).
 If we extend the Harish-Chandra central character to a
character $(\la,\chi)$ of the full center, the derived category of
modules is then identified with the derived category of coherent
sheaves on the formal neighborhood of the corresponding {\em
parabolic Springer fiber}, a subvariety of $\PP$ associated to
$\chi$ and $\la$.
In the most interesting case when $\la$ is integral, so $\chi$ is
nilpotent, the parabolic Springer fiber is the set of parabolic
subalgebras of a given type, containing $\chi$.
%; here $\PP=\{Ad(g)P \ |\ g\in G\}$  is the partial flag variety containing $P$ and 
%$\PP_\chi\subset \PP$ is the locus of parabolic subalgebras containing $\chi$.

The basic step is the Beilinson-Bernstein type localization
equivalence of derived categories of modules with (generalized)
infinitesimal character $\la$, and modules over a sheaf of algebras
$\tii\DD^\hla_\PP$ on the partial flag variety $\PP$. Here,
$\tii\DD^\hla_\PP$ is a central completion of the sheaf of algebras
$\tii \DD$; the latter is a version of twisted crystalline
differential operators, but for $\PP\ne \BB$ it is actually larger
than the TDO algebras. For example, in the extreme case  $P=G$,
$\PP=\pt$ we have $\tii \DD=\U(\g)$ is the enveloping algebra
$\U(\g)$. Next, $\tii\DD^\hla_\PP$ turns out to be an Azumaya algebra
that splits on the formal neighborhood of a parabolic Springer fiber,
and this gives an equivalence with coherent sheaves.
We also prove a localization theorem for twisted differential
operators on $\PP$, see
 \ref{DPsec}. The categories of modules with a fixed generalized central
character is then related to the smaller version
of the parabolic Springer fiber, namely, the set of parabolics
of a given type containing the given nilpotent in the radical.

The proofs of these results are, to a large extent, parallel to those
 of \cite{BMR}, though some new ideas are needed.
  The new aspect is that singular localization provides
a better understanding of  {\em translation functors} on
$\fg$-modules.
%   : they are related to push-forward and pull-back of
%   coherent sheaves.
We show that the closely related {\em intertwining functors}
generate an action of the {\em affine braid group} on the derived
category of modules with a fixed (generalized) central character.
This action intertwines different localization functors.

The action and its applications will be further discussed in a
future publication, see also \cite{ICM}.
An application discussed in the present paper
is concerned with a description
of  the usual duality on finite dimensional
Lie algebra modules (with a regular central character) in terms of
coherent sheaves; see remark \ref{onLu} for comments on  motivation
for this result.
%This is directly inspired by Lusztig's definition of an
%involution on the Grothendieck group of coherent sheaves on
% a Springer fiber
% \cite[part II]{Lu} (see remark \ref{onLu}).

%The paper is organized as follows ... .

\sss{} It is an honor for us to dedicate this paper to George
Lusztig. Our intellectual debt to him can not be overestimated. This
paper is a part of the project inspired by his conjectures; but also
much of the other work by the authors is an exploration of beautiful
landscapes whose very existence would 
probably 
remain completely
unknown if not for the hints found in
%hidden in the formulas in
Lusztig's papers.

\sss{}
The main part of this work was accomplished
when R.B. and I.M. were visiting
the Hebrew University in Jerusalem,
and while one of the
authors (I.M.) was visiting
Universite de Cergy-Pontoise
in Paris.
We thank both institutions for the wonderful working
atmosphere. % which this work vainly strives to reflect.

This work is a continuation of \cite{BMR} and is partly a byproduct
of communications acknowledged in {\em loc. cit.} We use the opportunity
to thank these mathematicians again. We are also grateful to 
Iain Gordon, Victor Protsak and Alexander Samokhin for  helpful discussions, and to 
the referees for an extraordinarily thorough job.

\sss*{
Notation
}
Most of the paper considers schemes over
an algebraically closed  field
$\k$
of characteristic $p>0$.
For a closed subscheme $\fY$ of  a scheme $\fX$
the category of modules
on $\fX$ supported on $\fY$
is introduced in
\ref{Derived categories of sheaves supported on a subscheme}.
The inverse image of sheaves  is denoted $f\inv$ and for
$\OO$-modules $f^*$ (both  direct images are denoted $f_*$).
We denote by $\TT_X$
and $\TT^*_X$
the
sheaves of sections of the (co)tangent bundles
$TX$ and $T^*X$.
We
denote by $X/G$ the invariant theory quotient of an affine scheme $X$ by
a group $G$.

\se{\bf
Localization
on partial flag varieties
}

For any partial flag variety $\PP$ we define
a sheaf of algebras
$\tii\DD_\PP$ on $\PP$. When $\PP$ is the full flag variety $\BB$ this is the
Beilinson-Bernstein deformation of differential operators on $\BB$
by twisted differential operators.
When $\PP=\pt$ this is just the enveloping algebra
$\U(\fg)$.
There will be three kinds of geometric objects related to
$\tii\DD_\PP$, besides of the base $\PP$ these are the classical limit
$\tii\fg^*_\PP$ of $\tii\DD_\PP$
and the spectrum
$\ZZ(\tii\DD_\PP)$
of the center
of $\tii\DD_\PP$.
For any sheaf of algebras $\AA$ we denote by $Z(\AA)$ and $\ZZ(\AA)$
the center of $\AA$ and the spectrum of the center.

For arbitrary $\la\in\fh^*$ with
``{\em singularity precisely $\PP$}'',\fttt{
By this we mean that the stabilizer (with respect to the dot action)
of $\la$ in the Weyl group
$W$ of $G$ is precisely the Weyl group $W_\PP$
of Levi groups
of parabolic  subgroups in $\PP$.
}
we establish some equivalences of derived categories.
The final result is that on the level of derived categories
$\fg$-modules with a given generalized central character
(with Harish-Chandra part $\la$)
are the same as coherent sheaves on
$\tii\fg^*_\PP$ supported
on the corresponding parabolic Springer fiber.
An intermediate step is
a Beilinson-Bernstein type localization of
$\U^\la$-modules
to modules for $\DD_\PP^\la\dff
\tii\DD_\PP\ten_{\OO(\fh^*)}\k_\la$.
(Actually we construct some equivalences for more general
$\la$.)

\sus{
Crystalline differential operators
on  flag varieties
}
%Here we recall the setting
%of \cite{BMR}.
We start by presenting settings and notations. 
They are slightly more general than in
 \cite{BMR}, as the group is now reductive and not
necessarily semisimple.

%{\Large Changed ``semisimple'' to reductive}

%\sss{ Semisimple group $G$ } \lab{Semisimple group G} Let $G$ be a
%semisimple simply-connected algebraic group over $\k$. Let $B=T\cdot$
\sss{ Reductive group $G$ } \lab{Semisimple group G} Let $G$ be a
reductive %simply-connected 
algebraic group over $\k$. Let $B=T\cdot
N$ be a Borel subgroup with the unipotent radical  $N$ and a Cartan
subgroup $T$. Let $H$ be the (abstract) Cartan group of $G$ so that
$B$ gives an isomorphism $\io_\fb=(T\con B/N\cong H)$. Let
$\fg,\fb,\ft,\fn,\fh$ be the corresponding Lie algebras. We call
elements of $\Lambda=X^*(H)$ {\em characters} and elements of
$\fh^*$ {\em weights}. The integral weights are differentials of
characters $\fh^*_\Fp\dff d\La\cong \La\ten_\Z\Fp$. The character
lattice $\La$ contains the set of roots $\De$ and  positive roots
$\De^+$ which are chosen as $T$-roots in $\fg/\fb$ via the above
``$\fb$-identification'' $\io_\fb$. Also, $\La$ contains the root
lattice $Q$ generated by $\De$, the dominant cone $\La^+\sub \La$
and the semi-group $Q^+$ generated by $\Delta^+$. Let $I\sub \De^+$
be the set of simple roots. For a root $\al\in\De$ let
$\ch\al\in\ch\De$ be the corresponding coroot. For
$\la,\mu\in\La$ we write $\la> \mu$ if $\la-\mu$ is a sum of
positive roots.
%\in\ \Z_{\geq 0}\cd\De^+$.

Similarly, $\io_\fb$ identifies
 $N_G(T)/T$
with  the Weyl group $W\sub \Aut(H)$. We have  the standard action
of $W$ on $\Lambda$ and on $\h^*=\La\otimes \k$ given by
$w:\lambda\mapsto w(\lambda)$. %=w\cd\la.

 Any character $\nu\in\La$ defines a line bundle
$\OO_{\BB,\nu}=\OO_\nu$ on the flag variety $ \BB\cong G/B $, and a
standard $G$-module $V_\nu\df\ \Hh^0(\BB,\OO_{\nu^+})$ with extremal
weight $\nu$, here $\nu^+$ denotes the dominant $W$-conjugate of
$\nu$ (notice that a dominant weight corresponds to a semi-ample
line bundle in our normalization). We will also write $\OO_\nu$
instead of $\pi^*(\OO_\nu)$ for a scheme $X$ equipped with a map
$\pi:X\to\BB$ (e.g. a subscheme of $\gtil^*$).

\sss{Affine Weyl groups and the dot-action} \label{dotac}
Along with the standard
action we will also use the {\em dot-action}, which differs from the
standard action by the $\rho$ shift: $w:\lambda \mapsto w\bu \lambda
\df\ w(\la+\rho) -\rho$, where $\rho$ is the half sum of positive
roots.

 We will indicate the dot-action by writing
$(W,\bu)$.
% (this is really the action of the $\rho$-conjugate $
%^\rho W$ of the subgroup $W\sub \Wap$).
 The  following rule will be observed in most cases:
    we will use the $\bu$-action
on $\fh^*$  (except in a few auxiliary calculations), while  the
action of $W$ on the Frobenius twist $\fh^*\tw$
%(\ref{})
of $\fh^*$ is the standard action of $W$.
% For instance, $\fh^*/W$
%will always denote the quotient $\fh^*/(W,\bu)$.

 Let  $\Waff\df W\ltimes Q\sub \Wap\df W\ltimes \La $ be
the affine Weyl group
%\fttt{
%This is affine Weyl group of the dual type because we use root lattice. 
%For the affine Weyl group of the same type the coroot lattice is used.}
and the extended affine Weyl group.

These groups act on $\La$ and on $\La_\R=\La\ten_\Z\R$. The group
$\Waff$ is generated by reflections in affine hyperplanes
$H_{\ch\al, n}$ given by
 $\langle
\ch\al,\ \rangle =n$; here $\ch\al$ is a coroot and $n\in \Zet$.
Thus $\Waff$ is the affine Weyl group of the {\em Langlands dual
group} in the standard
 terminology.

We extend the dot-action to $\Waff$ and  $\Wap $ so that $\mu\in
\La$ acts by the $p\mu$-translation: $w:\lambda \mapsto w\bu \lambda
\df\ p \,w(\frac{\la+\rho}{p}) -\rho$.

The hyperplanes $H_{\ch \al}$ define a stratification of $\La_\R$;
by a {\em facet} we mean a stratum of the stratification. Open
facets are called {\em alcoves} and codimension one facets are
called {\em faces}. The set of alcoves is a torsor for
$(\Waff,\bu)$. So, $(\Waff,\bu)$-orbits in the set of faces are
canonically identified with the  faces in the closure of the {\em
fundamental alcove} $A_0$ (the one that contains $(\epsilon-1)\rho$ 
for small positive $\epsilon$).

The group $\Waff$ is a Coxeter group, the Coxeter generators are
reflections in the faces of $A_0$; we denote the latter set by
$\Iaff$. Thus $\Iaff$ is the set of vertices of the affine Dynkin
diagram corresponding to the root system which is dual to the root
system of $G$.

\sss{
Restrictions on the group $G$ and the characteristic $p$
}
\lab{Restrictions on G and p}
$G$ is a connected
reductive algebraic group over a closed field
$\k$ of finite characteristic $p$.
We assume throughout the so-called standard, or
Jantzen-Premet, conditions
on $G$ and the prime $p$ :
\bi
(A)
The derived subgroup
$G'$ is simply connected.
\i
(B)
$p$ is a good prime for $G$,\fttt{
This excludes $G'$ having simple factors of type
$B,C,D,E,F,G$ if $p=2$,
$E,F,G$ if $p=3$ and
$E_8$ for $p=5$.
} and $p$ is odd.\fttt{This condition follows from the rest
unless  all simple factors of $G'$
are of  type $A$. 
It is not a part of standard assumptions, and is imposed  to conform with
the assumptions of \cite{BG}. We expect that it is, in fact,
redundant.}
\i
(C)
There exists a non-degenerate invariant bilinear form on $\g$.\fttt{
For semisimple $G$ satisfying (A) and (B) this amounts to $G$ not 
having simple factors $SL(mp)$, though $G=GL(n)$ satisfies
the assumption for any $n$.
}
\ei

If $p$ satisfies these conditions for a reductive group $G$,
then it also does so for a Levi subgroup $L\subset G$, see 
\cite[6.5]{Ja} or
\cite[Lemma 3.2]{BG}.\fttt{Here and below we use \cite{BG}
as a convenient reference for some basic statements about
Lie algebras. 
 We have not attempted
to identify the original source of the arguments; most of them
go back to Springer, Steinberg or Bardsley and Richardson.}
For a semisimple group $G$ these conditions amount to $p$ being very
good for $G$, some claims  in \cite{BMR} have been proved under this 
assumption.
For a reductive group $G$ they ensure
the standard structural results that we use such as the following

\slem
\label{connected centraliser}
\cite[Lemma 3.2]{BG}, \cite[7.4]{Ja}
%\fttt{In {\em loc. cit.} the assumption 
%$p\ne 2$ is also imposed. If $p=2$ is good then all simple
%factors of $G$ are of type $A$, then the statement of the Lemma
%is obvious.}
{\em
The centralizer in the group $G$
of any semisimple $h\in\fg$  is a Levi subgroup.
}

Starting from  chapter 2 we use  regular integral weights. So, in
order for such weights to exist we need the stronger assumption $p\geq
h$ where $h$ is the maximum of the Coxeter numbers of simple factors of $G^\prime$. For a simple factor $h= \lb \rho,\ch\al_0\rb+1$ where
$\ch\al_0$ is the highest coroot. In \cite{BMR} we mostly worked under 
the  assumption $p> h$.

\sss{ Crystalline differential operators } \lab{Crystalline
differential operators} The ring of {\em crystalline differential
operators} $\DD_Y$ on a smooth variety $Y$ is the sheaf of algebras
generated by functions $\OO_Y$, vector fields $\TT_Y$ and the
standard relations between them. These include
the module, commutator and Lie algebroid relations
$$
f\cd \del= f\del, \ \del\cd f-f\cd \del= \del(f),\ 
\del\cd \del'-\del'\cd \del=
[\del,\del'],\ \del,\del'\in\TT_Y,\ f\in\OO_Y
.$$
Locally $\TT_Y$ has a frame
$\del_i,\ 1\le i\le n$, and then $\DD_Y$ has a frame
$\del^I=\del_1^{I_1}\cddd  \del_n^{I_n},\ I\in\Z_+^n$. 
%According to
%\cite{BMR} (see also \cite{OV} for more information) 
The sheaf
$\DD_Y$ is an Azumaya algebra\fttt{
An Azumaya algebra is a sheaf of locally matrix algebras
in etale topology \cite{MI}.
}
on the Frobenius twist $T^*Y\tw$ of
the cotangent bundle \cite[Theorem 2.2.3]{BMR}\fttt{
See also \cite{OV} for more information.
}. 
Here, a vector field $\xi$ on $Y$ has $p\thh$
power $\xi^p$ in $\DD_Y$ and also another kind of $p\thh$ power
$\xi^{[p]}$, calculated this time in operators on $\OO_Y$, which
happens to be a vector field again. Now $\TT_Y\tw\sub\OO_{T^*Y\tw}$
maps to differential operators by $\TT_Y\tw\ni\xi\tw\mm\
\bio(\xi\tw)\dff \xi^p-\xi^{[p]}\in\DD_Y$, the difference of two
$p\thh$ powers. This gives an isomorphism $\bio:\OO_{T^*Y\tw}\to
Z(\DD_Y)$. When $Y$ is a torus $T$ then $\bio$ in particular maps
$S\ft\tw$ into $S\ft$, the corresponding map of spaces is the
Artin-Schreier map which we denote $\AS=\AS_T:\ft^*\to \ft^*\tw$. If
$\la$ is the differential of a character of $T$ and $c\in \k$, then
$\AS(c\la)=(c^p-c)\la\tw$; thus one can identify $\fh^*$ with the
affine space $\bbA^n$ so that $\AS$ is identified with the product of
$n$ copies of the standard Artin-Schreier map.

\sss{
The
Harish-Chandra and Frobenius ingredients
of
the center of $\U(\fg)$
}
\lab{center}
Let now $\U\dff\U(\fg)$ be the enveloping algebra of $\fg$.
The subalgebra
of $G$-invariants
$\zhc\df\ \U(\fg)^G$
is  central in $\U(\fg)$.
On the other hand
as the algebra of left invariant differential operators,
$\U=\DD(G)^{G\tim 1}$
has a central subalgebra
$\OO(T^*G\tw)^{G\tim 1}
=\OO(\fg^*\tw)$ which we call the
Frobenius center
$\zfr$.

\stheo
\cite{Ve,KW,MR}
{\em
(a)
There is a canonical isomorphism
$
\U^G@>{i_{{\mathrm{HC}}}}>> S(\fh)^{(W,\bu)}
$.\fttt{
A geometric description of $i_{\mathrm{HC}}$ is sketched in
\ref{Algebras}.
}

(b)
The center
$\fZ$
of  $\U$
is a combination of the Harish-Chandra part
and the Frobenius part
$$
\fZ\conn\
\zfr
\ten_{\zfr\cap\zhc}
\zhc,
\
i.e.,
\
\ZZ(\U)=\fg^*\tw\tim_{\fh^*\tw/W}\fh^*/W
.$$
Here,
$\fg^*\tw
@>>>
\fh^*\tw/W$
is the adjoint quotient, while the map
$\fh^*/W\
@>>>
\fh^*\tw/W
$
comes from the Artin-Schreier map
$\fh^*
@>{\AS}>>
\fh^*\tw
$, cf.
%$S(\fh\tw)@>{\AS}>>S\fh$ defined in
section \ref{Crystalline differential operators}.
}

For $\la\in\fh^*$ and $\chi\in\fg^*\tw$, we introduce the central
specializations of the universal enveloping algebra: $\U^\la\dff \U\ten_\zhc \k_\la$, $u_\chi\dff
\U\ten_\zfr \k_{\chi\tw}$, and $\U^\la_\chi\dff \U\ten_\fZ
\k_{(\chi,\la)}$.
These constructions can be applied to a parabolic subalgebra $\fp$ or its Levi factor $\bar\fp$ by restricting the reduction data to the corresponding subalgebras.
For instance, 
$\U^\la_\chi(\bar\fp) = \U(\bar\fp)\ten_{\fZ(\bar\fp)}\k_{(\chi|_{\bar\fp},\la|_{\bar\fp})}$
or $u_\chi(\fp)$ is the reduced enveloping algebra.
We also introduce the extended universal enveloping algebra
$\tii \U=\U\otimes _{\zhc} S(\fh)$, denoted
 $\tii \U (\bar\fp)$ when applied to a Levi factor.

%Do we use these?

\sss{ Unramified weights } \lab{Unramified weights} We will say that
$\la\in \fh^*$ is {\em unramified} if the map
${\AS}/W:\fh^*/W@>>>\fh^*\tw/W$ is unramified at $W\bu\la$. This is
equivalent to the equality of the stabilizers $W_\la\sub W_{{\AS}(\la)}$.
Notice that the stabilizers are Coxeter subgroups of $W$, in fact, 
they are the Weyl groups of the Levi factors
 in Lemma~\ref{connected centraliser}.
It is easy to observe that ${\AS}(w\la)={\AS}(w\bu\la)=w{\AS}(\la)$.
So, a weight $\la$ is {\em ramified} if and only if for some root $\al$ one has
$\langle \ch\alpha, \la+\rho\rangle\in \Fp^*=\Fp-\{0\}$\ie $\la$ is
integral and regular for $\al$. It is well known that such
ramification  produces reducibility of standard $\fg$-modules. In
the unramified case the representation theory simplifies, because
the corresponding central reductions (or completions) of the
enveloping algebra are Azumaya \cite[Theorems 2.6 and 3.10]{BG}, \cite{BrGo}.

 Let $\fh^*\unr \subset \fh^*$ be
the open set of all unramified weights.

\sss{
Azumaya locus
}
\lab{Azumaya}
The central variety $\ZZ(\U)$ contains an open part
$$
\ZZ(\U)\unr\dff
\ZZ(\U)\tim_{\fh^*/W}\ \fh^*\unr/W
=
\fg^*\tw\tim_{\fh^*\tw/W}\ \fh^*\unr/W
\ \sub\ \fg^*\tw\tim_{\fh^*\tw/W}\ \fh^*/W
=\
\ZZ(\U)
.$$
$\U$ is known to be generically Azumaya
and the Azumaya locus
in $\ZZ(\U)$ coincides
with the smooth locus
(see
Theorems 2.5, 2.6
in
\cite{BG}, or \cite{BrGo} for the semi-simple case). The variety
$\ZZ(\U)\unr$
is smooth, because it admits an etale map to $\fg^*\tw$;
thus  $\U$ is Azumaya
over $\ZZ(\U)\unr$.

\sus{
The algebras
$\tii\DD_\PP$
on partial flag varieties $\PP$
}
\label{redu}
In this section we define and discuss various sheaves of rings related
 to $\fg$; further sections will be devoted to relations between the
(derived) categories of sheaves of modules.
In \ref{Torsors} we start with a
general construction
that in particular
produces
algebras
$\tii\DD_\PP$.
A reader unfamiliar with torsors should consult \cite{Sk}.

\sss{
The algebras $\tii\DD_X$ associated to torsors
}
\lab{Torsors}
A torsor  $\tii X@>\pi>> X$ for a reductive group
$L$
defines a Lie algebroid
$\tii\TT_X\df\ \pi_*(\TT_{\tii X})^L$
with the enveloping algebra
$\tii\DD_X\df\ \pi_*(\DD_{\tii X})^L
$.
Let $\fl$ be the Lie algebra of $L$.
Locally, any trivialization of the torsor
splits the exact sequence
$0\ra \fl\ten\OO_X\ra\ \tii\TT_X\ra\TT_X\ra\ 0$
and gives
$\tii\DD_X\ \cong \DD_X\ten \ \U(\fl)$.
So the map
$\U(\fl)@>>>\tii\DD_X$
given by the $L$-action is an embedding,
$\tii\DD_X$ is a locally free $\U(\fl)$-module,
and the center $Z[\U(\fl)]$ of $\U(\fl)$
is central in $\tii\DD_X$. Recall that
$\ZZ(\U(\fl))
=
\fl^*\tw\tim_{\fh^*\tw/W_L}\fh^*/W_L
$
for the Cartan group
$H$ and the Weyl group $W_L$ of $L$
(Theorem
\ref{center}).
Since
$\U(\fl)$ is free over its Harish-Chandra center
$\zhc(\fl)=\U(\fl)^L\con \OO(\fh^*/W_L)$,
it follows that
$\tii\DD_X$ is  locally free over $\zhc(\fl)$.
We consider the specializations
$\DD_X^\la\df\ \tii\DD_X\ten_{\zhc}\k_\la$. Notice that  the
augmentation $\U(\fl)@>>>\k$ gives $\tii\DD_X \otimes_{\U(\fl)} \k
\con\DD_X$.

The center $\OO_{T^*\tii X\tw}$
of $\DD_{\tii X}$ gives  another central subalgebra
$(\pi_*\OO_{T^*{\tii X}\tw})^L
=\
\OO_{\tii T^* X\tw}
$ of $\tii\DD_X$.
We combine two central subalgebras
into a map from
functions on
$
\tii T^*X\tw\tim_{\fl^*\tw}\ZZ(\U(\fl))
=\
\tii T^*X\tw\tim_{\fh^*\tw/W_L}\fh^*/W_L
$
to the center of $\tii\DD_X$
(we use
the exact sequence
$0\ra
\fl\ten\OO_X
\ra\ \tii\TT_X\ra
\TT_X
\ra 0$).
The above local trivializations
now  show that this is an isomorphism\ie
$\ZZ(\tii\DD_X)=
\tii T^*X\tw\tim_{\fh^*\tw/W_L}\fh^*/W_L
$. Moreover we also see  that
the Azumaya locus of  $\tii\DD_X$ in $
\ZZ(\tii\DD_X)$
is the fibered product over ${\fl^*\tw}$
of
$
\tii T^*X\tw
$
and the Azumaya
locus
of $\U(\fl)$.
This implies that  the Azumaya locus in
$
\ZZ(\tii\DD_X)$
is precisely the
smooth part of $
\ZZ(\tii\DD_X)$
(use \ref{Azumaya}).
For any $\la\in\fh^*$
the center of the specialization
$\DD^\la_X$ contains the algebra of functions on
$\ZZ(\tii\DD_X)\tim_{\fh^*/W_L}\la=\
\tii T^*X\tw\tim_{\fh^*\tw/W_L}\AS(\la)
$.
This is a torsor for
$T^*X\tw$ over the quotient
$
X\tw\tim (\fl^*\tw\tim_{\fl^*\tw/W_L}\ {\AS(\la)})
$.
In particular, by using
\ref{Azumaya} for the group $L$ we see that
for any $L$-unramified weight $\la\in\fh^*\unr(L)$, the sheaf
$\DD_X^\la$ is an Azumaya algebra
on
$
\tii\ZZ\tim_{\fh^*/W_L}\la
$.

For instance if $\la=d(\phi)$ is the differential of a character
$\phi$ of   $L$ then  $\AS(\la)=0$ and $\DD_X^\la$ is identified with
the sheaf $^{\OO_\phi}\DD_{X}\dff \OO_\phi\ten\DD_X\ten\OO_\phi\inv$
of differential operators on sections of the line bundle $\OO_\phi$
on $X$, associated to the $L$-torsor $\tii X$ and  a character
$\phi$ of $L$. The following fact (which will not be used in the
sequel)
is proved in the same way
  as the version
in \cite[Lemma 2.3.1]{BMR}.

\slem
{\em
Let $\nu\in \fl ^*$ be the differential of some character of $L$.
Define a morphism $\tau_\nu$ from $\tii
T^*X\tw\times_{\ft^*\tw}\ft^*$ to itself by
$\tau_\nu(x,\la)=(x,\la+\nu)$. The Azumaya
loci of
algebras
$\tii\DD_X$ and $\tau_\nu^*(\tii \DD_X)$ are the same,
and the corresponding restrictions are equivalent Azumaya algebras.
More precisely, the choice of $\eta$ with 
$\nu=d\eta$ determines an equivalence between the two Azumaya algebras.
 The corresponding equivalence of the categories
of modules is given by $\MM\mapsto \OO_{\PP,\eta}\otimes _{\OO_\PP}\MM$.}

\srem The equivalence of Azumaya algebras
can also be deduced from Lemma \ref{lemsplit} below. 

\sss{
Algebras $\tii\DD_\PP$
}
\lab{Algebras}
Let $\PP$ be a partial flag variety and
$\BB
@>a>>
\PP
@>{b}>>
\pt$.
Here $P$ will always denote a
parabolic subgroup $P\in\PP$. Denote by
$J$
its unipotent radical
and let $\barr {P}\df\ P/J$ be the Levi quotient of $P$.
Its Lie algebra $\barr\fp$ has Cartan algebra $\fh$
and the Weyl group $W_{\PP}\sub W$.

Over $\PP=G/P$ there is a $\barr P$-torsor
$\pi=\pi_\PP:\tii\PP=G/J
@>>>
\PP$.
As in
\ref{Torsors}
it yields  an algebra
$\tii\DD_\PP\dff
(\pi_*\DD_{\tii\PP})^{\barr P}
$
and
its specializations
${\DD^\la}_\PP
\dff
\tii\DD_\PP \ten_{ S(\fh)^{W_\PP} }\ \k_\la
$
for
$\la\in\fh^*$.
The action of  $G\tim \barr P$  on  $\tii\PP=G/J$ by
$(g,p)\cd sJ
\df\ gsp\inv J
$, differentiates to
a map $\fg\pl\barr\fp@>>>\tii\TT_\PP$ which extends to
$\U(\fg)\ten \U(\barr \fp)
@>>>\tii\DD_\PP$.
When $\PP$ is  the full flag variety
$\BB$,  then $\tii\DD_\BB$
is defined by the  $H$-torsor
$\tii\BB\df\ G/N @>\pi>> \BB$
and it
is a deformation over $\fh^*$ of
the ring of differential operators $\DD_\BB=\DD_\BB^0$.\fttt{
We can use $\tii\DD_\BB$ to
describe the Harish-Chandra map
(\ref{center}).
The
$H$-action on $\tii\BB$
gives an isomorphism
$
\U(\fh)\con \Ga(\BB,\tii\DD_\BB)^G
$, and then
the $G$-action gives
the map
$\U^G
@>>>
\Ga(\BB,\tii\DD_\BB)^G\cong S(\fh)
$
which gives  an  isomorphism
$
\U^G@>{i_{\mathrm{HC}}}>> S(\fh)^{(W,\bu)}
$.
}
When
$P=G$ then  $\PP=\pt$ and $\tii \DD_\PP=\U$.

The sheaf of algebras 
$\tii \DD_\PP$ carries a filtration induced by the canonical 
filtration  on crystalline differential operators, cf. \cite[1.2]{BMR}.
The associated graded is canonically identified with the sheaf of regular
functions on the $T^*\PP$-torsor over $G\times_P  \barr\fp^*$,
%On the ``classical level'' 
%the $\barr P$-torsor $\tii\PP$ defines a
which is defined by: 
$\gtil^*_\PP
\dff \tii T^*\PP =\ \{(\fp,\chi)\in \PP\tim \fg^*,\ \chi|_{{\mathrm{nil}}(\fp)}=0\}$.
For $\PP=\BB$ we use a simplified notation $\gtil^* = \gtil^*_\BB$. 
We have projections 
$\PR_\fg:\gtil^*_\PP\to \fg^*$,
$\PR_\fg(\fp,\chi)=\chi$ and $\PR_{\barr \fp}: \gtil^*_\PP \to \barr
\fp^*/\barr P = \fh^*/W_\PP$ which sends $(\fp,\chi)$ to the
coadjoint image of $\chi|_\fp\in (\fp/{\mathrm{nil}}(\fp))^*=\barr \fp^*$; they
yield a map $\PR=\PR_\fg\times\PR_{\barr\fp}: \gtil^*_\PP\to
\fg^*\times_{\fh^*/W} \fh^*/W_\PP$ (we use the canonical
isomorphism $\fg^*/G @>>> \fh^*/W$).

\sss{
Cohomology of $\tii \DD_\PP$ and Azumaya property
}
\lab{D-cohomology and center} Define the open subset
$\fh^*\PPunr\sub\fh^* $ of {\em $\PP$-unramified} weights as
$\fh^*\PPunr \dff \fh^*\unr(\barr P) $ \ie weights that are
unramified for the  Levi group $\barr P$ of $P$. By definition in 
\ref{Unramified weights}, this means that  the map
${\fh^*/W_\PP}@>>>\fh^*/W$ is unramified at $W_\PP\bu\la$. The
condition is $(W_\PP)_\la=(W_\PP)_{{\AS}(\la)}$, or equivalently $\lb
\ch\alpha, \la+\rho\rb\not\in \Fp^*$ for roots $\al$ in $\De_\PP$.
It is clear that $\fh^*\PPunr$  contains the set  of
weights with singularity of type $\PP$\ie
with the stabilizer equal to $W_\PP$.
Let us also define the unramified part
of the central variety of $\tii\DD_\PP$
as the open subvariety
$$
\ZZ(\tii\DD_\PP)\unr
\dff
\ZZ(\tii\DD_{\PP})
\tim_{\fh^*/W_\PP}
\ (\fh^*\PPunr /W_\PP)
\sub\
\ZZ(\tii\DD_\PP)
.$$
In particular,
$
\ZZ(\U)\unr
\dff
\fg^*\tw
\tim_{\fh^*\tw/W}
\ \fh^*\unr
\sub\
\ZZ(\U)
$.

\spro
{\em
(a)
$\Rr a_*\tii\DD_\BB\cong \tii\DD_\PP\ten_{S(\fh)^{W_\PP}}
S(\fh)
$ canonically,
hence
$\tii\DD_\PP\cong\
\tii\DD_\BB^{ W_\PP }
$.
For each $\la\in\fh^*$,\
$\Rr a_*\DD_\BB^\la\cong \DD_\PP^\la
$.

(b)
The map
$\U(\fg)\ten S(\fh)^{W_\PP}@>>>\Ga(\PP,\tii\DD_\PP)$
factors through
$
\tii \U^{W_\PP}
=\ \U\ten_{S(\fh)^{W}} S(\fh)^{W_\PP}
$
and yields an isomorphism
$\tii \U^{W_\PP}
\con\RGa(\tii\DD_\PP)$.
Also $\U^\la\con\RGa(\DD^\la_\PP)
$ for $\la\in\fh^*$.

(c)
The
center of $\tii\DD_\PP$ is
the
sheaf
of functions on
$
\ZZ(\tii\DD_\PP)
=
\tii\fg^*_\PP\tw\tim_{\fh^*\tw/W_\PP}\ \fh^*/W_\PP$.
The algebra
$\tii\DD_\PP$ is
Azumaya over
the open subvariety
$
\ZZ(\tii\DD_\PP)\unr
\sub
\ZZ(\tii\DD_\PP)
$.

(d) Consider the canonical map 
$\PP=G/P \to \QQ=G/Q$, $P\sub Q$
between two partial flag varieties. It induces a  map of central
varieties 
$\vpi^\PP_\QQ:\ZZ(\tii\DD_\PP)@>>>\ZZ(\tii\DD_\QQ) 
$\ie 
$\vpi^\PP_\QQ:
\tii\fg^*_\PP\tw\tim_{\fh^*\tw/W_\PP}\ \fh^*/W_\PP @>>>\
\tii\fg^*_\QQ\tw\tim_{\fh^*\tw/W_\QQ}\ \fh^*/W_\QQ;
%\fg^*\tw\tim_{\fh^*\tw/W}\ \fh^*/W
$
in particular, for $\QQ=\pt$ we get
$$
\vpi^\PP\dff \vpi^\PP_{\pt}:
\ZZ(\tii\DD_\PP)=\tii\fg^*_\PP\tw\tim_{\fh^*\tw/W_\PP}\ \fh^*/W_\PP
@>>>\ \ZZ(\U)= \fg^*\tw\tim_{\fh^*\tw/W}\ \fh^*/W .$$ Then the
preimage of the Azumaya locus is contained in the Azumaya locus,
i.e.,
$$(\vpi^\PP_\QQ)^{-1}(\ZZ(\tii\DD_\QQ)\AZ)
\subset\ZZ(\tii\DD_\PP)\AZ;$$ and moreover the two Azumaya
algebras are canonically related via the pull-back :
\begin{equation}\label{where_to}
\tii \DD_\PP|_{(\vpi^\PP_\QQ)^{-1}(\ZZ(\tii\DD_\QQ)\AZ)}
\cong (\vpi^\PP_\QQ)^*(\tii \DD_\QQ)
.
\end{equation}

(e) The complement to 
$(\vpi^\PP)^{-1}(\ZZ(\U)\AZ)$
in $\ZZ(\tii \DD_\PP)$ has codimension at least two
for every~$\PP$.
}

\pf
The first claim in (a) implies the rest of (a),
because the derived global sections functor commutes
 with the derived tensor product over
a ring of global endomorphisms, while the ring $S(\h)$ acts locally freely
both on $\tii \DD_\BB$ and on $\tii \DD_\PP\otimes _{S(\h)^{W_\PP}}
S(\h)$ by \ref{Torsors}.

For a semisimple group $G$ and $\PP=\pt$ the first claim in (a)
has been established in \cite[Proposition 3.4.1]{BMR}. A very similar
argument shows 
that the first claim 
holds for $\PP=\pt$ and a reductive group $G$.
We proceed to deduce the general case.
For a fixed $P\in\PP$ choose a Levi factor
$L$ of $P$ and let $P_-=LJ_-$ be the corresponding
opposite parabolic.
The  maps
$\tii\BB
@>>> %\pi_\BB>>
\BB
@>a>>
\PP
@<<< %\pi_\PP<<
\tii\PP$,
restricted to the neighborhood
$\UU=
J_-\cd\fp$  of $\fp\in\PP$,
are naturally identified with
$
\UU
\tim
\tii\BB(\barr P)
@>>> %\si_\BB>>
\UU\tim \BB(\barr P)
@>>> %{id\times a_{\barr P}>>
\UU\tim \pt
@<<< %\si_\PP<<
\UU\tim \barr P$,
where
$\BB(\barr P)$ is flag variety of $\barr P$.

Therefore
$\DD_{\tii\PP}|\vpi^\PP\inv\UU
\cong
\DD_\UU\btim \DD_{\barr P}$
and
$\tii\DD_{\PP}|\UU
\cong
\DD_\UU\btim \DD(\barr P)^{\barr P\tim 1}
=\
\DD_\UU\btim \U(\barr\fp)
$.
Now, equality
$\DD_{\tii\BB}|a\inv\UU
\cong
\DD_\UU\btim \tii\DD_{\BB(\barr P)}$
shows that on $\UU$:
$$
\Rr a_*\tii\DD_\BB
\cong
\DD_\UU\ten_\k
\RGa[\BB(\barr P),\tii\DD_{\BB(\barr P)}]
=\
\DD_\UU\ten_\k\
\tii \U(\barr \fp)
=\
[\DD_\UU\ten_\k\
\U(\barr \fp)]\ten_{S(\fh)^{W_\PP}}S(\fh)
\cong
\tii\DD_\PP\ten_{S(\fh)^{W_\PP}}S(\fh)
.$$

Similarly,
in (b) the
second sentence
follows from the first one.
The case $\PP=\BB$ follows from \cite[Proposition 3.4.1]{BMR},
where it is stated for semisimple $G$. Thus we have:
$\RGa(\PP;\Rr a_*\tii\DD_\BB)
=
\Rr(ba)_*\tii\DD_\BB
\cong
\U\ten_{S(\fh)^W}S(\fh)
$, where $b$ denotes the projection $\PP\to \pt$.
Substituting the description of $\Rr(ba)_*\tii\DD_\BB$ from (a) we see
that the morphism $\U\otimes _{S(\h)^W}S(\h)^{W_\PP}\to \Rr b_*\tii\DD_\PP$
becomes an isomorphism after tensoring with $S(\h)$ over
$S(\h)^{W_\PP}$. Since $S(\h)$ is free over $S(\h)^{W_\PP}$ by \cite{De},
this implies the first sentence in (b).

%Assumptions on  characteristic $p$ 
%imply that the functor of $W_\PP$-invariants
%is exact, so
%$
%\RGa(\tii\DD_\PP)
%=
%\RGa([a_*\tii\DD_\BB]^{W_\PP})
%\cong
%[\U\ten_{S(\fh)^W}S(\fh)]^{W_\PP}
%$.

(c) has already been observed in \ref{Torsors}.

Claim (d) in the particular 
 case $\PP=\BB$, $\QQ=\pt$ and $G$ semisimple is in 
\cite[Proposition 5.2.1]{BMR}.
The general case reduces to this one as in the proof of (b).

(e) follows from the fact that $(\chi,\la)\in \ZZ(\U)$ lies in $\ZZ(\U)\AZ$
provided that $\chi$ is regular (not necessarily semisimple) \cite[Theorems 2.5, 2.6, Corollary 3.4]{BG},
and that the preimage of regular elements in $\gt^*$ has complement of
codimension two.
The latter claim is equivalent to non-regular
elements forming a subset of codimension two in
a fixed Borel subalgebra $\fb\subset \fg$. Writing $\fb=\fh \oplus 
\fn$ we see that the set of non-regular elements
is contained in the union of $\fh_\alpha \times \fn$, where $\fh_\alpha$
runs over the set of root hyperplanes. Moreover, $\fh_\alpha \times \fn$
is not contained in the set of non-regular elements for any root hyperplane
$\fh_\alpha$, which implies the bound on the codimension.
\epf

\rems
(0)
$
\ZZ(\tii\DD_\PP)$ consists of all
$
({\fp}\tw,\chi\tw,W_\PP\bu\la)\in
\PP\tw\tim\fg^*\tw\tim\fh^*/W_\PP
$,
such that
(i)
$
\chi\perp\fu
$
(so that
the restriction
$\chi|\fp$
factors
to
$\chi_{\fp}\in\ \barr{\fp}^*$),
and\
(ii)
the image of
$
\chi_{\fp}\tw
$
in
$
(\barr{\fp}^*/\barr {P})\tw
\cong
\fh^*\tw/W_\PP
$ is the orbit
$W_\PP \nu $ for $\nu=\AS(\la)
$.
Here,
(i) means that
$\fp$ lies in the (generalized)
Springer fiber
$\PP_\chi\dff
\tii\fg^*_\PP\tim_{\fg^*}\chi$.

(1) The fibers of $\tii\DD_\PP$ and $\DD^\la_\PP$ at $P\in\PP$ are
the induced
$\fg$-modules $\Ind_{\U(\fp)}^{\U(\fg)}\ \tii\U(\barr\fp)$ and
$\Ind_{\U^\la(\fp)}^{\U^\la(\fg)}\ \U^\la(\barr\fp)$ (see the proof of
(a) above). The fiber of ${\tii\DD_\PP}$ at a point
$(\fp\tw,\chi\tw,W_\PP\bu\la)$ of $\ZZ(\tii\DD_\PP)$ is $
\Ind^{u_\chi(\fg)}_{u_\chi(\fp)} \ \U^\la_\chi(\bar\fp) $.

(2) $\ZZ(\tii\DD_\PP)$  is a torsor for the twisted cotangent bundle
$T^*\PP\tw=\ (G\tim_{P}\ \fp^\perp)\tw$\ over the quotient $
(G\tim_{P}\ \barr\fp^*) \tw \tim_{\fh^*\tw/W_\PP} \ \fh^*/W_\PP =\
G\tim_{P}\ [\barr\fp^*\tw \tim_{\fh^*\tw/W_\PP} \ \fh^*/W_\PP] =\
G\tim_{P}\ \ZZ(\U(\barr\fp)) $.

(3) The classical limit of the claims (a) and (b) above is the
observation that the canonical maps
$\tii\fg^*_\BB@>>>\tii\fg^*_\PP\tim_{\fh^*/W_\PP}\fh^*$ and
$\tii\fg^*_\PP@>>>\fg^*\tim_{\fh^*/W_\PP}\fh^*/W_\PP$ are
affinizations over $\PP$ and over the point respectively. On the
level of centers we get that $ \ZZ(\tii\DD_\PP) @>>> \ZZ(\tii
\U)/W_\PP=\ \fg^*\tw\tim_{\fh^*\tw/W}\fh^*/W_\PP $ is an
affinization.
%%% (All three  maps are proper,
%isomorphisms off codimension $3$ in the target, and targets are  normal, so
%%%    it suffices to notice that the targets are indeed affine over
%%%  $\PP$ and  a point.)

\sus{
Splitting on parabolic Springer fibers
}
\lab{Splitting on Springer fibers}

Here we generalize the observation from \cite{BMR} that the {\em
equivalence class} of the Azumaya algebra on the (twisted) cotangent
bundle to $\BB$ is the pull-back of a certain class under the Springer map.
This implies, in particular,
 that when
$\DD_\PP^{\la} $
is an Azumaya algebra\ie
$\la\in\fh^*\PPunr$,
then this Azumaya algebra splits
on formal neighborhoods of parabolic Springer fibers
$
\PP_{\chi,\AS(\la)}\tw $ in $ \ZZ(\tii\DD_\PP) $.

\sss{
An auxiliary etale covering
of $\fh^*$
}
\lab{An auxiliary etale covering}
%Lemma \ref{lemsplit} is the key result of this subsection;
%to state it
We need some notations. Consider the addition map $\fh^*\unr\times
\fh^*_\Fp=\cupl_{\la_\Fp\in \fh^*_\Fp} \fh^*\unr\to \fh^*$,
$(\la\unr, \la_\Fp)\mapsto \la\unr +\la_{\Fp}$ and the induced
maps on the quotients: $a_\PP: (\fh^*\unr\times \fh^*_\Fp)/W_\PP
\to \fh^*/W_\PP$. 
Notice that the action of $W$ on $\fh^*\unr$ and  $\fh^*$ is the dot-action
but the action on $\fh^*_\Fp$ is the usual action.
When
%It will be convenient to allow the case
$P=G$\ie  $\PP=\pt$ then $W_{\pt}=W$ and $a\dff
a_{\pt}:(\fh^*\unr\times \fh^*_\Fp)/W\to \fh^*/W$.

\slem
{\em a) For every partial flag variety
$\PP$ the map $a_\PP$ is an etale covering.
%for every $\PP=G/P$ (including $\PP=\pt$).

b) Consider the two maps $\prr_1, \, a_\PP: (\fh^*\unr\times
\fh^*_\Fp)/W_\PP \to \fh^*/W_\PP$. Their compositions with the
Artin-Schreier map $\AS/W_\PP$ are equal.}

\pf a) To check that the map is etale it suffices to check the equality
of stabilizers  ${\mathrm{Stab}}_W(\la\unr, \la_\Fp) =
{\mathrm{Stab}}_W(\la\unr+\la_\Fp)$  for every $(\la\unr, \la_\Fp) \in
\fh^*\unr \times \fh^*_{\Fp}$.  In view of Lemma
\ref{Restrictions on G and p}, the
stabilizer of $\la\in \fh^*$ is generated by reflections $s_\al$
where $\al$ runs over coroots  such that $\langle \al,
\la+\rho\rangle=0$. By the definition of an unramified element of
$\fh^*$ for every coroot $\al$ we either have  $\langle \al,
\la+\rho\rangle\not \in \Fp$, in which case $\langle \al,
\la\unr+\la_{\Fp}+\rho\rangle\ne 0$; or $\langle \al,
\la\unr+\rho\rangle=0$. Thus ${\mathrm{Stab}}_W(\la\unr)\supseteq
{\mathrm{Stab}}_W(\la\unr+\la_\Fp)$, hence
 ${\mathrm{Stab}}_W(\la\unr, \la_\Fp) =
{\mathrm{Stab}}_W(\la\unr+\la_\Fp)$.

It remains to show that  $a_\PP$ is a covering, thus we need to check
that  $\fh^*=\fh^*\unr+\fh^*_\Fp$. Pick $\la\in \fh^*$. We claim
that there exists $\la_\Fp\in \fh^*_\Fp$ such that $\langle
\al,\la+\rho\rangle=\langle \al,\la_\Fp\rangle$ for every coroot $\al$
such that $\langle \al,\la\rangle \in \Fp$; then $\la\unr \dff
\la-\la_\Fp$ is unramified. To check the existence of  $\la_\Fp$
notice that the conditions on $\la_\Fp$ constitute a finite
collection of (non-homogeneous) linear equations on an element of
$\fh^*_\Fp$; the equations are defined over $\Fp$ and have a
solution $\la$ in $\fh^*$, i.e. they have a solution over a larger
field of coefficients. By standard linear algebra they also have a
solution over $\Fp$.

b) is clear,  because
$\AS(\la+\la_\Fp)=\AS(\la)$ for $\la_\Fp\in \fh^*_\Fp$. \epf

\medskip

\sss{
Azumaya algebras on an etale covering
of central varieties
}
We now use $a_\PP$ to get coverings of the central varieties
described in Proposition \ref{D-cohomology and center}.c
 and Theorem \ref{center}.
We set
$$
\tii\ZZ(\tii \DD_\PP)
\dff
\ZZ(\tii\DD_\PP)\times_{\fh^*/W_\PP}
(\fh^*\unr\times \fh^*_\Fp)/W_\PP
=\gt^*_\PP\tw\times_{\fh^*\tw/W_\PP}
 \left(( \fh^*\unr\times \fh^*_{\Fp})/W_\PP \right)
.$$
We also let $\tii \ZZ(\tii \DD_\PP)\AZ=(\id\times
a_\PP)^{-1}(\ZZ(\tii \DD_\PP)\AZ)$ be the preimage of the Azumaya
locus of $\tii \DD_\PP$.
%We define its  unramified part
%$\tii \ZZ(\tii \DD_\PP)\unr$
%as the inverse of $\fh^*\PPunr/W_\PP$
%under the map
%$\tii \ZZ(\tii \DD_\PP)\to(\fh^*\unr\times \fh^*_\Fp)/W\to \fh^*/W_\PP$
%where the first map is the projection and the second
%map is constructed in
%Lemma
%\ref{An auxiliary etale covering}.b.
%So,
%$$
%\tii \ZZ(\tii \DD_\PP)\unr
%=\
%( \id\times a_\PP)^{-1}
%\ZZ(\tii \DD_\PP) \to {\fh^*/W_\PP})\inv \ \fh^*\PPunr /W_\PP =\ (
%\id\times a_\PP)^{-1} \tii \ZZ(\tii \DD_\PP)\unr .$$

 Recall that
$\vpi^\PP\dff \vpi^\PP_{\pt}: \ZZ(\tii\DD_\PP) @>>>\ \ZZ(\U)$ and
consider
$$
\ZZ(\tii \DD_\PP)\subb\ \ZZ(\tii \DD_\PP)\AZ @<{\id\times a_\PP}<<
\tii \ZZ(\tii \DD_\PP)\AZ@>{\vpi^\PP \times \prr_1}>>
\ZZ(\U)\tim_{\fh^*/W} \fh^*\unr/W \dff \ZZ(\U)\unr\ \sub \ZZ(\U)
,$$ where $\prr_1:( \fh^*\unr\times \fh^*_{\Fp})/W_\PP\to
\fh^*\unr/W $ is the projection.

\slem
{\em We have a Morita equivalence of
 Azumaya algebras on $\tii \ZZ(\tii \DD_\PP)\AZ$
(in particular, on $\tii \ZZ(\U)\AZ$):
$$
(\id\times a_\PP)^*(\tii \DD_\PP)
\sim
(\vpi^\PP \times \prr_1)^*(\U)
.$$
}

Before proving the lemma we provide a slightly
more precise statement.  Notice that the variety $(\fh^*
\times \fh^*_\Fp)/W_\PP$ is a union of connected components
$H^{
%\la_\Fp
\th
}_\PP\dff (\fh^*\times 
%(W_\PP\bu \la_\Fp)
\th
)/W_\PP\cong
 \fh^*/{\mathrm{Stab}}_{W_\PP}(\la_\Fp)$
indexed by $W_\PP$-orbits 
$
\th=
W_\PP \la_\Fp\in \fh^*_\Fp/W_\PP$ in $\fh^*_\Fp$.

We also get the corresponding components in the coverings of the central
varieties:
$\tii \ZZ(\tii \DD_\PP)^{
\th
%\la_\Fp
}= \ZZ(\tii \DD_\PP)\times _{\fh^*/W_\PP}
H^{
%\la_\Fp
\th
}_\PP$, $
%\la_\Fp 
\th
\in \fh^*_\Fp/W_\PP$.

\lem \label{lemsplit}
{\em For every $
%\la_\Fp 
\th
\in \fh^*_\Fp/W_\PP$
the choice of $\tii\la\in \La$ such that
$\Wap\bu\tii\la $ equals $
%W\bu\la_\Fp
W\bu\th
$
%$\tii\la \mod \Wap$ equals $\la \mod W$
in 
$
\La/\Wap
=
\fh^*_{\Fp}/W$, defines an equivalence of Azumaya
algebras on 
$
{\tii \ZZ(\tii \DD_\PP)^{
%\la_\Fp
\th
}}\AZ
\dff
{\tii \ZZ(\tii \DD_\PP)^
%{\la_\Fp}
\th
}
\cap
{ \tii \ZZ(\tii \DD_\PP) }\AZ
$,
\begin{equation}\label{eqeq}
(\id\times a_\PP)^*(\tii \DD_\PP)\sim (\vpi^\PP \times \prr_1)^*(\U)
,
\end{equation}
which satisfies the following properties:
\newline 
%\ben 
(i) For    any
character 
 $\mu$ of the Levi group $L$ the equivalence corresponding
to $\tii\la'=\tii \la+p\mu$ is the composition of the equivalence
corresponding to $\tii \la$ and the twist by $\OO_{\PP\tw,\mu}$. 
\newline
(ii)
Given a pair of partial flag varieties $\PP=G/P \to \QQ=G/Q$ for
parabolics $P\sub Q$,
the equivalence \eqref{eqeq} is compatible
with the isomorphism \eqref{where_to} from Proposition
\ref{D-cohomology and center}.d. 
%\een 
}

\pf To establish
the equivalence of Azumaya algebras recall a well known fact
(cf., e.g., \cite{MI} IV.2, Corollary 2.6) that
given
two Azumaya algebras on a smooth variety over a field,
every equivalence between their restrictions to
 a dense open subvariety extends to the whole variety;
moreover, such an extension is unique provided that the open
subvariety has complement of codimension at least two (on every
component).

Thus it suffices to construct an equivalence between the
two Azumaya algebras on such an  open subset.

 By Proposition \ref{D-cohomology and center}.e  the complement to the preimage of $\ZZ(\U)\AZ$
under the projection $\ZZ(\tii \DD_\PP)
\to \ZZ(\U)$ has codimension at
least two.  Thus it suffices to define \eqref{eqeq} on 
$(\id\times a_\PP)^{-1} 
((\vpi^\PP_{\pt}) ^{-1}
(\ZZ(\U)\AZ ))
$. 
In view of
Proposition \ref{D-cohomology and center}.d, the construction of the
equivalence reduces to the case $\PP=\pt$.
Here
$$
\ZZ(\U)\subb\ \ZZ(\U)\AZ @<{\id\times a}<<
\tii \ZZ(\U)\AZ
@>{\vpi \times \prr_1}>>
\ZZ(\U)
\tim_{\fh^*/W} 
\fh^*\unr/W 
\dff \ZZ(\U)\unr\ \sub \ZZ(\U)
,$$ 
for
the projection
$\prr_1:( \fh^*\unr\times \fh^*_{\Fp})/W\to
\fh^*\unr/W 
$ 
and
$$
\tii\ZZ(\U)=\ZZ(\U)\times_{\ZZ_{\mathrm{HC}}}(\fh^*\unr\times \fh^*_\Fp)/W
=\fg^*\tw\times_{\fh^*\tw/W} \left(( \fh^*\unr\times \fh^*_{\Fp})/W \right)
.$$

To define the required equivalence in this case we need to define
for every $
%\la_\Fp
\th
\in \fh^*_\Fp/W$ with a fixed representative $\tii
\la\in \La$ a bimodule $M=M(\tii\la)$ over $\U$ such that:

i) As a $\zfr$ bimodule $M$ is 
%scheme theoretically 
concentrated on the diagonal, and
as a $\zhc$ bimodule it is concentrated on   
$\{(\mu,
\mu+
d\tii \la
%\th
)\mod W\ |\ 
\mu\in \fh^*\}\inj
(\fh^*/W)^2=\Spec(\zhc)^2$.

ii) The restriction of $M$ as a left $\U$ module
to $\ZZ(\U)\unr$ is  free of rank
$
%|W/{\mathrm{Stab}}_W(\la_\Fp)|
|\th|
$ (notice that the cardinality $|\th|$ of the orbit equals
the degree
of the map
$H^\th
%{\la_\Fp}
\to
\fh^*/W$).

To construct such a module we 
observe that
$^{\OO_{\tii \la}}
\tii\DD
\cong \tii\DD
$
%(with 
so that  
$
\tii\DD\ot_{\OO_\BB}\OO_{-\tii \la}
$
is a
$\tii\DD$-bimodule.
Then the $\U$-bimodule
$\tii M=\
\Ga(\BB,\tii\DD\ot\OO_{-\tii \la})
$ can be viewed as
sections
of
$\tii\DD\ot\OO_{-\tii \la}
$
on
$\ZZ(\tii\DD)$, and also (by removing a codimension two subvariety),
as sections on
$
\gt^*\tw_{\mathrm{reg}}\times_{\fh^*\tw/W}\fh^*/W 
$,
where $\gt^*_{\mathrm{reg}}$
 is the preimage of the subset of regular elements $\g^*_{\mathrm{reg}} \subset \g^*$.

The natural map $\gt^*_{\mathrm{reg}}\to
 \g^*_{\mathrm{reg}}\times _{\fh^*/W}\fh^*$ is an isomorphism,
%since the target space is smooth in view of \cite[Corollary 3.4]{BG},
%and the map induces a bijection on closed points.
thus
the Weyl group $W$ acts on $\gt^*\tw_{\mathrm{reg}}$.
It is not hard to see
that this action
is compatible with the natural action on $\La=\Pic(\BB)\iso
\Pic(\gt^*_{\mathrm{reg}})$. Moreover,  a line bundle
$\OO_{\gt^*_{\mathrm{reg}}\tw,\tii\la}$ carries an equivariant structure with
respect to the subgroup
$
{\mathrm{Stab}}_W(\tii\la)\sub\ W$.
Here,
${\mathrm{Stab}}_W(\tii\la)\cong {\mathrm{Stab}}_W(\la_\Fp)$ in $W$
for some $\la_\Fp\in\th$.
%God knows why
  This gives an action of
${\mathrm{Stab}}_W(\tii\la)$ on $\tii M$, and then we
 set $M=\tii M^{{\mathrm{Stab}}_W(\tii\la)}$. Properties (i) and (ii) are easy to check.
This yields the required equivalence; the compatibilities then follow
from the construction. \epf

\sss{
Parabolic Springer fibers
}
\label{Spring is here}
We proceed to describe the parabolic analogues of the Springer fibers
and splitting of the Azumaya algebra on them.

Recall the maps
$\fg^*@<\PR_\fg<<
\gtil^*_\PP
@>\PR>>
\fg^*\times_{\fh^*/W}
\
\fh^*/W_\PP$
from
\ref{Algebras} (here the $W$ action is {\em not} the dot action, but
the standard one).
For  $(\chi,\nu)\in \fg^*\times_{\fh^*/W} \fh^*$
define parabolic Springer fibers
$\PP_\chi$ and $\PP_{\chi,\nu}
$ as  subschemes
 $\PR_\fg\inv(\chi)$
and  $\PR\inv(\chi,W_\PP \nu)$
of $\gtil^*_\PP$.
%%%
%%% What is the image of
%%% $\BB_{\chi,\nu}@>>>\BB_{\chi,\nu}$ ? (It is connected.)
%% What? (Roma)
Via the projection
$
\gtil^*_\PP @>\tau>> \PP
$ the parabolic Springer fiber
$
\PP_{\chi,\nu}
$
can be identified
with a subscheme
$
\tau(\PP_{\chi,\nu})
$
of
$\PP$,
and
$\PP_{\chi,\nu}$ is a section
of
$\gtil^*_\PP$ over
$\tau(\PP_{\chi,\nu})$.

For $\la\in\fh^*$  denote $\nu\dff \AS(\la)\in\fh^*\tw$.
%A point
%$(\chi,\la)\in \ZZ(\U)=\fg^*\tw\tim_{\fh^*\tw/W}\fh^*/W$
%gives in this way a point
%$(\chi,\nu)\in (\fg^*\tim_{\fh^*/W}\fh^*)\tw$.
%Now
%$\PP_{\chi,\nu}\tw
%\cong
%\PP_{\chi}\tw\times_{\fh^*\tw}\la
%$.

The fiber at $ (\chi,W_\PP\bu\la) $ of the  map of central varieties
$\ZZ(\tii\DD_\PP) @>>> \ZZ[\Ga(\tii\DD_\PP)]$\ie
of the map $
\tii\fg^*_\PP\tw\tim_{\fh^*\tw/W_\PP}\ \fh^*/W_\PP @>>>
\fg^*\tw\tim_{\fh^*\tw/W}\ \fh^*/W_\PP $, is the same as  the parabolic
Springer fiber $\PP_{\chi,\nu}\tw$. Moreover, the formal
neighborhoods of the fibers of the two maps are canonically identified. In
other words, the formal neighborhood for $
\ZZ(\tii\DD_\PP)\tim_{\ZZ[\Ga(\tii\DD_\PP)]}\ (\chi,W_\PP\bu\la)
\ \sub\ \ZZ(\tii\DD_\PP)$ can be identified with the formal
neighborhood of
$ \tii\fg^*_\PP\tw 
\tim_{ \fg^*\tw\times %_{\fh^*\tw/W}\
\fh^*/W_\PP }\ (\chi,\nu) =\PP_{\chi,\nu}\tw \ \sub\
\tii\fg^*_\PP\tw$.

\scor
{\em
If  $\la\in\fh^*\PPunr $ then the Azumaya algebra obtained by restricting
$\tii
\DD_\PP$
to the formal neighborhood
${\FN(\tii \fg^*_\PP\tw \times _{\fh^*\tw/W_\PP} W_\PP\bu \la
) }$,
is equivalent to the pull-back of  an Azumaya algebra on
$\FN(\fg^*\tw \times _{\fh^*\tw/W} W\bu \la )$.

In particular, for any $\chi$ with $(\chi,W\bu\la)\in\ZZ(\U)$, the
restriction of $\tii\DD_\PP $ to the formal neighborhood of $
\PP_{\chi,\nu}\tw $ in $ \ZZ(\tii\DD_\PP) $  is a split Azumaya
algebra (here  $\nu=\AS(\la)$ as above).

An equivalence and a splitting as above can be assigned in a
canonical way to a choice of  $\mu\in \La$ such that $\la -d\mu \in
\fh^*\unr$.

The equivalence and the splitting are compatible with pull-backs
from  $\ZZ(\tii \DD_\QQ)$ to  $\ZZ(\tii \DD_\PP)$
for a map of partial flag varieties $\PP\to \QQ$,
and with replacing $\mu$ by $\mu'\in \mu +p\La$, as in Lemma
\ref{lemsplit}.
%The equivalence and a splitting are compatible with pull-backs
%to $\ZZ(\tii \DD_\QQ)$ for a map of partial flag varieties
%and with replacing $\mu$ by $\mu'\in \mu +p\La$ as in Lemma
%\ref{lemsplit}.
 \epf}

\rem\label{Spring is there} Assume that $\la$ is integral.
Then
 a splitting
bundle on the formal neighborhood of a parabolic Springer fiber 
can be constructed by a simpler method, which is completely parallel
to the corresponding construction in \cite{BMR}.

Namely, if $\la 
\in \h^*_\Fp$ is $\PP$-unramified, then
we have ${\mathrm{Stab}}_W(\la)\supset W_\PP$. We can find a character $\eta\in \La$
such that $d\eta=\la$ and ${\mathrm{Stab}}_W(\eta)\supset W_\PP$.
Each such $\eta$ defines a splitting bundle $\MM=\MM^\PP_{\chi,\eta}$
 as follows (cf. \cite{BMR}, section 5).

Consider first the case  $\eta=-\rho$.
Proposition \ref{D-cohomology and center}.d shows that
 the restriction of the
Azumaya algebra $\tii \DD_\PP$ to $\FN_{\ZZ(\tii \DD)}
[(\vpi^\PP)^{-1}(\chi,-\rho)]
\cong\
(\FN_{\gt_\PP}(\PR^{-1}(\chi,0)))\tw$ 
is canonically isomorphic to the pull-back under 
$\vpi^\PP$
of the Azumaya algebra $\U$
on $\FN_{\ZZ(\U)}(\chi,-\rho)$. We fix a splitting bundle
$\MM^\PP_{\chi,-\rho}$ of the form $\MM^\PP_{\chi,-\rho}=\
(\vpi^\PP)^*\MM_0$
where $\MM_0$ is a splitting bundle for $\U$ on
the formal neighborhood ${\FN_{\ZZ(\U)}(\chi,-\rho)}$.

Let now $\eta\in \La$ be any character as above. Then $\eta+\rho$
is a character of the Levi $L\subset P$, thus we have a line
bundle $\OO_{\PP,\eta+\rho}$ on $\PP$. We set
$\MM^\PP_{\chi,\eta}=\OO_{\PP,\eta+\rho}\otimes _{\OO_\PP} 
\MM^\PP_{\chi,-\rho}$. In view of Lemma \ref{Torsors} this is indeed
the required splitting bundle.

It is clear that given two partial flag varieties 
$\PP=G/P\to \QQ=G/Q$,
 the constructed splitting bundles
enjoy the following compatibility:
\begin{equation}\label{PQcomp}
(\pi^\PP_\QQ)^*\
\MM^\QQ_{\chi,\nu}
\ \cong\ 
\MM^\PP_{\chi,\nu}
\end{equation} 
for any $\nu\in \La$ such that ${\mathrm{Stab}}_W(\nu)\supset W_\QQ\supset W_\PP$.
This isomorphism is compatible with the actions of $\tii \DD_\PP$,
$\tii \DD_\QQ$ via the isomorphism \eqref{where_to}.
 
\medskip

We leave it as an exercise to the reader to check that the splitting
bundles produced by this construction
are also obtained by the construction of Corollary \ref{Spring is here}.
This compatibility will not be used in the paper.

\sus{
Derived categories
}
\sss{
Derived categories of sheaves supported on a subscheme
}
\lab{Derived categories of sheaves supported on a subscheme}

 Let $\AA$ be  a coherent sheaf on a Noetherian scheme $\fX$
equipped with an associative $\OO_\fX$-algebra structure.
We denote by $\mod^c(\AA)$ the abelian
category of  coherent %(i.e. locally finitely generated)
sheaves of $\AA$-modules. We also use notations
$\Coh(\fX)$ if $\AA=\OO_\fX$
and
$\mod^{fg}(\AA)$ if $\fX$ is affine.

For a subscheme $\fY$ we denote by $\mod^c_{\hattt\fY}(\AA)$ the full
subcategory of coherent $\AA$-modules supported set-theoretically in
$\fY$, i.e., killed by some power of the ideal sheaf $\II_\fY$. Then
the tautological functor $ \Db[\mod^c_{\hattt\fY}(\AA)] @>>>
\Db[\mod^c(\AA)] $ is a full embedding and the image consists of all
$\F\in \Db[\mod^c(\AA)]$ that satisfy the equivalent  conditions:\ i)
$\F$ is killed by a power of the ideal  sheaf $\II_\fY$, i.e. the
tautological arrow $\II_\fY^n\otimes _\OO \F\to \F$ is zero for some
$n$; and \ ii)  the cohomology sheaves of $\F$ lie in
$\mod^c_{\hattt\fY}(\AA)$ (see, e.g., \cite{BMR}).

We will use this
for the sheaves of algebras 
$\tii\DD_\PP$ and $\U$
over
$\tii \fg^*_\PP\tw$ and $\fg^*\tw$ respectively (or
$\ZZ(\tii\DD_\PP)$
and $\ZZ(\U)$ respectively).
For $\la\in\fh^*$ we have the
abelian categories
 $
\mod^c(\DD^\la_\PP) \sub \mod^c_{\la}(\tii\DD_\PP) \sub
\mod^c(\tii\DD_\PP) $; the second category is a particular case of the above situation where
$\fX=\ZZ(\tii \DD_\PP)$ and $\fY$ is the preimage
of $\la\mod W_\PP$ under the projection $\ZZ(\tii \DD_\PP)\to \fh^*/W_\PP$.
We also  have $ \mod^{fg}(\U^\la)\sub \mod^{fg}_{\hattt
\la}(\U) \sub \mod^{fg}(\U) $, where the second category is obtained as above
for $\fX=\ZZ(U)$ and $\fY$ being the preimage of $\la \mod W$ under the projection to $\ZZ_{HC}$.
For the corresponding triangulated categories we get
$
\Db[\mod^{c}(\DD^\la_\PP)] @>>> $ $\Db[\mod^{c}_{ \la}(\tii \DD_\PP)] $ $
\sub \Db[\mod^{c}(\tii \DD_\PP)] $,
$
\Db[\mod^{fg}(\U^\la)] @>>> $ $\Db[\mod^{fg}_{\hattt \la}(\U)] $ $
\sub \Db[\mod^{fg}(\U)] $.
We get similar categories for $(\chi,\la)\in\ZZ(\U)$.

\sss{
The global section
functors on $D$-modules
}
\lab{global section functors}
We choose appropriate derived global section functors for
$\mod^c(\DD^\la_\PP)
\sub
\mod^c_{\la}(\tii\DD_\PP)
\sub
\mod^c(\tii\DD_\PP)
$.
We start with
the functor on quasi-coherent sheaves
$
\Ga:\mod^{qc}(\OO_\PP)
@>>>\Vect
$.
The
map
$
\tii \U
@>>>
\Ga(\tii \DD)
$
gives functors
$$
\mod^{qc}(\tii\DD_\PP)
@>{\Ga_{\tii\DD_\PP}}>>
\mod(\U)
,\
\mod^{qc}_{\hattt \la}(\tii\DD_\PP)
@>{\Gamma_{\tii\DD_\PP,\la}}>>
\mod_{\hattt \la}(\U),
\aand
\mod^{qc}(\DD^\la_\PP)
@>{\Gamma_{\DD^\la_\PP}}>>
\mod(\U^\la)
.$$
We  derive
these functors
in the above abelian categories
(this can be done because the categories
of modules have direct limits), and this gives
$$
\RGa:\
 \Dd[\mod^{qc}(\OO_\PP)]
@>>>
\Dd[\Vect]
,\ \ \ \
\RGa_{\tii\DD_\PP}:\ \Dd(\mod^{qc}( \tii\DD_\PP))
@>{ }>>
\Dd(\mod ( \U))
,
$$
$$
\RGa_{\tii\DD_\PP,\la}:\ \Dd(\mod^{qc}_{\la}
(\tii\DD_\PP))
@>{ }>>
\Dd(\mod_{\la}( \U))
,\ \ \
\RGa_{\DD^\la_\PP}:\ \Dd(\mod^{qc} (\DD^\la_\PP))
@>{}>>
\Dd(\mod( \U^\la))
.$$
The following  properties
 were explained in \cite[Section 3.1.9]{BMR} for the particular
case $\PP=\BB$; the same arguments apply in our present generality.

(0)
All of these functors have
finite homological dimension, so they give functors between
bounded derived categories.

(1)
The  derived
functors commute with
 the  forgetful functors, for instance
$
\mbox{Forg}^{\tii \U}_\k
\ci
\RGa_{\tii\DD}
\cong \
 \RGa  \ci \mbox{Forg}^{\tii\DD}_\OO
$
canonically
for
the  forgetful functors
$
{\mbox{Forg}^{\tii\DD}_\OO}:\mod^{qc}(\tii\DD)
\to
\mod^{qc}(\OO)
$,
$
\mbox{Forg}^{\tii \U}_\k:
\mod(\tii \U)\to \Vect _\k
$.

(2)
The above  (derived)  functors of global sections
preserve coherence, i.e.
$\RGa_{\tii\DD}$ sends the full subcategory  $\Db[\mod^c(\tii\DD)]
\subset \Db[\mod^{qc}(\tii \DD)]$ into the full subcategory
$ \Db[\mod^{fg}(\tii \U)]\subset
\Db[\mod(\tii \U)]$, etc.

\sus{
Equivalences
}
\lab{Equivalences}
The following Theorem
%\ref{Localization theorem}
 is the main result of this section.
We analyze the necessary conditions on $\la$
in Lemma
\ref{Applicability of equivalence}.
In particular, we find that for each Harish-Chandra character
the corresponding category of $\U$-modules has at least one localization\ie
there is some $\la$ in the corresponding $W$-orbit
in $\fh^*$,  and some partial flag variety
$\PP$ such that
all parts of the theorem below
apply (Lemma \ref{Applicability of equivalence}c).

We say that $\la$ is
{\em $\PP$-regular} if
the stabilizer $(W,\bu)_\la$ lies in $W_\PP$.
We say that a Harish-Chandra character has a singularity of $\PP$-type
if the corresponding $W$-orbit in $\fh^*$ contains $\la$ with
the stabilizer $W_\la=W_\PP$.

\theo
\lab{Localization theorem}
{\em
Consider a $\PP$-regular
$\la\in\fh^*$.

(a)
The global section functors
 provide equivalences
 of triangulated categories
\begin{equation}
\label{equivD}
\RGa_{\tii\DD_\PP,\la}:
\Db[\mod^c_\la(\tii\DD_\PP)]
\con
\Db[\mod^{fg}_\la( \U)].
\end{equation}
\begin{equation}
\label{equivD0}
\RGa_{\DD^\la_\PP}:
\Db[\mod^c(\DD^\la_\PP)]
\con
\Db[\mod^{fg}(\U^\la)].
\end{equation}

(b)
For any $\chi\in\fg^*\tw$, compatible with
$\la$
(i.e., $(\chi,\la)\in \ZZ(\U)$),
the  above equivalences
restrict to equivalences of
full subcategories with a generalized Frobenius character $\chi$
$$
\Db[\mod^c_{ {\la,\chi}}(\tii \DD_\PP)]
\cong\
\Db[\mod^{fg}_{ {\la,\chi}}( \U)]
\aand
\Db[\mod^c_{\chi}(\DD^\la_\PP)]
\cong\
\Db[\mod^{fg}_{\chi}(\U^{\la})]
.$$

(c)
If $\la$ is also a
$\PP$-unramified weight, then
there are
equivalences  (set $\nu=\AS(\la)$)
\fttt{
The equivalences depend on the choice of the splitting bundle
for $\U$ on the formal neighborhood
of the point $(\chi,\la)$ in
$
\ZZ(\U)$.
}
$$
\Db[\mod^{fg}_{\hattt (\la,\chi)}(\U)]
\cong\
\Db[\mod^{c}_{\hattt (\la,\chi)}(\tii \DD_\PP)]
\cong
\Db[\CC oh_{\PP_{\chi,\nu}\tw}
(\tii \fg^*_\PP\tw)]
;$$
$$
\Db[\mod^{fg}_{\hattt\chi}(\U^\la)]
\cong\
\Db[\mod^c_{\hattt\chi}(\DD^\la_\PP)]
\cong\
\Db[\CC oh_{\PP_{\chi,\nu}\tw}
(\tii \fg^*_\PP\tw\tim_{\fh^*\tw}\nu)]
.$$
}

The proof will be given in section \ref{19}.

\srem It would be interesting to describe
 explicitly objects in $\Db[\CC oh_{\PP_{\chi,\nu}\tw}
(\tii \fg^*_\PP\tw\tim_{\fh^*\tw}\nu)]$
corresponding to (at least some)  irreducible $\U$ modules,
and also the locally free sheaves
on the formal neighborhood of $\PP_{\chi,\nu}\tw$
corresponding to indecomposable projective pro-objects in
$\mod^{fg}_{\hattt\chi}(\U^\la)$, $\mod^{fg}_{\hattt (\la,\chi)}(\U)$.
 For $\PP=\BB$ this
has been done in a few cases in \cite[5.3.3]{BMR}. More examples can 
be computed using Lemma and Remark
\ref{coherent translation of intertwining functors},
or
Remark \ref{pushpull} below.

\lem
\lab{Applicability of equivalence}
{\em
(a)
All parts of the theorem
apply precisely when
$\la$ is $\PP$-regular and $\PP$-unramified. This is equivalent to
the following relation of  stabilizers:
$$
W_\la\ =\ W_\PP\cap W_{\AS(\la)}
.$$
(b)
A
sufficient condition
is given by $W_\la=W_\PP$.
For integral
$\la$ this is an equivalent condition.

(c)
Each $W$-orbit in $\fh^*$
contains $\la$ such that the singularity of $\la$ is of type
$\PP$
for some partial flag variety $\PP$\ie
$W_\la=W_\PP$.

(d) If $\la$ is the  differential of a character $\nu\in\La$ that
lies in the closure of the fundamental alcove $A_0$, but not on any
face  corresponding to an affine coroot $\ch\al\in \Iaff\setminus
\Sigma$, then $W_\la=W_\PP$ is satisfied when $\PP$ is the partial
flag variety corresponding to the set  of hyperplanes that contain
$\nu$. }

\pf
(a) is clear from definitions since
$\la$ is $\PP$-regular if
$(W,\bu)_\la\sub W_\PP$,
and
$\PP$-unramified if
$(W_\PP,\bu)_\la=(W_\PP)_{{\AS}(\la)}
$\ie
$(W,\bu)_\la\cap W_\PP=W_{{\AS}(\la)}\cap W_\PP$.

(b) The first claim follows from $W_\la\sub W_{{\AS}(\la)}$. For an
integral weight $\la\in\fh^*_\Fp$ the stabilizer of ${\AS}(\la)=0$ is
$W$, so the condition in (a) reduces to $W_\la=W_\PP$.

(c) follows  from Lemma \ref{connected centraliser}.

%is a standard observation. A choice of a Cartan  subgroup $T\sub
%G$ identifies $\Th$ with $\Th'\sub \ft^*$. The centralizer of any
%$\la'\in\Th'$ in $G$ is a Levi subgroup $L\subb T$ by Lemma
%\ref{Restrictions on G and p}. Now we can choose a Borel subgroup
%$T\sub B\sub G$ so that $P=LB$ is a parabolic subgroup, and let
%$\la\in\fh^*$ correspond to $\la'\in\ft^*$ via the identification
%$\io_\fb$ (\ref{Semisimple group G}).

(d)
For any character $\nu$  the quotient map
$\Wap@>>>W$
gives
an isomorphism of stabilizers
$(\Wap,\bu)_\nu\con
(W,\bu)_{d\nu}$.
Because of the assumption on $p$ one has
$(\Waff,\bu)_\nu
=
(\Wap,\bu)_\nu
$.
Finally, by Chevalley's
theorem
$(\Waff,\bu)_\nu=W_\PP$.
Now, the inclusion
$W_\PP\sub (W,\bu)_{d\nu}$
is equality since the orders of groups are the same.

\scor
{\em
If in the conditions of theorem \ref{Localization theorem}
the weight $\la$ is $\PP$-regular and $\PP$-unramified,
then
there is a natural isomorphism of
Grothendieck groups $ K(\U^\la_\chi)\cong K(\PP_{\chi,\nu}) $. In
particular, the number of irreducible $\U^\la_\chi$-modules is the
rank of $K(\PP_{\chi,\nu})$.
}

\srem The method of \cite[Section 7]{BMR} can be used to show
that for large $p$ the rank of $K(\PP_{\chi,\nu})$ equals the sum
of Betti numbers of the corresponding parabolic Springer fiber
$\PP_{\chi',\nu'}$ in characteristic zero; here
$\nu'$ is an element in the Lie algebra over a characteristic zero field whose
centralizer is the Levi subgroup of the same type as the centralizer
of $\nu$, and $\chi'$ corresponds to $\chi$ under the bijection
arising from the Bala-Carter classification.

\sss{
Reduction of the theorem to the part (a)
}
\label{reduction to a}
A proof of the part (a) will be given in section \ref{19}.
Here we assume (a) and prove (b) and (c).

(b)
$\OO(\fg^*\tw)$ acts on the categories $\mod^c(\tii\DD_\PP)$,
$\mod^{fg}(\U)$
etc.,
and on their derived categories. This means an algebra homomorphism
from $\OO(\fg^*\tw)$ to the algebra of natural endomorphisms of
the identity functor (often called the center of the category).
The homomorphism takes $f\in\OO(\fg^*\tw)$ into a natural transformation
defined by $m\mapsto fm$ for an object $M$ and $m\in M$. 
The  equivalences in (a) are equivariant under $\OO(\fg^*\tw)$
and therefore they restrict to the full subcategories of
objects on which the
Frobenius center acts by the
generalized character $\chi$
(see \ref{Derived categories of sheaves supported on a subscheme}).

(c) Consider $\tii\DD_\PP$-modules as sheaves on $\ZZ(\tii\DD_\PP)$.
They have  a generalized central character $(\chi,\la)$ precisely if
they are set theoretically supported on
$\ZZ(\tii\DD_\PP)\tim_{\ZZ(\U)}(\chi,\la) =\PP\tw_{\chi,\nu}$. Since
$\la$ is $\PP$-unramified,
%was $\PP$-split!
 the restriction of
$\tii\DD_\PP$ to the formal neighborhood of the Springer fiber
$\PP_{\chi,\nu}\tw$ in $\ZZ(\tii\DD_\PP)$ (the same as the formal
neighborhood $\FN_{\PP_{\chi,\nu}}(\tii \fg^*_\PP)\tw$ of
$\PP_{\chi,\nu}\tw$ in $\tii \fg^*_\PP\tw$), is a trivial Azumaya
algebra by the Corollary \ref{Spring is here}.
 Therefore any choice
of a splitting vector bundle $\MM^\PP_{\chi,\la}$ on this formal
neighborhood provides equivalences
%\fttt{ Vector bundle $\MM^\PP_{\chi,\la}$ carries a
%$\tii\DD_\PP$-action and this identifies the restriction of
%$\tii\DD_\PP$ to $\FN_{\PP_{\chi,\nu}}(\tii \fg^*_\PP)\tw$ with
%$\EEnd_\OO(\MM^\PP_{\chi,\la})$. }
 of abelian
categories
$$
\Coh_{\PP_{\chi,\nu}\tw}(\tii\fg^*_\PP\tw) \con\
\Coh_{\PP_{\chi,\nu}\tw}(\ZZ(\tii\DD_\PP))
@>\MM^\PP_{\chi,\la}\ten->\cong>\ \mod^c_{\chi,\la}(\tii\DD_\PP) .$$
This proves the first claim. We get the second one similarly, as
$\DD_\PP^\la$-modules are sheaves on
$\ZZ(\tii\DD_\PP)\tim_{\fh^*/W_\PP} \la$.

\sus{
Localization as the left adjoint of global sections
}
\lab{Localization as the left adjoint of global sections}

\sss{
Localization functors
}
\lab{Localization}
We start with the localization functor
${\mathrm{Loc}}_\PP$ from finitely generated
$\U$-modules to coherent $\tii\DD_\PP$ modules,
$$
{\mathrm{Loc}}_\PP(M)\
\dff
\tii \DD_\PP
\otimes_{\U}
M
.
$$
For each $\la\in\fh^*$ it
restricts to a functor
$$
{\mathrm{Loc}}^\lambda_\PP: \mod^{fg}(\U^\lambda) \to
\mod^c(\DD^\lambda_\PP),\
{\mathrm{Loc}}^\lambda_\PP(M)
\dff
\DD^\la_\PP\ten_{\U^\la}M
.
$$

Since $\U$ has finite homological dimension, the functor 
${\mathrm{Loc}}_\PP$ has a left derived functor
$
\Db[\mod^{fg}(\U)]
@>{\LL_\PP}>>
\Db[\mod^c(\tii\DD_\PP)]
$.
Fix $\lambda\in \fh^*$, for any $M\in
\Db[\mod_\lambda^{fg}(\U)]$ 
the action of $\zhc=S(\fh)^W$ on
$\LL_\PP(M)$ factors through an ideal $I$,
the power of the maximal ideal corresponding
to the orbit $W\bu\la$.
The finite-dimensional algebra $S(\fh)^{W_\PP}/IS(\fh)^{W_\PP}$,
which also acts on  $\LL_\PP(M)$,
has a primitive idempotent for each orbit
$W_\PP\bu \mu \sub W\bu \lambda$.
These primitive idempotents give
a canonical decomposition
$\LL_\PP(M)=
\pl_{W_\PP\bu \mu \sub W\bu \lambda}\
\LL^{W\bu\lambda\to W_\PP\bu\mu}_\PP(M)$
with $\LL^{W\bu\lambda\to W_\PP\bu\mu}_\PP(M)\in
\Db[\mod^c_\mu (\tii \DD_\PP)]$.
Localization with the generalized character $\la$
is the functor
$
\LL^{\hatt \lambda}_\PP
\df
\LL^{W\bu\lambda\to W_\PP\bu\lambda}_\PP:
\Db[\mod^{fg}_\la(\U)]
@>>>
\Db[\mod^c_\la(\tii\DD_\PP)]
$.

The functor 
${\mathrm{Loc}}^\la_\PP$ also  has a left derived functor
$\LL^\la_\PP:\ \Dmin(\mod^{fg}(\U^\la))\to \Dmin(\mod^c(\DD^\la_\PP)),\
\LL^\lambda_\PP(M)=\DD^\la_\PP\Lten_{\U^\la}M
$. The algebra $\U^\lambda$ may have infinite homological
dimension\fttt{
For regular $\la$  finiteness of homological
dimension follows from Theorem~\ref{Localization theorem}.
%was observed in \cite{BMR}. 
For singular $\la$ the homological dimension is infinite. 
For instance, for $\la =-\rho$ the algebra $\U^\lambda$ is a split Azumaya algebra 
over the Frobenius twisted nilpotent cone.} % \cite{BMR}.}
so a'priori $\LL^\la_\PP$ need not  preserve  the
bounded derived categories.

\lem
\lab{adjoints and forgettings}
{\em
(a)
The functor $\LL_\PP$ is left adjoint to
$\RGa_{\tii\DD_\PP}$.

(b)
The functor
$\LL^{\hatt \la}_\PP $ is left adjoint to
$\RGa_{\tii\DD_\PP,\la}$.

(c)
The functor $\LL^\la_\PP$ is left adjoint to the functor
$
\Dmin (\mod^c(\DD^\la_\PP))
@>{\RGa_{\DD^\la_\PP}}>>
\Dmin(\mod^{fg}(\U^\la))$.

(d) For $\PP$-regular $\la$ the localizations at
$\la$ and at the generalized
character $\la$ are compatible, i.e., for the obvious functors
$\Dmin(\mod^{fg}(\U^\la))
@>{i}>>
\Dmin(\mod_\la^{fg}(\U))$ and
$
\Dmin(\mod^c(\DD^\la_\PP))
@>{\io}>>
\Dmin(\mod_\la^c(\tii\DD_\PP))$, there is a  canonical isomorphism
$\iota \circ \LL^\la\cong \LL^{\hatt \la} \circ i$.
This isomorphism is compatible with the adjunction arrows, i.e.
for $M\in \Dmin(\mod^{fg}(\U^\la))$, $\F\in \Dmin (\mod^c(\DD^\la_\PP))$
the composition:
\begin{multline*}
 \Hom(M,R\Gamma_\la\F)\cong \Hom(\LL^\la M,\F)\to
\Hom (\iota \LL^\la M, \iota \F) \cong \Hom (\LL^{\hatt \la}iM,
\iota \F)\\ 
\cong \Hom (iM, R\Gamma_{\hatt \la}\iota \F)\cong \Hom (iM, i
R\Gamma_{\la}\F)
\end{multline*}
coincides with the map induced by functoriality of $i$.
%in the
%obvious sense.
}

\pf
(a)
One checks from the definitions
that the functors between abelian
categories  form adjoint pairs.
Since $\mod^{qc}(\tii
\DD_\PP)$ (respectively, $\mod(\U)$) has enough injective (respectively,
projective) objects, and the functors $\Gamma_\PP$, ${\mathrm{Loc}}_\PP$ have bounded
homological dimension,
 it follows that their derived functors form
an adjoint pair of functors between bounded derived categories.
This adjunction restricts to the required adjunction for
$\LL_\PP$ and $\RGa_{\tii\DD_\PP}$.

(b)  follows from (a)
%the adjunction for
%$\LL^{\hatt \la}_\PP $,
%$\RGa_{\tii\DD_\PP,\la}$ 
%follows 
by passing to summands.

(c)
The proof 
%for  the last pair
%$\LL^\la_\PP$,
%${\RGa_{\DD^\la_\PP}}
%$ 
is analogous to the proof for the first pair.

(d) For a module $M$ with character $\la$, the $S\fh$-module
$\LL_\PP M$ is supported on the subscheme
$\fh^*/W_\PP\tim_{\fh^*/W}\la\sub\fh^*$. Since $\la$ is
$\PP$-regular the map ${\fh^*/W_\PP}@>>>\fh^*/W$ is unramified at
$W_\PP\bu\la$. Then the point $\la$ is a connected component of
$\fh^*/W_\PP\tim_{\fh^*/W}\la$, and the corresponding  summand of
$\LL_\PP M$ is ${\mathrm{Loc}}^\la_\PP M$.\ \ \epf

%{\bf ADD WHY THE FUNCTORS ON DERIVED CATEGORIES ARE ADJOINT!!!!!!!!!!!}

\cor
{\em
If $\la$ is
$\PP$-regular
the functor $\LL^\la_\PP$ sends the bounded derived category $
\Db[\mod^{fg}(\U^\la)]$ to $ \Db[\mod^{c}(\DD_\PP^\la)]$. \epf
}

\sss{
Localization is fully faithful
}
This is one way to view the following claim
(if a functor $L$ between triangulated categories
has a right adjoint $R$ then $L$ is fully faithful if and only if 
the map  $\id@>>>R\ci L$ is an isomorphism).

\spro
\lab{compid}
{\em
(a)
The composition $\RGa_{\tii\DD_\PP} \circ \LL_\PP:
\Db[\mod^{fg}(\U)] \to \Db[\mod^{fg}( \U)]$ is isomorphic to the
functor $
M\mapsto
M\ten_{ S(\fh)^{W}   } S(\fh)^{W_\PP}$.

(b)
For $\PP$-regular $\la$ the
adjunction map
$\id\to \RGa_{\tii\DD_\PP,\la} \circ
\LL^{\hatt\la}$ is an isomorphism
on $\Db[\mod^{fg}_\la(\U)]$.

(c)
For any $\la$,
the adjunction map is an isomorphism
$\id@>>>\RGa_{\DD^\la_\PP}\ci \LL^\la_\PP$
on $\Dmin(\mod^{fg}(\U^\la))$.

}

\pf
(a)
For any $\U$-module $M$ the action of $\U$ on
$\Gamma_{\tii\DD_\PP} (\LL_\PP(M))$
extends to the action of
$\Gamma (\tii \DD_\PP)= \tii \U^{W_\PP}$.
So, the adjunction map
$M\to \Gamma_{\tii\DD_\PP} (\LL_\PP(M))$
extends to
$S(\fh)^{W_\PP}\otimes _{S(\fh)^W} M
=\tii \U^{W_\PP} \ten_\U M
\to \Gamma_{\tii\DD_\PP} \circ \LL_\PP (M)$.
Proposition \ref{D-cohomology and center}.b
implies that if $M$ is a free module
then this map is an isomorphism, while higher derived functors
$\Rr^i\Ga_{ \tii\DD_\PP } ( \LL_\PP(M) )$,\ $i>0$, vanish.
This yields  (a). Statement
(c) is proved in the same way using the second claim in
Proposition \ref{D-cohomology and center}.b.

To deduce (b) observe that for a $\PP$-regular
$\la$ and  $M\in \Db[\mod_\la^{fg}(\U)]$, we have 
a canonical decomposition of $M\otimes _{S(\fh)^W} S(\fh)^{W_\PP}$
into summands $M\ten_\k\OO(C)$
over
 connected components $C$ of
$\fh^*/W_\PP\tim_{\fh^*/W}\la$.
Since $\la$ is $\PP$-regular, one such component is $\la$.
The corresponding component of the
adjunction morphism viewed as a map
$M\to \pl_{C}\ M\ten_\k\OO(C)$
is $\id_M$.
Now the claim follows
since the corresponding summand can be viewed as
$\RGa_{\tii \DD_\PP,\la}
(\LL^{\hatt\la}_\PP(M))
$.
\epf

\sus{
Calabi-Yau categories
}
The remaining ingredient of the proof of the localization theorem
is the use of properties of
the Calabi-Yau subclass
of triangulated
categories. These we recall here,
for more detail see \cite{BMR}, \cite{BK}
or the  original paper \cite{MuMu}.

\sss{Relative Serre functors}
\lab{Relative Serre functors}
Let $\O$ be a finite type commutative algebra over a field; and
let $D$ be an $\O$-linear triangulated category.
A structure of an $\OO$-triangulated category on $D$
is a functor $\RHom_{D/\OO}: D^{op}\times D\to \Db[\mod^{fg}(\OO)]$,
together with a functorial isomorphism
$\Hom_D(X,Y)\cong
\Hh^0(\RHom_{D/\OO}(X,Y))$.
By an {\it $\OO$-Serre functor } on $D$
we will mean an
auto-equivalence $S:D\to D$
together with
a natural (functorial) isomorphism
$\D_\OO[\RHom_{D/\OO}(X,Y)]\cong
\RHom_{D/\OO}(Y,SX)$
for all $X,Y\in D$.
($\D_\O$ denotes the Grothendieck
duality for $\OO$-modules.)
An  $\OO$-triangulated category will be called {\it Calabi-Yau} if
for some $n\in \Z$
the shift functor $X\mapsto X[n]$
admits a structure of an $\OO$-Serre functor.

\sss{
Calabi-Yau categories have no
retracts
}
The following lemma appears in our first paper
\cite[Lemma 3.5.2]{BMR}, a similar argument has appeared earlier
\cite[Theorem 2.3]{MuMu}.
One can summarize it by:  
{\em a retract of a Calabi-Yau category
is a summand.} 

\slemm
\lab{abstract}
{\em
Let
$D$ be a Calabi-Yau $\OO$-triangulated category
for  some
commutative finitely generated algebra $\OO$.
Let $C$ be a triangulated category. 
Then a sufficient condition for
a triangulated functor $L:C\to D$
to be  an equivalence is given by

(i) $L$ has a right adjoint functor $R$ and the adjunction
morphism $\id\to R\circ L$ is an isomorphism, and

(ii) $R$ is not zero on any summand of $D$.
}

\pf
(i) implies that $L$ is a full embedding\ie we can consider $C$ as a full subcategory of $D$. 
Then $\id\cong R\ci L$ can be thought of
as a retraction from $D$ to $C$, however a retract of a Calabi-Yau
category is a summand:  the  retract structure
shows that any $d\in D$ is in an exact triangle
$c
@>>>
d
@>>>
c'$ with $c\in C$ and $c'\in C^\perp$.
Recall that $C^\perp$ is a full subcategory of $D$ with
objects $x$ such that $\Hom_D(c,x)=0$ for all $c\in C$.
Similarly, ${\,^\perp}C$ is a full subcategory of $D$ with
objects $x$ such that $\Hom_D(x,c)=0$ for all $c\in C$.
Now the Calabi-Yau property exchanges left and right
$C^\perp=\ ^\perp C$,
and then $D\cong C\pl C^\perp$.
Now (ii) implies $C^\perp=0$.

Another useful simple fact is the following lemma, which has been known for
$\CA= \OO$ \cite[Lemma 4.2]{MuMu}.

\lemm
\lab{indecompo}
{\em Let $X$ be a connected scheme quasiprojective over a field $\k$,
and let $\CA$ be an algebra vector bundle over  $X$ which is
generically Azumaya. Then the category
$\Db[\mod^c(\CA)]$ is indecomposable. Moreover, if $Y\subset X$ is
a connected closed subset then $\Db[\mod^c_Y(X,\CA)]$ is
indecomposable.
}

\pf
This is already proved when $\CA$ is Azumaya\cite[Lemma 3.5.3]{BMR}.
The proof is valid in this more general case.
\epf

\sus{
Frobenius algebras
}
The main goal here  is to prove that
$\Db[\mod^c(\U)]$ is
Calabi-Yau
over $\fg^*\tw$
(see \ref{U is Calabi-Yau}.c).
A  {\em Frobenius} algebra  structure
on
a locally free  algebra sheaf $\AA$ over
a scheme $X$ is a functional
$
\tau\in\HHom_{\OO_X}(\AA,\OO_X)
$
which is nondegenerate
in the sense that
the map
$\tii\tau:\AA\to \HHom_{\OO_X}(\AA,\OO_X),\ (\tii\tau a)b=\tau(ab)$, is an
isomorphism.\fttt{
%Frobenius structures are the same as isomorphisms
%$T:\AA\con\HHom_{\OO_X}(\AA,\OO_X)$ since such $T$
%gives  a  map $\tau:\pi_*\AA\to\OO_X$
%(as the adjoint of the inclusion $\OO_X\sub\AA$),
%such that  $\tii\tau=T$.
Frobenius structures are the same as 
right $\AA$-module isomorphisms
$T:\AA\con\HHom_{\OO_X}(\AA,\OO_X)$ since such $T$
gives  a  map $\tau:\AA\to\OO_X$
via $\tau a = (Ta) 1_\AA$
such that  $\tii\tau=T$.
}
We say that
$\AA$ is Frobenius if it has a Frobenius structure.
We say that a map of schemes $X@>f>>Y$ is Frobenius if it is finite and
$f_*\OO_X$ is a Frobenius algebra over $\OO_\YY$.

\lem \lab{Frobenius property} {\em (a) Let $\AA$ be a sheaf of
algebras over a  normal  variety $X$ which is locally free of finite
rank  over $\OO_X$. If $\AA$ is Frobenius on an open set 
with a complement of codimension at least two (on each component),
then $\AA$ is Frobenius everywhere.

(b)
Let $\OO$ be a finite type commutative algebra over $\k$.
Let $Y$ be a variety over $\k$ equipped with a
projective morphism
$\pi:Y\to \Spec(\OO)$ then
$D=\Db[\Coh_Y]$
is $\OO$-triangulated by
$\RHom_{D/\OO}(\F,\G) \dff
\Rr\pi_* \RHHom(\F,\G)$.
If $Y$ is also smooth and quasiprojective,
then for any Frobenius  algebra
$(\BB,\tau)$ on
$Y$,
$\Db[\mod^c(\BB)]$ has a
natural structure of an $\O$-triangulated category, 
and the functor  $\F\mapsto
\F\otimes \om_Y[\dim Y]$ has a natural structure of an $\OO$-Serre functor.
In particular, if $Y$ is a Calabi-Yau manifold
(i.e.,
$\om_Y\cong \OO_Y$)
then
the
$\O$-triangulated category
$\Db[\mod^c(\BB)]$ is Calabi-Yau.

(c)
Any Azumaya algebra $\CC$ over  a scheme $Z$ has a canonical
Frobenius structure given by
the reduced trace.

(d)
Let
$X@>\pi>>Y$
be
a  finite flat map of smooth varieties with trivial canonical
classes.
For any $\OO_X$-algebra
$\AA$
there is a canonical bijection between the set
of Frobenius structures on $\AA$ over $X$ and 
the set of Frobenius structures on $\pi_*\AA$
over $Y$.

(e)
Any finite map of smooth varieties with trivial canonical classes
$X@>f>>Y$ has a canonical
Frobenius structure.
}

\pf
(a)
Restriction of $\AA$ to  an open $X_0\sub X$ 
with the complement of
codimension $\ge 2$, carries a non-degenerate
Frobenius functional
$\tau_0$.
Since
$X$ is
normal
the section $\tau_0$ of $\HHom(\AA,\OO_X)$ on $X_0$
extends to
$\AA
@>\tau>>\OO_{X}
$ on $X$. Since $\tau$ is non-degenerate
outside of codimension two, it is
 nondegenerate on $X$.

(b)
The first claim is obvious.
To prove the second claim, let us suppose for a moment
that $\FF$, $\GG$ are  locally projective
$\BB$-modules of finite rank.
The trace  map
$\EEnd_\BB(\FF)@>tr_\FF>>\BB$
leads to the pairing
$$
\HHom_\BB(\FF,\GG)\ten \HHom_\BB(\GG,\FF)
@>>>
\EEnd_\BB(\FF)
@>\tau\ci tr_\FF>>
\OO_Y
$$
which is non-degenerate\ie it gives
$
\HHom_\BB(\GG,\FF)
\con
\HHom_{\OO_Y}[\HHom_\BB(\FF,\GG),\OO_Y]
$.
>From now on let $\FF$, $\GG$
be arbitrary objects of $\Db[\mod^c(\BB)]$.
Under our assumptions
$\Db[\Coh(\BB)]$ is generated by locally projective modules.
We conclude that there is
an isomorphism
$
\RHHom_\BB(\GG,\FF)
\con
\RHHom_{\OO_Y}[\RHHom_\BB(\FF,\GG),\OO_Y]
$
for any $\FF,\GG\in\ \Db[\mod^c(\BB)]$.
Since
 $\DD_Y\cong \RHHom_{\OO_Y}(-,\om_Y[\dim(Y)])$
for a smooth $Y$
we obtain
$$
\D_Y\RHHom_\BB(\FF,\GG)
\cong
\RHHom_{\OO_Y}[\RHHom_\BB(\FF,\GG),\OO_Y)] \ten \om_Y[\dim(Y)]
$$
$$
\cong
\RHHom_\BB(\GG,\FF)\ten\om_Y[\dim(Y)]
\cong
\RHHom_\BB(\GG,\FF\ten\om_Y[\dim(Y)])
.$$
Finally, since
Grothendieck-Serre duality commutes with proper direct images
$$
\D_\OO [\Rr\pi_*\RHHom_\BB(\FF,\GG)]
\cong
\Rr\pi_*[\D_Y \RHHom_\BB(\FF,\GG)]
\cong
\Rr\pi_*\RHHom (\GG, \FF\otimes \om_Y[\dim Y])
.$$

(c)
After an  etale
base change
$\tii Z\to Z$, the algebra
$\CC$ is
isomorphic to a matrix algebra
 $\EEnd_{\OO_{\tii Z}}(\EE)$
for some vector bundle
$\EE$.
The trace
$
\EEnd_{\OO_{\tii Z}}(\EE)@>>>\OO_{\tii Z}
$
descends to
the so called {\em reduced trace}
$
\CC
@>>>
\OO_{Z}
$,
which is clearly non-degenerate.
(Actually all invariant polynomials descend.)

(d) Since the canonical classes of $X,Y$ are trivial,
 the Grothendieck-Serre duality functor for $X,Y$ is given by
 $\F\mapsto \RHom(\F,\OO)[d]$ where $d=\dim X=\dim Y$ (we assume
 without loss of generality that $X$ and $Y$ are equidimensional).
 Since  the Grothendieck-Serre duality
 commutes with proper, in
particular, finite, push-forwards we have canonical isomorphisms

$$
\pi_*\HHom_{\OO_X}(\AA,\OO_X) =\HHom_{\OO_X}(\pi_*\AA,\OO_Y) .$$ The
isomorphism is compatible with (both left and right)  $\pi_*\AA$
action by functoriality of the isomorphism between the two
compositions of proper push-forward and the Grothendieck-Serre
duality.
 This gives a bijection
of (iso)morphisms $\AA\to \HHom_{\OO_X}(\AA,\OO_X)$ and $\pi_*\AA\to
\HHom_{\OO_Y}(\pi_*\AA,\OO_Y)$.

(e) is the  case $\AA=\OO_X$ of (d). \epf

\cor
\lab{U is Calabi-Yau}
{\em
(a)
$\U$ is Frobenius over $\fg^*\tw$. In particular,  $\Db[\mod^c(\U)]$ is
Calabi-Yau over $\fg^*\tw$.

(b)
The algebras
$\tii D_X$ constructed  from  torsors
in \ref{Torsors} are Frobenius over their
$p$-central varieties.
In the particular case
of the algebras
$\tii\DD_\PP$
associated to partial flag varieties $\PP$,
the categories
$\Db[\mod^c(\tii\DD_\PP)]$ are Calabi-Yau over $\fg^*\tw$.

(c) $\U^\la$ is Frobenius over $\fg^*\tw\tim_{\fh^*\tw/W}{\AS}(\la)$
when  $\la\in\fh^*$ is unramified for ${\AS}:\fh^*/W\to\fh^*\tw/W$.
Also, $\DD^\la_\PP $ is Frobenius over its $p$-center when $\la$ is
unramified for $\fh^*/W_\PP\to\fh^*\tw/W_\PP$. In both cases the
derived categories of coherent modules are Calabi-Yau over
$\fg^*\tw\tim_{\fh^*\tw/W}{\AS}(\la)$. }

\pf
(a) The ``regular part''
$\ZZ_r=\fg^*_{\mathrm{reg}}\tw\tim_{\fh^*\tw/W}\fh^*/W$
of the  central variety $\ZZ(\U)$ is smooth since
$\fh^*/W$ is smooth by \cite{De} and the map
$\fg^*\tw_{\mathrm{reg}}\to {\fh^*_r\tw/W}$ is also smooth
\cite[Corollary 3.4]{BG}.
Since in
$\ZZ(\U)$
the Azumaya locus coincides with the
smooth locus
\cite[Theorems 2.5, 2.6]{BG}, \cite{BrGo},
we see that 
the algebra $\U$ is Azumaya over $\ZZ_r$,  
and therefore also Frobenius.
The map ${\AS}/W:\fh^*/W\to\fh^*\tw/W$
is canonically Frobenius by
part (e) of Lemma \ref{Frobenius property},
and then
$\U$ is also Frobenius over $\fg^*_{\mathrm{reg}}\tw$
by the part (d) of the same lemma.
The complement of $
\fg^*_{\mathrm{reg}}
\sub \fg^*$ has codimension  at least three
\cite[Proposition 3.2]{BG},
so  $\U$ is Frobenius over
all of $\fg^*\tw$.
Then
$\Db[\mod^{fg}(\U)]$ is Calabi-Yau
over
$\fg^*\tw$ by the Lemma \ref{Frobenius property}.b.

(b)
The algebras
$\tii D_X$  are (locally in $X$),
tensor products of
differential operators and enveloping algebras,
so Frobenius structures come from the canonical Frobenius structures
for
differential operators from Lemma \ref{Frobenius property}.c,
and enveloping algebras
from the claim (c). All of the $p$-central varieties involved here
are smooth and Calabi-Yau.
The algebras
$\tii \DD_\PP$ are in this class and have the additional property that
their $p$-central variety $\fg^*_\PP\tw$ are  proper
over $\fg^*\tw$, so we can use
Lemma \ref{Frobenius property}.b.

(c)
We only check the Frobenius claim for $\U^\la$, the rest then follows as
in the proof of (b).
Since $\U$ is Frobenius over $\fg^*\tw$ by the part (c),
the restriction
$$
\U\ten_{\OO(\fg^*\tw}\
[\OO(\fg^*\tw)\ten_{\OO(\fh^*\tw/W)}\k_{{\AS}(\la)}] =\
\U\ten_{\OO(\fh^*/W)}\
[\OO(\fh^*/W)\ten_{\OO(\fh^*\tw/W)}\k_{{\AS}(\la)}]
$$
is Frobenius over
$\fg^*\tw\tim_{\fh^*\tw/W}{\AS}(\la)$.
This restriction carries the action of
the algebra of functions
on the fiber of  $\fh^*/W@>>>\fh^*\tw/W$ at ${\AS}(\la)$.
So each connected component
of the fiber of $\fh^*/W@>>>\fh^*\tw/W$ at ${\AS}(\la)$
gives a summand of the restriction,
and these summands are again
Frobenius.
When the map is unramified at $\la$ then
$\U^\la$ is one of these summands.

\sus{
Proof of Theorem \ref{Localization theorem}
}\label{19}
The theorem has already been reduced to the part
(a) (see \ref{reduction to a}).
Since $\Db[\mod^c(\tii \DD_\PP)]$ is  Calabi-Yau with respect to
$\O(\fg^*\tw)$
(corollary \ref{U is Calabi-Yau}),
 its full triangulated
subcategory
$
\Db[\mod^c_\la(\tii\DD_\PP)]
$ is also Calabi-Yau with respect to
$\O(\fg^*\tw)$.
Therefore, the equivalence \eqref{equivD} follows from
Proposition \ref{compid}.b and
Lemmas \ref{abstract}, \ref{indecompo}.

To deduce \eqref{equivD0} from \eqref{equivD} we use Lemma
\ref{adjoints and forgettings}(d).
It says that the functors $i$, $\iota$ send the
adjunction arrows into adjunction arrows; since $i$, $\iota$ kill
no objects, and the adjunction arrows in
$\Db[\mod^c_\la(\tii\DD)]$, $\Db[\mod_\la^{fg} (\U)]$ are
isomorphisms, we conclude that the adjunction arrows in
$\Db[\mod^c(\DD^\la)]$, $\Db[\mod^{fg} (\U^\la)]$ are isomorphisms,
which implies \eqref{equivD0}.
\epf

\sus{
Twisted $\DD$-modules on parabolic flags.
} \label{DPsec}
%{$\PP$ is derived $\DD^\lambda$-affine.} 

Fix a parabolic $P=LJ$ with a partial flag variety
$\PP$.
In this subsection 
 we switch  attention from  $\DD_\PP^\la,\ \la\in\fh^*$, 
to  the more traditional ``smaller'' sheaf
 of rings,  the sheaf of twisted  differential operators on $\PP$.
Possible twists are indexed by  $\la\in \Pic(\PP)\otimes \k$;
recall  that $\Pic(\PP)$ can be identified with 
 the sublattice $\La_\PP\subset \La$
consisting of the elements  $\la$ such that
$\langle\la,\alpha\rangle =0$ if $\alpha$ is a coroot of the Levi.\fttt{
Notice that  $\La_\PP$ contains some non-singular weights, due
to the $\rho$ shift in the definition of regularity.
}
For $\la \in \La_\PP \otimes \k$
the corresponding sheaf of rings will be denoted
$
\fD^\la_\PP
$.
Thus $\fD^\la_\PP$ is a central reduction 
of the sheaf of rings $\tii\DD_\PP$ associated with the
$P/[P,P]$ torsor $G/[P,P]$
over $\PP$ as in \ref{Torsors}.
 The sheaf $\fD^\la_\PP$ is related to the sheaf $\DD^\la_\PP$
considered in the other sections via: 
$\fD^\la_\PP=
\DD^\la_\PP
\ten_{\U^\la(\barr\fp)}\ \k$. For example, for $\PP=\BB$ 
we have $\fD_\BB^\la=\DD^\la_\BB$, while for $\PP=\pt$,
$\La_\PP=\{0\}$ and we have: $\DD^0_{\pt}=\U^0$, $\fD_{\pt}^0=\k$.

\sss{$\fD^\la_\PP$ is derived affine
}

The famous result of \cite{BB}
asserts that the flag variety $\BB$ over a characteristic zero field
is $\DD_\BB$-affine, and also $\DD^\la_\BB$-affine for any semi-ample line
bundle $\OO_\la$. This implies that the same is true for
the partial flag variety $\PP$: the global sections of a $\fD^\la_\PP$-module
coincide with the global sections of its pull-back to $\BB$;
since the pull-back functor is exact and faithful, it is clear that the functor
of global sections on the category of $\fD^\la_\PP$-modules is exact and 
faithful,
provided that this is known for $\DD^\la_\BB$ 
modules.
However, for a sheaf of rings $\RR$ over an algebraic variety
$X$
exactness and faithfulness of the global sections
functor for $\RR$-modules is equivalent to $X$ being $\RR$-affine. 
%UP TO SOME FINITENESS CONDITION?

 In this subsection 
we show that a similar fact holds on the level of
derived categories in positive characteristic. % $p> h$.
The material of this section is not used elsewhere in the
paper.
For simplicity we treat the case of an integral regular weight only,
 leaving
the general case as an exercise to the interested reader.
Fix $\la\in \La_\PP$, we will use the same notation for 
$\la\ot 1\in \La_\PP\ot \k$. 

We set $\U^\la_\PP=\Gamma(\fD^\la_\PP)$, and let $\pi$
denote the projection $\BB\to \PP$. 
We assume that $p>h$ till the end of the section.  

\medskip

\spro \label{DPequiv} {\em  Let $\la\in \La_\PP\otimes \Fp$ be regular, and assume
that $\Rr^i\Gamma (\fD^\la_\PP)=0$ for $i>0$.
 Then
$\fD_\PP^\la$ is derived affine, i.e. we have an equivalence
$\RGa:\Db[\mod^c(\fD_\PP^\la)]\iso \Db[\mod^{fg}(\U_\PP^\la)]$.}

 {\em Sketch of proof.}\  To show that a coherent  sheaf $\RR$ 
of $\OO_X$-algebras 
on a Noetherian scheme $X$ is ``unbounded derived affine'',
i.e., that the derived global sections functor $\RGa$ induces an
equivalence between $\Dmin[\mod(\RR)]$ and $\Dmin[\mod(\Ga(\RR))]$  it
suffices to check that the composition of the global sections functor and the left adjoint 
localization functor is isomorphic to identity, and 
that the global sections functor
$\RGa$ does not kill any complex of $\RR$-modules. 
To prove
a similar statement   for the bounded
derived category it is enough to verify also that
$\RGa(\F)\in \Db[{\mathrm{Vect}}] \iff \F\in \Db[\mod(\RR)]$ for $\F\in
\Dmin[\mod(\RR)]$.

%Cohomology vanishing for $\OO_{T^*\PP}$ is well-known \cite{LK}.
%\cite[6.18]{Ja0}. 
%By a standard argument \cite[3.4.1]{BMR} it implies
%cohomology vanishing for  
%$\fD_\PP^\la$.
Vanishing of $\Rr^i\Gamma(\fD^\la_\PP)$ for $i>0$ implies that
 $\RR=\fD_\PP^\la$
satisfies the first of the above conditions. The remaining ones
follow directly from the established derived affinity of
$\DD_\BB^\la$, together with the fact
that % for $\pi:\BB\to \PP$ 
the functor
$\pi^*:\mod^c(\fD_\PP^\la) \to \mod^c(\DD_\BB^\la)$ is exact and
commutes with derived global sections. The last claim  follows by the
projection formula  from $\Rr\pi_*(\OO_\BB)=\OO_\PP$. \epf

\rem The higher cohomology vanishing for $\fD^\la_\PP$  follows 
from the similar vanishing for $\OO_{T^*\PP}$ 
by a standard argument \cite[3.4.1]{BMR}. 
The latter vanishing holds in characteristic zero
by the Grauert-Riemenschneider vanishing theorem \cite{Broer}.
It is expected to hold in characteristic $p$ \cite[5.C]{BrKu}.
For a fixed dominant $\la$ and large $p$ the higher cohomology vanishing
for $\fD^\la_\PP$ follows from Lemma \ref{whenonto} below.

\sss{
Equivalence class of the Azumaya algebra
$\fD_\PP^\la$ 
}
Set $X= T^*\PP\tw\times _{\PP\tw} \BB\tw$, and 
consider the maps $T^*\PP\tw @<{\PR}<< X
%T^*\PP\tw\times _{\PP\tw} \BB\tw
\aa{\iota}\inj
T^*\BB\tw$. 
We also have the moment maps $\mu_\PP:T^*\PP\tw\to \N\tw$,
$\mu_\BB:T^*\BB\tw\to \N\tw$. 

Recall that the restriction
$\U^{-\rho}$ of  $\U$  
under the embedding $\N\tw\imbed \ZZ(\U)=\g^*\tw\times_{\h^*\tw/W} \h^*/W$,
$\chi\mapsto (\chi,-\rho)$ is an 
Azumaya algebra. An isomorphism $\DD^{-\rho}_\BB
=\DD^{(p-1)\rho}_\BB\cong \mu_\BB^*\U^{-\rho}$
of Azumaya algebras on $T^*\BB\tw$ has been explained in \cite[5.2]{BMR}.

%Let $\pi:\BB\to \PP$ be the projection. 
The sheaf 
$\fD_{\BB\to \PP}\dff
\pi^*(\fD_\PP)$ is a left module for $\DD_{\BB}$ and a right module for
$\pi^{-1}\fD_\PP$. It is not hard to check that as 
a module over the center $\ZZ(\DD_\BB)=\OO_{T^*\BB\tw}$ the sheaf
$\DD_{\BB\to \PP}$ is supported
on  the image of $\iota$. Moreover, the action of the center 
$\pi^{-1}\OO_{T^*\PP\tw}\subset \pi^{-1}\DD_\PP$ is compatible with the
action of $\OO_{T^*\BB\tw}$ via $\PR$, thus $\DD_{\BB\to\PP}$ is a module
over $\iota^*(\DD_\BB)\otimes _{\OO_X}
\PR^*(\fD_\PP^{\mathrm op})$. 
According to \cite[Proposition 3.7]{BurBur} or
\cite[Theorem 2.4]{OV}, this bimodule provides an equivalence
between the two Azumaya algebras 
$\iota^*(\DD_\BB)$ and 
$\PR^*(\fD_\PP)$.
Hence, for line bundles $\OO_{\BB,\nu}$ on $\BB$ and
$\OO_{\PP,\la}$ on $\PP$ the sheaf 
$$
^\nu\fD_{\BB\to \PP}^\la
\dff
\O_{\BB,\nu}\otimes _{\OO_\BB} 
\fD_{\BB\to \PP} \otimes 
_{\pi^{-1}\OO_\PP} \pi^{-1}\OO_{\PP,-\la}
$$ 
is a left module
for $\iota^*\DD^\nu_\BB$ and a right module for $\PR^*\fD_\PP^\la$,
providing an equivalence between the two Azumaya algebras.
In particular, for $\nu=(p-1)\rho$ we get an equivalence 
between $\DD^{(p-1)\rho}_\BB\cong \mu_\BB^*(\U^{-\rho})$ and
$\PR^*(\fD^\la_\PP)$. 

\slem
\label{splDP} 
{\em 
For every $\la\in\La_\PP$
there exists
a unique (up to a unique isomorphism) module $\MM$ for
$ \mu_\PP^*(\U^{-\rho})
\otimes_{\OO_{T^*\PP\tw}}
(\fD^\la_\PP)^{\operatorname {op}}
$
together with an isomorphism of bimodules
$\PR^*(\MM)\cong {\ }
 ^{(p-1)\rho}\fD_{\BB\to \PP}^{\la}$,\ 
where we have used the identification
$\PR^*\mu_\PP^*\U^{-\rho}=\mu_\BB^*\U^{-\rho}\cong \DD_\BB^{(p-1)\rho}$.
This module provides an equivalence
between Azumaya algebras $\mu_\PP^*(\U^{-\rho})$ and $\fD^\la_\PP$.
}

\pf First we show that there is some equivalence
between $\mu_\PP^*(\U^{-\rho})$ and $\fD_\PP^\la$.
It follows from  \cite[IV.2, Corollary 2.6]{MI}
 that such an equivalence exists,
provided that it exists on an open dense subvariety of $T^*\PP\tw$.
However, over an open dense subvariety the map $\PR$ admits a section,
thus the sought for equivalence comes from the bimodule 
$ ^{(p-1)\rho}\fD_{\BB\to \PP}^{\la}$.

Let now $\MM'$ be some bimodule providing an equivalence 
between $\mu_\PP^*(\U^{-\rho})$ and 
$\fD_\PP^\la$.
Then 
$\PR^*\MM'\cong \LL \otimes _{\OO_{X}}   {
^{ (p-1)\rho}\fD^\la_{\BB\to\PP} } $ 
for some line bundle
$\LL$ on $X$, because
a splitting bundle for an Azumaya algebra is defined uniquely
up to twist by a line bundle. Existence of $\MM$ would follow once
we show that $\LL\cong \PR^*\LL'$ for some line bundle 
$\LL'$ on $T^*\PP\tw$.
Since the Picard group of the total space of a vector bundle
is isomorphic to the Picard group of the base, it is sufficient to see
that the restriction of $\LL$ to the zero section $\BB\tw\subset
X$
is a pull-back under $\pi\tw$. Furthermore,
since the quotient $\La/\La_\PP$ is torsion free, it is enough to check that 
the vector bundle 
$
\OO_{\BB\tw}
\otimes_{\OO_X}\
^{(p-1)\rho}\fD_{\BB\to\PP}^\la
$ 
on $\BB\tw$
is a pull-back under $\pi\tw$. The latter bundle is readily
identified with $\HHom \left( (\pi\tw)^*
\Fr_*(\OO_{\PP,\la}),
\Fr_*(\OO_{\BB,(p-1)\rho}) \right)$, cf. \cite[2.2.5]{BMR}.
Since the vector bundle
$\Fr_*(\OO_{\BB,(p-1)\rho})$ on $\BB\tw$ is trivial, we have proved
the existence of $\MM$. Its uniqueness up to a unique
isomorphism follows from the fact that 
 $\PR ^*$ is fully faithful
on the abelian category of quasi-coherent sheaves, which is clear from
$\pi_*\OO_\BB=\OO_\PP$. \epf

\cor
\label{eqcohP} 
{\em Let $\la\in \La_\PP$ 
 satisfy the assumptions of Proposition \ref{DPequiv}, and
 $\chi\in \N\tw$ be in the image of $\mu_\PP$. Then 
we have a canonical equivalence of triangulated categories
$$
\Db [\mod^{fg}_\chi(\U_\PP^\la)]
\
\cong 
\
\Db [\Coh _\chi(T^*\PP\tw)]
,$$
where the subindex $\chi$ denotes, respectively, 
the full subcategory of modules with the action of the central
subalgebra $\Gamma (\OO_{T^*\PP\tw})$ having generalized central character
defined by $\chi$
and
the full subcategory
of sheaves set-theoretically supported on $\mu_\PP^{-1}(\chi)$.
} \epf

We denote the inverse of this equivalence by $\F\mapsto M_\F^\la$.

\rem  It follows from the construction that for $\chi=0$
the splitting bundle for 
$\fD^\la_\PP$ on the formal neighborhood of $\mu_\PP^{-1}(0)$, that we
used to construct the equivalence, 
restricts 
on the zero section $\PP\tw\subset T^*\PP\tw$
to the bundle $\Fr_*(\OO_{\PP,\la})$.

It is explained  in Remark \ref{Spring is there}
 above that for $\PP=\BB$ an equivalence
as in the last Corollary can be fixed by fixing a weight $\nu$ 
such that $\nu+\la$
is unramified. Comparing the definitions we see that the equivalence
of Corollary corresponds to the choice $\nu=(p-1)\rho-\lambda$.

\ex \label{SLn} Let $G=SL(n+1)$, $\PP=\Pn$. It is easy to see that
$\Rr^i\Gamma (\OO_{T^*\Pn})=0$ for $i>0$, thus Proposition \ref{DPequiv}
applies for line bundles corresponding to weights, which are regular
modulo $p$.

A line bundle $\OO(i)$
on $\Pn$ corresponds to the weight $i\omega$, where $\omega$
is a fundamental weight of $SL(n+1)$. The weight $i\omega$
modulo $p$ is regular if and only if  $i\ne j\ \mod \ p$
for $-n\leq j \leq -1$. Without loss of generality assume that
$0\leq i< p-n$. 

Fix $\chi=0$. For $\F\in \Db \Coh _{\mu_\PP^{-1}(0)}(T^*\Pn\tw)$
we have %by definition
 $M_\F= R\Gamma(\E\otimes \F)$, where $\E$ is the splitting bundle
for $\DD_\Pn^{i\omega}$ on the formal neighborhood
of the zero section in $T^*\Pn\tw$. 
According to the previous remark, $\E|_{\Pn\tw}
\cong \Fr_*(\OO_\Pn(i))$.
%By \cite[III.6.4]{Ha},  $\Fr_*(\OO(i))$ is a sum of line bundles. 
%It is also not hard to show that
By \cite[Proposition 4.1]{HKR}, 
$\Fr_*(\OO(i))\cong \oplusl_{j=0}^n\OO_\Pn(-j)^
{\oplus d_j}$ for some $d_j>0$. We claim that $\E\cong
\oplusl_{j=0}^n \hatt\OO(-j)^{\oplus d_j}$,
where $\hatt \OO$ stands for the structure sheaf of the
formal neighborhood of the zero section in $T^*\Pn\tw$.
To show this one can use induction in $m$ to construct an isomorphism
on the $m$-th neighborhood of the zero section. The obstruction
on the $m$-th step lies in 
$\Hh^1\left( \Pn,\operatorname{Sym}
^m(\TT_{\Pn})\otimes
\EEnd(\Fr_*(\OO(i)) ) \right)$. It is not hard to show that 
$\Hh^{>0}(\OO_{T^*\Pn}(j))=0$
for $-n\leq j\leq n$, thus the obstruction vanishes.

It follows that the coherent sheaves corresponding to indecomposable
projective pro-objects in 
$\mod_0(\U^{i\omega}_\Pn)$
are the line bundles $P_j=\hatt\OO(j)$, $j=0,\dots,
n$. The objects $L_k$ in the derived category of coherent sheaves 
corresponding to irreducible $\U^{i\omega}_{\Pn}$ modules are characterized by
$\RHom(P_j, L_k)=\k^{\delta_{jk}}$. 
Thus we have $L_j= {\mathfrak i}_*\Omega_{\Pn}^j(j)[j]$,
where $\mathfrak i$ stands for the embedding of the zero section.

\srem
 The same objects $\Omega_{\Pn}^j(j)[j]$ correspond
to irreducible objects in $\Coh^{\Zet_{n+1}}({\mathbb A}^{n+1})$
under an equivalence between $\Db[\Coh ^{\Zet_{n+1}}({\mathbb A}^{n+1})]
$ and the derived category of coherent sheaves on the total
space of canonical bundle on $\Pn$, cf., e.g.,  \cite[4.3]{Reid}.

\sss{Compatibility} 
For $\la\in \La_\PP$ we have the following exact functors between
abelian categories:
$$
\Coh_{\mu_\PP^{-1}(\chi)}(T^*\PP\tw)\cong
\mod_\chi(\fD_\PP^\la) 
@>{\pi^*}>> \mod_\chi(\DD_\BB^\la)
\cong 
\Coh_{\BB_\chi}( T^*\BB\tw)
,$$
where the equivalences come from the canonical splitting arising from
Lemma \ref{splDP}. The definition of equivalences
implies that the composition functor is identified with
$\iota_*\PR^*:\Coh_{\mu_\PP^{-1}(\chi)}(T^*\PP\tw)\to 
 \Coh_{\BB_\chi}(T^*\BB\tw)$.

Since $\Rr \Gamma \circ \pi^* (\MM)\cong \Rr\Gamma(\MM)$ for $\M\in 
\DD^\la_\PP$ we arrive at the following

\spro {\em Fix $\PP$, $\la\in \La_\PP$ satisfying the assumptions of Proposition 
\ref{DPequiv}.
Let $\phi^\la_\PP:\U^\la \to \U^\la_\PP$ be the map coming
from the action of $\g$ on $\PP$. Then for 
 $\F\in \Db[
\Coh_{\mu_\PP^{-1}(\chi)}(T^*\PP\tw)]$ we have a canonical isomorphism
$(\phi^\la_\PP)^*(M_\F^\la)\cong M^\la_{\iota_*\PR^*\F}$.}

\rem \label{pushpull} Assume that the map $\phi^\la_\PP$ is surjective.
Then for an irreducible $\U^\la_\PP$-module $M$, the module
$(\phi^\la_\PP)^*(M)$ is also irreducible.  Thus in this case the previous
proposition  can be used to describe explicitly the image
of an irreducible $\U^\la$-module under the equivalence
of \cite{BMR}, once the image
of an irreducible $\U^\la_\PP$ module in $\Db[\Coh(T^*\PP\tw)]$
is known.
For example, if $G=SL(n+1)$, $\PP=\Pn$,
then Lemma \ref{whenonto}.a below implies surjectivity, thus 
Example \ref{SLn} yields an explicit
description of $n+1$ out of the total $(n+1)!$ irreducible objects.

We finish the section by giving a sufficient criterion for surjectivity
of $\phi^\la_\PP$.

\lem \label{whenonto}{\em
a) Suppose that the map $\mu_\PP$ is birational onto its image,
and the image is a normal variety. Then the map $\phi_\PP^\la$
is surjective for any $\la\in \La_\PP$.

b) Fix  $\la\in \La_\PP$, such that $\lambda+\rho$ is dominant.
 There exists
(an explicitly computable)
$N\in \Zet$ depending on the type of $G$ and on $\la$, such that
$\Rr^i\Gamma(\fD^\la_\PP)=0$ for $i>0$ and
$\phi^\la_\PP$ is surjective, provided that $p=\chara(\k)>N$.}

{\em Sketch of proof.}\ The assumptions in (a) imply that the map
$\mu_\PP^*:\operatorname{Sym}(\g)\to \Gamma (\OO_{T^*\PP})$ 
is surjective. This yields surjectivity of $\phi^\la_\PP$ by 
a well-known argument, cf., e.g.,  \cite[3.4.1]{BMR}.

To show (b) observe that 
$Rr^i\Gamma (\fD^\la_\PP)=Rr^i\Gamma (\pi^*\fD^\la_\PP)$, while
the map $\phi_\PP^\la$ 
can be obtained from the canonical surjection of sheaves
$\DD^\la_\BB\to \pi^*(\fD^\la_\PP)$ by applying the functor
of  global sections. This map of sheaves can be extended
to a resolution of $\pi^*(\fD^\la_\PP)$
by locally projective $\DD^\la_\BB$ modules, coming
from the relative De Rham complex of the map $\pi$.
The terms of this resolution are of the form 
$\DD^\la_\BB\otimes
_{\OO_\BB}\Lambda^i(\TT_{\BB\to \PP}) $ placed in degree $-i$,
where $\TT_{\BB\to \PP} $ is the kernel of the differential
$d\pi$. Thus desired surjectivity and cohomology vanishing 
 follow from cohomology vanishing
$\Hh^{>0}(\DD^\la_\BB\otimes
_{\OO_\BB}\Lambda^i(\TT_{\BB\to \PP}))=0$ for all $i$,
which can be established by the method of \cite{BB}.
Namely, we
choose a dominant $\nu\in \Lambda$, such that 
$\Hh^{>0}
(\OO_{T^*\BB}\otimes _{\OO_\BB}\Lambda^i(\TT_{\BB\to \PP})
\otimes _{\OO_\BB}\OO_{\BB,\nu})=0$, it suffices that $\nu+\alpha$
is dominant when $\alpha$ is a sum of some subset of positive
roots in the Levi $\bar P$. Then 
$\Hh^{>0}[\O_{\BB,\nu}
\otimes _{\OO_{\BB}}\DD^\la_\BB\otimes
_{\OO_\BB}\Lambda^i(\TT_{\BB\to \PP})]=0$. Thus we will be done
if we show that for a $\DD^\la_\BB$-module $\MM$, the sheaf $\MM$ is a direct
summand in the sheaf $V\otimes _\k \O_{\BB,\nu} \otimes_{\OO_\BB} \MM$,
where $V$ is the irreducible representation of $G$ with lowest weight 
$-\nu$. We have an action of $\g$ on this sheaf, preserving the filtration
arising from the canonical filtration of $V\otimes \O_\BB$
by line bundles. The associated graded sheaves of the filtration
are of the form $ \O_{\BB,\eta+\nu} \otimes_{\OO_\BB} \MM$, $\eta\in 
\WT(V)$, where $\WT(V)$ is the set of weights of $V$.
 The Harish-Chandra center acts
on such a subquotient via the character corresponding to the
$W'\aff$ orbit of $\eta+\nu+\la$. Thus if 
$\la$ is the only element in 
$(\la+\nu +\WT(V))\cap W'\aff\bu \la$,
%of this form in the orbit $W_{aff}'\bu \la$, 
then $\MM$ is indeed
a direct summand in $V\otimes _\k \O_{\BB,\nu} \otimes_{\OO_\BB} \MM$.
Since $\la+\rho$ is dominant, we have $w\bu(\la)\preceq \la$
for any $w\in W$, while $\eta\succeq -\nu$ for any $\eta\in \WT(V)$.
Thus $(\la+\nu +\WT(V))\cap W\bu \la=\{\la\}$. There is only
a finite number of primes dividing one of the finite
number of nonzero weights $\eta+\nu+\la
-w\bu\la$, $w\in W$, $\eta\in \WT(V)$, $\eta\ne -\nu$. 
For $p$ outside of this finite set we have
$\{\la\}=(\la+\nu +\WT(V))\cap W'\aff\bu \la$. \epf

\se{\bf
Affine braid group action on $\fg$-modules by intertwining
functors 
}

Here we develop a characteristic $p$ version of the theory of
intertwining functors of
 Beilinson and Bernstein
\cite{BeBe}, see also \cite{BeGi}, \cite{Milicic}
 for results in characteristic
zero which are parallel to many of the results in this section.

For each regular integral weight $\la$ we construct a canonical
action of an incarnation $^\Th {B'}\aff$ of the (extended) affine
braid group $B'\aff$ on the category  $\Db[\mod_{\la}^{fg}(\U)]$.
The group $^\Th B'\aff$  depends  on a regular orbit
$\Th=\Wap\bu \la$;
%but a choice
%of an isomorphism $^\Th B'\aff\cong B'\aff$ depends on $\la$.
it is non canonically isomorphic to  $B'\aff$.
The
action is generated by {\em intertwining functors}; they intertwine
different localization functors.

The localization equivalences allow us to
 translate the  action of\
$^\Th B'\aff$ from $\Db[\mod_{\la}^{fg}(\U)]$ to $\Db[
\mod_\la^c(\tii\DD)]$. Here  the group that naturally acts on
$\Db[\mod_\la^c(\tii \DD)]$ is always $B'\aff$. Moreover, the
categories $\mod_\la^c(\tii \DD)$ for different weights $\la$ are
canonically identified by tensoring with line bundles; this
identification is compatible with the  $B'\aff$ action (Corollary
\ref{bbbI_acts}).

In the remainder we assume that $p>h$; recall that regular weights
exist only for $p\geq h$, the case $p=h$ is excluded to avoid 
$G=SL(p)$ which violates assumption (C) from \ref{Restrictions on G and p}.
 Since exponents of $W$ do not exceed $h$, this also implies that $p$ does not divide $\mid W\mid$.

%{\bf{Disclaimer.}} 
\setcounter{subsubsection}{0}
\sss{
Further extensions
}
In this section we do not work in the maximal possible
and perhaps even not in the maximal reasonable generality. For example,
translation functors could be defined and probably described
geometrically for modules
with non-integral Harish-Chandra central characters. This would, however, 
further complicate the notations without yielding new ideas or 
obvious
%potential
applications, 
%in sight, 
thus we did not 
pursue 
%include 
it.
 
%Also, 
In this paper the action of the braid group is introduced on derived categories
of modules with a generalized Harish-Chandra central character.
One can, in fact, define it also for the derived categories of modules
with a fixed  Harish-Chandra central character. This technical variation
will be discussed in a future publication.

\sus{Affine braid group actions: statement}

We start with recalling some standard material on
 affine Weyl groups.

 \sss{Affine Weyl and
braid groups }

The groups $W$, $\Waff$ and some related objects were introduced in
\ref{dotac}. These groups are Coxeter groups (with the sets of simple
reflections $\Sigma$, $\Iaff$ respectively), so we can consider the
corresponding braid groups $B$, $B\aff$.

These algebraic constructions have a topological interpretation:\
$B=\pi_1((\fh^*_{\Ce})^{\mathrm{reg}})/W)$, $B\aff=\pi_1 ((\ch
H_{sc})^{\mathrm{reg}}/W)$ where $\ch H_{sc}=\ch H_{sc} (\Ce)$ is the torus
dual to the Cartan subgroup $H_{adj}$  of the adjoint group
$G_{adj}=G/Z(G)$.
% (i.e. the lattices $X_*(H)=\Hom(\Gm,H)$ and
% $X_*(\ch H)=\Hom(\Gm,H^*)$
% are dual).

 For $\alpha\in \Sigma\aff$
let $s_\alpha\in W\aff$, $\st_\alpha\in B\aff$ denote the
corresponding standard generator.
%; the braid monoid
%$B^+\aff}\subset B\aff$ is
% the subsemigroup (with unit) generated by $\st_\alpha$,
%  $\alpha\in \Sigma\aff$ .

We have an embedding of sets $C:W\aff\inj B\aff$, also denoted
by $w\mapsto \tii w$, where for a minimal length decomposition
$w=s_{\alpha_1} \cdots s_{\alpha_{\ell(w)}}$ we have $C(w)=\tii
w=\st_{\alpha_1} \cdots \st_{\alpha_{\ell (w)}}$. The map $C$ is not a
group homomorphism; however, we have $C(w_1 \cdot w_2)=C(w_1)\cdot
C(w_2)$ when the lengths add up, \ie $\ell(w_1w_2)=\ell(w_1)+\ell(w_2)$,
where $\ell:\Waff \to \Z_{\geq 0}$ is the length function sending
$w$ to the length of its minimal decomposition.
 This happens, in particular, when
 $w_1=\lambda_1$, $w_2=\lambda_2$ lie in the semigroup
$ Q\dom$ of dominant weights in $Q$ (we have $Q\dom\subset Q \subset
W\aff$). The image of $C$ generates $\Baff$, and it is known that
$\Baff$ can be defined by generators $\tii w,\ w\in \Waff$, and
relations $\tii u\tii v=\tii{uv}$ when lengths add up:
 \begin{equation}\label{lplus}
 \widetilde{uv}=\tilde u \cdot \tilde v
 \ \ \ \ \ \ {\rm{when}}\ \
\ell(u)+\ell(v)=\ell(uv). \end{equation}

The length function admits a natural extension to the extended
affine Weyl group $\Wap$, which is characterized by
$\ell(w\om)=\ell(w)$ if $w\in \Waff$, and $\om \in \Om\dff
{\mathrm{Stab}}_{\Wap}(A_0)$.
 The extended affine braid group $B'\aff$ is
defined
 by generators
$\tilde w,\ w\in \Wap$, and
 relations \eqref{lplus}.

 Then one finds that
$B'\aff$ is the semidirect product $\Baff\rtim\Om$, where the finite
abelian group $\Om ={\mathrm{Stab}}_{\Wap}(A_0)\cong \La /Q =\Wap/W\aff$
acts on $\Baff$ via its action  on the set of simple reflection
$\Iaff$.

\sss{Local action and the affine braid group} For each $\al\in
\Iaff$, each alcove $A$ has a unique face of type $\al$. We denote
by $*$ the right action of $\Waff$ on the set of alcoves such that
for $\al\in \Iaff$ the alcove $A*s_\al$ is the reflection $s_F(A)$
of $A$ in the unique face $F$ of $A$ of type $\al$ (cf., e.g.,
\cite[1.1]{LuAdv}).
We also get a free $*$-action on the set of regular characters
$\La_{\mathrm{reg}}$\ie the ones which lie in alcoves, the quotient is
identified with the set of characters in the open fundamental
alcove.

We will need to twist
these constructions by certain
torsors $\Th$ for  the affine Weyl group as
follows.

 Let $\Theta\subset
\Lambda_\RE = \Lambda \otimes \RE$
 be a free orbit of $\Wap$.
We define the {\em local } extended affine Weyl group $\LW'\aff$
 as the group of permutations of the set
$\Theta$ which commute with the action of $\Wap$. We
endow $\LW'\aff$ with the composition law opposite to the
composition of permutations; thus the action of $\LW'\aff$ on $\Theta$
is a right action, we denote it by $w:\la\mapsto \la \ast w$.
Every choice of $\lambda\in \Theta$ defines an isomorphism
$\isla:\Wap\cong \LW'\aff$. The image of the normal
subgroup $W\aff\subset \Wap$ under $\isla$
 does not depend on $\lambda\in \Theta$, thus we get a normal subgroup
$\LW\aff=\isla(W\aff)\subset \LW'\aff$.

We use the isomorphism $\isla$ for $\la$ in the fundamental
alcove $A_0$ to endow $\LW\aff$ with a structure of a Coxeter
group; an element $s\in \LW\aff$ is a simple reflection if and only if it
sends every $\lambda\in \Theta$ to an element which is symmetric to
$\lambda$ with respect to a  face of the  alcove containing
$\lambda$. The set of simple reflections $\LI\aff\subset
\LW\aff$ is in a canonical bijection with the set of pairs
$(\lambda, F)$ where $\lambda \in \Theta$, and $F$ is a  face of the
alcove of $\lambda$, modulo the action of $\Wap$. A choice of
$\lambda\in \Theta$ defines a bijection $\LI\aff\cong
\Sigma\aff$; in fact,  this bijection depends only on the
$W\aff$ orbit of $\la$, and one easily sees that bijections
corresponding to $\la, \mu \in \Theta$ lying in different $W\aff$
orbits differ by composition with an element
 of the abelian group
$\Omega \cong \Wap/W\aff$ acting on $\Iaff$.

One can define the length function $\ell:\LW'\aff\to \Zet_{\geq
0}$ in such a way that $\ell(\isla(w))=\ell(w)$ for any
$\lambda\in A_0\cap \Theta$. It can be described in a more intrinsic
way as follows: for any $\lambda\in \Theta$ and $w\in \LW'\aff$
the number of affine coroot hyperplanes $H_{\ch \al, n}$ separating
$\lambda$ and $w(\lambda)$ equals $\ell(w)$. We can then define the
braid groups $\LB\aff$, $\LB'\aff$ by generators $\tii w$, $w\in
\LW\aff$ and relations \eqref{lplus}. Thus $\LB'\aff$  comes with a
canonical lifting map $\LC:\LW'\aff\to \LB'\aff$, $w\mapsto
\tilde w$.

\srem
 Notice that, although our description involved the choice of
an initial alcove $A_0$, the  group $\LW\aff$ with its Coxeter
structure, and the groups
 $\LW'\aff$, $\LB\aff$, $\LB'\aff$
can be constructed in a canonical way from the affine Euclidean space
$\Lambda \otimes \RE$ subdivided into alcoves with a fixed subset
$\Theta \subset \Lambda \otimes \RE$ (unlike the Coxeter structure
on the group $W\aff$, which depends on the choice of an alcove).

\sss{} We proceed to state the main results of this chapter.

%For $w\in \Wap$ and $\la\in \La$ we will say that  $w$ {\em
%increases} $\lambda$ if $w=\sigma s_1\cdots s_n \om$, where $\om\bu
%\la$ and  $\la$ lie in the same alcove, $s_i\in \LI\aff$, and
%$s_{i+1}\dots s_n\omega\bu \la
% < s_i s_{i+1}\dots s_n\omega\bu \la $ for every
%$i=1,\dots, n$.

%Similarly,
For $w\in \LW'\aff$ and $\lambda \in \Theta$ we will
say that $w$
{\em increases} $\lambda$ if $w=s_1\cdots s_n \om$
where $\ell(\om)=0$,  $s_i\in \LI\aff$, and $\la \ast s_1 \cdots
s_i < \la *s_1 \cdots s_{i+1}$ for every $i=1,\dots, n$.

The set of  integral  Harish-Chandra central characters is identified
with
$\fh^*_{\Fp}/W=\Lambda/W'\aff$. We fix an integral regular
Harish-Chandra central character  % $\fz$
and let
$\Theta$ %=\Theta(\fz)
 be the corresponding free orbit of
$\Wap$ in $\Lambda$. We set $\mod_{\Theta}(\U)=\mod_{\la}(\U)$,
$\la\in \Theta$.

For $\lambda \in \Theta$ we have the category of
twisted $D$-modules
 $\mod^c_\la(\tii \DD _\BB)$, the category $\mod^{fg}_{ \la}(\U)$,
the equivalence 
$\RGa_{\tii\DD_\BB, \la}:\Db[\mod^c_\la(\tii \DD_\BB)]\to \Db[\mod^{fg}_{\Theta}(\U)]$ 
and the  inverse equivalence
$\LL^{\hatt\la}$, see 
\cite[Theorem 3.2]{BMR}.

For $\la,\mu\in \La$ we  have an equivalence $\mod^c_\la(\tii \DD
_\BB) \to \mod^c_\mu(\tii \DD _\BB)$ sending $\F $ to
$\OO_{\mu-\la}\otimes _{\O_{\BB}} \F$; we will denote it by
${\mathrm{Tens}}_{\mu,\la}$.

\theo\lab{braid_act_mod} {\em
 There is a  unique action\fttt{Here by an action of a group $H$ on
a category 
 we mean a weak action, i.e., a homomorphism from $H$ to the group
of isomorphism classes of auto-equivalences of the category. 
We expect that the action we construct can be enhanced to a strong action,
where for each $g\in H$  a functor $F_g$ is given and isomorphisms
$F_{gh}\cong F_g\circ F_h$ are fixed and satisfy the natural compatibilities.
We do not attempt to construct this richer structure in the present paper.}
 of $\LB'\aff$
on $\Db[\mod_\Theta (\U)]$, $b\mapsto \bI_b$,
%:\Db[\mod_\Theta (\U)] \to
%\Db[mod_\Theta (\U)],
 such that whenever $w\in \LW'\aff$ increases
some $\lambda\in\Th$
we have
\begin{equation}\label{Tens}
 {\mathrm{Tens}}_{\la,\la\ast w} \circ
\LL^{\hatt{\la \ast w}}
  \cong  \LL^{\hatt{\la}}\circ \bI_{\tilde w} .
\end{equation}
  }

The proof appears in section  \ref{proof_of_act_mod}.

 \rems 1) One can define the action in a slightly different way:
 one could require that \eqref{Tens} holds when $w$ {\em decreases}
 $\la$. Then the canonical generators would act by functors
 $\bI_\al^!$ rather than $\bI_\al^*$ (see
section \ref{Tranintert} below). That choice
 would yield standard formulas for the  action of
 the affine Hecke algebra which have
origin in the action on Iwahori invariant 
functions on the Langlands dual $p$-adic group, or on 
the corresponding categories of constructible sheaves
(see \cite{Lu},
 \cite{CG}). The present normalization
yields more convenient
formulas for the action on the derived categories of coherent sheaves,
cf. \cite{ICM}.

 2)  It follows from Lemma \ref{Up and down}.a
  below that the finite abelian
subgroup $\Omega \subset \LB'\aff$ acts by
{\em translation functors}, which are {\em exact} functors on
the abelian category $\U_{\hat \chi}^{\lambda}$. E.g. for $G=SL(2)$,
$\lambda=0$ and $\chi=0$ one can show that the non identity element in
$\Omega\cong \Zet/2\Zet$ interchanges the trivial and the $(p-1)$
dimensional module.

\cor
\label{bbbI_acts}
{\em
For any regular integral $\la$ there exists an (obviously
unique) action $b\mapsto \bbI_b$ of $B'\aff$ on
$\Db[\mod_c^\la(\tii\DD)]$ satisfying the following requirements:

i) Assume that $\la$ lies in the fundamental alcove $A_0$, so that
the isomorphism $\isla:\Wap\cong \LW'\aff$ respects the Coxeter
structure, and therefore induces an isomorphism
of braid groups ${\bf isom}_\la: B'\aff\cong
 \LB'\aff$.
Then the equivalence
 $\RGa_{\tii\DD_\BB,\la}:\Db[\mod_\la^c(\tii \DD )]\to \Db[\mod_\la^{fg}(\U)]$
intertwines the actions of
 $B'\aff$ and $\LB'\aff$.
%
%   (here we have used ${\bf isom}_\la$
%   to identify  $B'\aff$ with $\LB'\aff$).

ii) If $\la$, $\mu$ are two integral regular weights, then
the equivalence ${\mathrm{Tens}}_{\mu,\la}$ commutes with the $B'\aff$ action. }

See section \ref{action on D} for a proof of the
Corollary and  a brief discussion of the properties of
this action.

\sus{
Translation and intertwining functors for
$\fg$-modules
and $\DD$-modules
} \label{22}

\sss{
Translation and intertwining functors
on $\fg$-modules
}\label{Tranintert}
For $\mu\in\La$
we denote by $M
\mm
[M]_\mu$ the  projection of the
category of finitely generated
$\fg$-modules with a
{\em locally finite action of $\zhc$}
to its direct summand $ \mod^{fg}_\mu(\U)
\dff
\mod^{fg}_{d\mu}(\U)
$ (see section \ref{Derived categories of sheaves supported on a subscheme}
 for notations).
For $\mu,\nu\in \La$ the translation functor
$T_{\nu}^\mu:
\mod^{fg}_{\nu}(\U) @>>> \mod^{fg}_{\mu}(\U) $ is defined by
(see
\ref{Semisimple group G} for the representations
$V_\eta$ and \cite[Section 6]{BMR} for more details)
$$
T_{\nu}^\mu(M)\df\ [V_{\mu-\nu}\ten M]_\mu .$$

It is obvious from the definition that
%\begin{equation}\label{change-by-w}
$$
T_{\nu}^\mu=T_{w\cdot \mu}^{ w\cdot
\nu}
,\ \ \
w\in \Wap
.$$
%\end{equation}
In particular, for $\la\in \Theta$ and  $w\in
\LW'\aff$ the functor $T_{\la}^{\la\ast w}$  is independent of the
choice of $\la\in \Theta$, we denote this functor by $T_w$.

In the
situation
where
%$$
%(\star)
%\ \ \
%\imat
%\tx{ $\mu,\nu\in\La$ lie in the same closed alcove $\AA\sub\ \La_\R$},
%\\
%\tx{
%$\mu$ is regular and $\nu$ lies in a single wall $H$ of $\AA$,
%of type $\al$;
%}
%\eimat
%$$
$\mu$  is regular and $\nu$ lies on a codimension one face of the
alcove of $\mu$ one says that
%$T_{\nu\mu}$
$T_\mu^\nu$ is a {\em down} functor
and $T_{\nu}^\mu$ is an {\em up} functor. The composition $up\ \ci\
down$ is called a {\em reflection} functor or a {\em wall crossing}
functor
$$
R_{\mu|\nu}: \mod^{fg}_\mu(\U) @>>> \mod^{fg}_\mu(\U) ,\ \ \
R_{\mu|\nu}\dff T_{\nu}^\mu\ci T_{\mu}^\nu. $$
These functors are exact, and we use the same notation
for the corresponding functors on the derived categories.
 It is easy to show
that $T_{\mu}^\nu$ and $T_{\nu}^\mu$ are always adjoint (in either
order), hence $R_{\mu|\nu}$ is self adjoint and there are canonical
maps $\id@>>>R_{\mu |\nu}@>>>\id$.
%%% \fttt{
%%% The self-adjointness structure is compatible with
%%% the canonical maps \cite{Ja0}.
%%% }
Define  functors
$\tI_{\mu|\nu}    ^*,\tI_{\mu|\nu    }^!$
as follows:
$$
\tI^!_{\mu|\nu}=\cone(R_{\mu|\nu}@>>>\id)[-1]
\aand
\tI^*_{\mu|\nu} = \cone(\id@>>>R_{\mu|\nu});$$
here $\cone$ denotes the cone of a map between exact functors
on the category of complexes, thus the expressions in the right
hand sides define exact endofunctors of the categories of complexes.
We use the same notations $\tI^!_{\mu|\nu}$, $\tI^*_{\mu|\nu}$ for the induced
endofunctors of the derived category ${\mathrm{D}}(\mod^{fg}_\mu(\U))$.
It is clear that these functors preserve the subcategories $\Dmin$ and $\Db$.

\sss{ Intertwining functors and $\DD$-modules } \lab{Intertwining
functors and D-modules} Let us start with some notation. Recall that for
$\mu\in\La$ we have a canonical isomorphism between
$\DD^\mu\dff ^{\OO_\mu}\DD$ and
 the  specialization $\DD^{d\mu}$ of $\tii\DD\dff\ \tii\DD_\BB $
to the differential $d\mu\in\fh^*$.

%We also denote  $\U^\mu\df \U^{d\mu}$ etc., and use the short form
 %$\Ga_{\eta}$ for $\Ga_{\tii\DD,\eta}$.

For two characters $\mu',\mu''\in\La$
define a functor
$
\tI_{\mu''\mu'}:\
\Dm[\mod^{fg}_{\mu'}(\U)]
@>>>
\Dm[\mod^{fg}_{\mu''}(\U)]
$
by
$$
\tI_{\mu''\mu'}(M)
\dff
\RGa_{\tii\DD_\BB, \mu''}[\LL^{\mu'}M\ten_{\OO_\BB}\OO_{\mu''-\mu'}]
%\ =\
%\RGa_{\mu''}[(\LL^{\mu'}M)(\mu''-\mu')]
.$$
In the case when the differentials
%   $d\mu',d\mu''$
are regular in
$\fh^*$ (i.e.,
%$\mu',\mu''$
characters are $(\Waff,\bu)$-regular),  $\tI_{\mu''\mu'}$ is an equivalence
and one has
$
\tI_{\mu''\mu'}\ci \tI_{\mu'\mu}\cong \tI_{\mu''\mu} $ and
$\tI_{\mu\mu}=\id$. In the
case when $ \Ga(\DD^{\mu'})=\Ga(\DD^{\mu''})$ (this is equivalent to
$d\mu''\in W\bu d\mu'$\ie to $\mu''\in \Wap\bu\mu'$), the
categories $ \Db[\mod^{fg}_{\mu'}(\U)]$ and $
\Db[\mod^{fg}_{\mu''}(\U)] $ coincide and $ \tI_{\mu''\mu'} $ is an
autoequivalence.
%If one thinks of  categories of D-modules
%at $\mu',\mu''$
%as identified
%by tensoring with a line bundle,
%then $
%\tI_{\mu''\mu'}
%$
%is a
%{\em difference of localizations}
%at $\mu'$ and $\mu''$.

\lem {\em
\lab{Up and down}
Let
$\mu,\nu\in\La$ be characters such that
$\nu$ is in the closure of the facet of $\mu$,
and let
$\MM\in {\mathrm{D}}[\mod^c_\mu(\tii\DD)]$
and
$\NN\in {\mathrm{D}}[\mod^c_\nu(\tii\DD)]$.

(a) (``Down''.) We have a functorial isomorphism $
T_{\mu}^\nu(\RGa_{\tii\DD_\BB, \mu}\MM) \ \cong\ 
\RGa_{\tii\DD_\BB, \nu} ( \MM \ten_{\OO_\BB}
\OO_{\nu-\mu} ) $.

%In particular we have:
%$T_{\nu\mu}\cong \RGa_\nu(\OO_{\nu-\mu}\ten\LL_\mu)$.

(b) (``Up''.)
Assume also that
$\mu$ is regular and $\nu$ lies in the (codimension one) face
$H$ of the $\mu$-alcove.
Let $\mud$ denote the reflection of $\mu$ in the wall $H$.
If $\mud<\mu$, then  there is an exact triangle
$$
\RGa_{\tii\DD_\BB, \mud} ( \NN \ten_\OO \OO_{\mud-\nu} ) @>>>
T_{\nu}^\mu(\RGa_{\tii\DD_\BB, \nu}\ \NN) @>>> \RGa_{\tii\DD_\BB, \mu} ( \NN \ten_\OO
\OO_{\mu-\nu} ) .$$ 
(c) 
Keep the assumptions of (b); if ${\mud} <\mu $
then we have natural isomorphisms\
$ \tI^!_{\nu|\mu} \cong \tI_{\mud\mu} $ and $\tI^*_{\nu|\mu} \cong
\tI_{\mu\mud} $. }
%;  if $ {\mud}>\mu $ then  we have natural isomorphisms
%$ \tI^!_{\nu|\mu} \cong
%\tI_{\mu\mud} \aand \tI^*_{\nu|\mu} \cong \tI_{\mud\mu} .$ }

\pf
Claims (a)
and (b) are in \cite[Lemma 6.1.2]{BMR}. For the reader's convenience we recall
 that they follow, respectively,
 from an isomorphism of sheaves
of $\fg$-modules on $\BB$:
 $$T_\mu^\nu(\F)\cong \F\otimes \O_{\nu-\mu}$$
for $\F \in \mod^{c}_\mu(\tii\DD)$; and from the short exact sequence
of sheaves of $\fg$ modules:
$$0\to \F\otimes \O_{\mud - \nu} \to [\F\otimes V_{\mu-\nu}]_\mu
\to \F\otimes \O_{\mu-\nu} \to 0$$
where  $\F \in \mod^{c}_\nu(\tii\DD)$; $\mu$, $\nu$ are as in (b)
and $\mud< \mu$
 (cf. \cite{BMR}, section 6.1.1). Here we have used the same
notation $T_\mu^\nu$, $\F\mapsto \F_\mu$ for the functors
on the categories of sheaves of $\g$-modules as for
the corresponding functors on the categories of $\g$-modules,
which were  introduced earlier. The functors for sheaves are obtained
by applying the corresponding functors for $\g$-modules to the $\g$-module
of sections on each open subset. 

It remains to prove
 (c). Assume first that  $\mud<\mu$; fix
 $M\in  \Db[\mod^{fg}_\mu(\U)]$, and set $\MM=\LL^\mu M$.
We can represent $\MM$ as a bounded complex $C_\MM^\bu$
such that all
$C_\MM^i\in  \mod^{c}_\mu(\tii\DD)$ are $\Gamma$-acyclic where
$\Gamma$ is the functor
of global sections. Part (a) shows that the sheaves
 $C_\MM^i\otimes \OO_{\nu-\mu}$ are also $\Gamma$-acyclic.
Hence the same is true for the sheaf  $C_\MM^i\otimes \OO_{\nu-\mu}
\otimes V_{\mu-\nu}$ and for its direct summand
 $[C_\MM^i\otimes \OO_{\nu-\mu}
\otimes V_{\mu-\nu}]_\mu$. Using the above short exact sequence
we see that the complex of sheaves of $\fg$ modules $C_\MM^\bu\otimes
\O_{\mud-\mu}$ is quasiisomorphic to the following complex
of $\Gamma$-acyclic sheaves of $\fg$-modules:
$$\cone\left( [C_\MM^\bu\otimes \OO_{\nu-\mu}
\otimes V_{\mu-\nu}]_\mu\to C_\MM^\bu\right)[-1].$$
Applying the functor of global sections to this complex we get
the first isomorphism in (c); the second one is
 established in a similar way. \epf

\cor
\lab{properties of translation functors}
{\em
%a)
%Let $\mu,\nu
%\in\La$; if $\nu$ lies in the closure of the $\mu$-facet then for
%$M\in \Db[\mod^{fg}_{\mu}(\U)]$ there is a canonical functorial
%isomorphism $T_{\mu}^\nu(M)\cong
%\RGa_\nu(\OO_{\nu-\mu}\ten\LL^{\hatt\mu}(M))$.
%b)
 Assume that $\mu$ is regular and $\nu$ is on a codimension
one face of the alcove of $\mu$.
Then
$\tI_{\mu|\nu}^*$ and $\tI_{\mu|\nu}^!$ are mutually inverse equivalences.
 \epf}

%\pf
%a)  is immediate from
%lemma \ref{Up and down}.a and the fact that
%$\RGa_\la \circ \LL^{\la}\cong \id$ (even if the composition in the
%other direction is not identity). Part (b) follows from
%Lemma \ref{Up and down}.c. \epf

\sss{ 
Translation  functors for coherent sheaves 
}
\lab{Intertwining  functors for coherent sheaves}
We now explain how to express the translation functors
$T_\mu^\nu$ for some pairs $\mu, \nu\in \La$ in terms of coherent sheaves.
Consider   a  partial flag variety $\PP$ and $\la\in\fh^*$
such that $\la$ is  $\PP$-regular and $\PP$-unramified, so that the
equivalences of Theorem \ref{Localization theorem}.c apply for
any  $\chi$ with $(\chi,\la)\in\ZZ(\U)$.

Recall that the choice of an equivalence in Theorem
 \ref{Localization theorem}.c is determined by the choice of a
splitting vector bundle for the Azumaya algebra $\tii \DD_\PP$
on the formal neighborhood of $\ZZ(\tii \DD_\PP)
\times _{\ZZ(\U)} (\chi, W\bu \la)$ in $\ZZ(\tii \DD_\PP)$.
 Given such a bundle $\MM$
we get an equivalence

$$
\ga_{\MM}:\ \Db[\Coh_{\ZZ(\tii \DD_\PP)
\times _{\ZZ(\U)} (\chi,W\bu \la)}(\ZZ(\tii \DD_\PP))]
 \con\ \Db[\mod^{fg}_{\la,\chi}(\U)] ,\ \ \
\ga_{\MM} (\FF)
\dff
\RGa_{\tii\DD_\PP,\la}[\MM\ten_{\OO_{\ZZ(\tii \DD_\PP)}}\FF] .$$ 

 The splitting bundle
$\MM$ is unique up to tensoring with a line bundle.

We now assume that $\la \in \h^*_\Fp$ is integral.
Then the formal neighborhood of $\ZZ(\tii \DD)
\times _{\ZZ(\U)} (\chi, W\bu \la)$
can be identified with the formal
neighborhood of $\PP\tw_{\chi}=\PP_{\chi,0}\tw$ in $\gt^*\tw_\PP$.

We have $\la =d\eta $
for some $\eta\in \La$; moreover, since $\la$ is $\PP$-regular
and $\PP$-unramified, we have ${\mathrm{Stab}}_W(\la)=W_\PP$, thus we
can (and will) assume that ${\mathrm{Stab}}_W(\eta)=W_\PP$
(cf. Lemma \ref{Applicability of equivalence}.c). 

 Fixing such $\eta$ we get a splitting
bundle $\MM=\MM^\PP_{\chi,\eta}$ described in Remark \ref{Spring is there}.
We let $\ga^\PP_{\chi,\nu}$ denote the corresponding equivalence.

\medskip

In the next lemma we fix two  parabolic subgroups  
$P\subset Q$, let
  $\pi^\PP_\QQ:\PP=G/P\to
\QQ=G/Q$ be the map between the partial flag varieties and 
$\tii \pi^\PP_\QQ$ be the natural map $\gt^*_\PP\to \gt^*_\QQ$. 
% We also let $\pi'$ denote the projection $\BB\to \PP$.

\slem 
\lab{coherent translation of intertwining functors}
{\em
Let
$\mu,\nu\in\La$ be characters such that $\mu$ is $\PP$-regular
and $\PP$-unramified, $\nu$ is $\QQ$-regular and $\QQ$-unramified
and $\nu$ lies in the closure of the facet of $\mu$. 
 
For any nilpotent $\chi \in \g^*\tw$ the  down and up functors
between $\mod_{\mu,\chi}(\U)$ and $\mod_{\nu,\chi}(\U)$ 
 correspond
under the equivalences with the derived categories of coherent sheaves
to the push-forward and pull-back functors under 
$\tii \pi^\PP_\QQ\tw$.
% for the canonical map $\tii
%\pi\tw:\tii\fg^*_\QQ\tw @>>>\tii\fg^*_\PP\tw $.
In other words,  we have 
natural isomorphisms:
$$
T_{\mu}^\nu \ci \ga^\PP_{\chi,\mu} \ \cong\ \ga^\QQ_{\chi,\nu} \ci
\Rr(\tii \pi^\PP_\QQ\ \tw)_* 
\aand 
T_{\nu}^\mu\ci \ga^\QQ_{\chi,\nu} \ \cong\
\ga^\PP_{\chi,\mu} \ci L(\tii \pi^\PP_\QQ\ \tw)^*\cong
\ga^\PP_{\chi,\mu}\ci  L(\tii \pi^\PP_\QQ\ \tw)^!
\ .$$
% or
%$ \mod_\nu^\chi(\U)$.
}

% For simplicity we will
%sometimes omit the Frobenius twist from the notation.

\pf  An isomorphism between 
$(T_{\mu}^\nu\ci \ga^\PP_{\chi,\mu}) 
\FF
 \dff \ T_{\mu}^\nu
[\RGa (\MM^\PP_{\chi,\mu}\ten_{\OO_{\ZZ(\tii \DD_\PP)}}\FF)]
$ 
and  
$
\ga^\QQ_{\chi,\nu} \tii\pi_*\tw\FF
=
\RGa
(\MM^\QQ_{\chi,\nu}\ten_{\OO_{\ZZ(\tii \DD_\QQ)}}
 \tii\pi_*\tw\FF)]
$ 
can be obtained as the following composition:
\begin{multline*}
\ T_{\mu}^\nu[\RGa
(\MM^\PP_{\chi,\mu}\ten_{\OO_{\ZZ(\tii \DD_\PP)}}\FF)] \cong
\ T_{\mu}^\nu[\RGa (\pi^\BB_\PP)^*
(\MM^\PP_{\chi,\mu}\ten_{\OO_{\ZZ(\tii \DD_\PP)}}\FF)]
\\
\cong 
\RGa [ \OO_{\BB,\nu -\mu}\ot _{\OO_\BB} (\pi^\BB_\PP)^*
(\MM^\PP_{\chi,\mu}\ten_{\OO_{\ZZ(\tii \DD_\PP)}}\FF)]
\cong 
\RGa [ (\pi^\BB_\PP)^*(\OO_{\PP,\nu-\mu}\ot _{\OO_\PP} 
(\MM^\PP_{\chi,\mu}\ten_{\OO_{\ZZ(\tii \DD_\PP)}}\FF))]
\\
\cong
\RGa [ (\pi^\BB_\PP)^*
(\MM^\PP_{\chi,\nu}\ten_{\OO_{\ZZ(\tii \DD_\PP)}}\FF)]
\cong 
\RGa   
(\MM^\PP_{\chi,\nu}\ten_{\OO_{\ZZ(\tii \DD_\PP)}}\FF)
\cong
\Rr\Ga [  
(\tii\pi\tw)^*(\MM^\QQ_{\chi,\nu})\ten_{\OO_{\ZZ(\tii \DD_\PP)}}\FF]
\\ \cong
\RGa [  
\MM^\QQ_{\chi,\nu}\ten_{\OO_{\ZZ(\tii \DD_\QQ)}}
\tii \pi_*\tw\FF].
\end{multline*}

Here the key step is the second isomorphism, which is provided
by Lemma \ref{Up and down}.a. The first one is clear from 
the projection formula and the well-known fact that 
$\Rr(\pi^\BB_\PP)_*\OO=\OO$.
The third isomorphism is a triviality (notice that $\nu-\mu$
is a character of the Levi of $P$, so $\OO_{\PP,\nu-\mu}$
is defined). The fourth one 
 follows from $\MM^\PP_{\chi,\nu}=
\OO_{\PP,\nu-\mu}\ot _{\OO_\PP} 
\MM^\PP_{\chi,\mu}$ which is clear from the construction of these
bundles in Remark \ref{Spring is there}. The fifth one again
follows from 
$\Rr(\pi^\BB_\PP)_*\OO=\OO$. The sixth one is obtained by
substituting \eqref{PQcomp}
from Remark \ref{Spring is there}. The last one is the projection formula.

This proves the first isomorphism in the lemma.
The second isomorphism follows since $T_{\mu}^\nu$ is both left and
right adjoint to
 $T_{\nu}^\mu $, while $L\tii \pi^*$ is the left
adjoint of $\Rr\tii \pi_*$ and $L\tii\pi^!$
is the right adjoint of $\Rr\tii \pi_*$.  \epf

\srems 1) A variation of a well-known argument \cite[2.5]{BeGi}
can be used to
show that 
the functor $T_\mu^\nu$ sends the full subcategory 
 $\mod(\U^\mu)\subset \mod_\mu(\U)$ to $\mod(\U^\nu)\subset
\mod_\nu (\U)$ and identifies  $\mod(\U^\nu)$ with a Serre 
quotient category of  $\mod(\U^\mu)$. It follows that
 $T_\mu^\nu$ sends
an irreducible object to zero or to an
irreducible object, and that every irreducible object in $\mod(\U_\chi^\nu)$
is isomorphic to $T_\mu^\nu(L)$ for a unique irreducible object
in  $\mod(\U_\chi^\mu)$. By adjunction, the  functor
 $T_\nu^\mu$ sends the projective cover of an irreducible
object  in $\mod_{\nu,\chi}(\U)$
 to the projective cover of the corresponding irreducible
in  $\mod_{\mu,\chi}(\U)$.
 (Here the projective
cover is understood to be a pro-object in $ \mod_{\nu,\chi}(\U)$.)

Together with the previous lemma  this provides an effective
tool for computation of coherent sheaves corresponding to some
$\g$ modules of interest.

2) The isomorphism  $L\tii \pi^*\cong L\tii \pi^!$
which was just deduced from the properties of translation functors
also follows easily from the fact that $\gt^*_\PP$ and $\gt^*_\QQ$ 
have trivial
canonical classes and their dimensions are equal.

\sss{Composition of translation functors} We need some more
information on the behavior of translation functors. Recall the
functors \fttt{Notice that these are {\em not} the functors used
above to get the localization
 theorem: the latter are direct summands of the former.}
$\RGa:\Db[\mod^c(\tii \DD)]\to
\Db[\mod^{fg}(\U)]$, $\LL:
\Db[\mod^{fg}(\U)] \to
\Db[\mod^c(\tii \DD)]$, $\LL:M\to \tii \DD
\Lotimes _{\U}M$. The functor $\RGa$ can be "upgraded" to the derived functor
of global sections
$\Rr\tii \Ga:\Db[\mod^c(\tii \DD)]\to
\Db[\mod^{fg}(\tii \U)]$.
% here
%$\tii \U=\U\otimes _{\zhc} S(\fh)= \Ga (\tii \DD)$ and 
%$\Rr\tii
%\Ga$ is the derived functor of global sections.

For $\la\in \fh^*$ let
$\mod_\la(\tii \U)\subset
\mod(\tii \U)$ be the
subcategory of modules where the central subalgebra $S(\fh)\subset
\tii \U$ acts through the quotient by a power of the maximal ideal
defined by $\la$. Then 
$\Rr\tii \Ga$ restricts to 
$\Rr\tii \Ga_\la:
\Db[\mod_\la(\tii \DD)]\to \Db[\mod_\la(\tii \U)]$. Let also
$\U^{\hatt\la}$, $\tii \U^{\hatt \la}$ denote the completions of $\U$,
$\tii \U$ at the corresponding character of $\zhc$, $S(\fh)$
respectively; thus $\mod_\la(\U)$,  $\mod_\la(\tii \U)$ are identified
with full subcategories in, respectively, $\mod(\U^{\hatt\la})$, $\tii
\mod(\U^{\hatt \la})$.

\slem\label{tilde_Ga} {\em   For an integral $\la\in \h^*_{\Fp}$ we
have canonical isomorphisms
%an isomorphism of endo-functors
%of $\Db[\mod(\U)]$:
%is  isomorphic to a direct summand of the composition
$$
\Rr\tii \Ga \circ \LL\cong \Ind_\U^{\tii \U},$$
$$
\Rr\tii \Ga \circ \LL^{\hatt \la}\cong \Ind_{\U^{\hatt \la}}^{\tii \U^{\hatt
\la}}.$$
 }

\pf The isomorphisms of functors on derived categories follow from
the corresponding isomorphisms of functors on the categories of
bounded complexes of projective objects. The latter are established
directly, cf. \cite{BMR},
Corollary 3.4.2.
\epf

\spro \lab{properties of translation functors prop} {\em If $\mu$
lies in the closure of the $\la$ facet, and $\nu$ lies in the
closure of the $\mu$-facet, then
 there is a canonical isomorphism
of functors
$T_{\mu}^\nu \circ T_{\la}^\mu\cong
T_{\la}^\nu$
and also
$T_{\mu}^\la \circ T_{\nu}^\mu\cong
T_{\nu}^\la$.

In particular, if $\mu$ and $\nu$ are in the same facet then
$T_{\mu}^\nu$ and $T_{\nu}^\mu$ are inverse equivalences.}

\pf We construct the first isomorphism, the second one follows
by adjointness. Lemma
 \ref{Up and down}.a implies that
%$$
\begin{equation}\label{TTG}
T_{\mu}^\nu \circ T_{\la}^\mu\circ \RGa_{\tii\DD_\BB, \la} \cong
T_{\la}^\nu\circ \RGa_{\tii\DD_\BB, \la}
%$$
\end{equation}
 canonically. By Lemma \ref{tilde_Ga}
this implies that
$$T_{\mu}^\nu \circ T_{\la}^\mu\circ \Res^{\tii \U^{\hatt \la}}_{\U^{\hatt
\la}}
 \Ind^{\tii \U^{\hatt \la}}_{\U^{\hatt \la}} \cong
T_{\la}^\nu\circ  \Res^{\tii \U^{\hatt \la}}_{\U^{\hatt \la}}
 \Ind^{\tii \U^{\hatt \la}}_{\U^{\hatt \la}}.$$

The assumption $p>h$ implies that $p$ does not divide $|W|$, thus the subspace
of $W$ invariants in $S(\fh)$ is  canonically a direct summand
as a module for $S(\fh)^W$.
Similarly, the completion $S(\fh)^{\hatt \la}$ contains
$(S(\fh)^W)^{\hatt \la}$ as a canonical direct summand. Thus
$T_{\mu}^\nu \circ T_{\la}^\mu$, $T_\nu^\mu$ are canonical direct
summands in the left hand side (respectively, right hand side) of
the last isomorphism. This yields a morphism $T_{\mu}^\nu \circ
T_{\la}^\mu\to T_\la^\nu$ defined as a composition of the embedding
of a direct summand with the last isomorphism and the projection to
a direct summand.

 Applying 
both  $T_{\mu}^\nu \circ
T_{\la}^\mu$ and $ T_\la^\nu$
 to the object $\U^\la=
R\Gamma_{\tii \DD_{\BB,\la}} (\DD^\la)$ we get 
isomorphic objects by \eqref{TTG}. 
It is straightforward to check that this isomorphism 
coincides with the map coming from the morphism of functors constructed
in the previous paragraph.
 Since both functors are exact,
and $\U^\la$ generates $\mod^\la(\U)$ under extensions and taking
cokernels, the desired isomorphism follows. \epf

\cor
{\em
Assume that $\mu$ is regular, and $\nu$, $\nu'$ lie
on the same codimension one face of the alcove of $\mu$.
Then we have a canonical isomorphism $R_{\mu|\nu}\cong R_{\mu|\nu'}$.
\epf}

The corollary allows one to define
a reflection functor
$R_\alpha:\Db[\mod^{fg}_\la(\U)]\to \Db[\mod^{fg}_\la(\U)]$
for
any  $\alpha\in \LI\aff$,
as follows.
Recall that $\alpha\in \LI\aff$ is an orbit of $\Wap$
on the set of pairs
$(\lambda, F)$ with $\la\in \Theta$ and $F$ a face of the $\la$-alcove.
Choose a representative
$(\lambda, F)$ of $\al$ and a character $\nu \in F$.
We then set
$R_\alpha=R_{\lambda|\nu}$. The corollary
 shows that $R_\alpha$ is independent
 of these choices. Similarly, we get functors $\bI_\al^*$, $\bI_\al^!$,
$\al\in \LI\aff$.

\pro
{\em  For $\alpha \in \LI\aff$ and $\mu\in \Theta$ such that
$\mud=\mu\ast s_\al>\mu$
we have $\LL^{\hatt\mu}\circ \bI_\alpha^*\cong \LL^{\hatt{\mud}}$.}

{\em Proof.}\ This is a restatement of Lemma \ref{Up and down}.c. \epf

%{\em Proof.}  $\tI_{\al}^*$ is invertible in view of the Proposition,
%and it's easy to see that  $\tI_{\al}^*$ and $\tI_{\al}^!$ are adjoint.
%\epf

\sus{Affine braid group actions: proofs} \sss{Proof of Theorem
\ref{braid_act_mod}}\lab{proof_of_act_mod} It suffices to check
that:

(1) If $\omega\in \LB'\aff$, $\ell(\omega)=0$, then  for all $\la
\in \Theta$ we have a canonical isomorphism $\LL^{\hatt{\la\ast
\om}}\cong \LL^{\hatt\la} \circ T_\omega$ (see the first paragraph
of section \ref{Tranintert} for the notation).

(2) Let $\al\in \LI\aff$ and let
$s=s_\al\in \LB\aff$ be the corresponding  simple
reflection.
If $s$ increases  $\la\in \Theta$, then
$\LL^{\hatt{\la\ast s}}\cong \LL^{\hatt\la} \circ \bI_\al^*$.

(3) The map $\tilde\omega\mapsto T_\omega$, $\st_\al\mapsto
\bI_\al^*$ extends to an action of $\LB'\aff$.

We claim  that (3) follows from (1) and (2). Indeed, for $w\in
\LW\aff$ it is easy to find an element $\la\in \Theta$, such that
$w$ increases $\la$. By induction on $\ell(w)$ it follows from (1)
and (2) that for any decomposition $w=\omega s_{\al_1}\dots
s_{\al_{\ell(w)}} $ where $\ell(\omega)=0$ and $\al_i \in \LI\aff$
we have
$$\LL^{\hatt{\mu}}\cong \LL^{\hatt{\la}}\circ T_\om\circ \bI_{\al_1}^*
 \circ \cdots \circ \bI_{\al_{\ell(w)}}^*,$$
where $\mu=\la*\om*s_{\al_1}*\dots *s_{\al_{\ell(w)}}$. This shows
that the composition of functors on the right hand side does not
depend on the reduced decomposition of $w$, which implies (3).

It remains to check (1) and (2).
 In fact, (1) is immediate from  Lemma \ref{Up and down}.a,
while (2) follows from Lemma \ref{Up and down}.c. \epf

\sss{The action of  $B'\aff$ on $\tii \DD$ modules: proof
of Corollary \ref{bbbI_acts}}
\label{action on D}
It suffices to check the following: if $\la$, $\mu$ lie in the fundamental
alcove $A_0$, then the translation functor $T_\la^\mu$ intertwines
the action of $^{\Theta(\la)}B'\aff$ with the action of $^{\Theta(\mu)}
B'\aff$; here we have identified
 $^{\Theta(\la)}B'\aff$ with $^{\Theta(\mu)}
B'\aff$ by means of the composed isomorphism:
 $$^{\Theta(\la)}B'\aff
@>{{\bf isom}_\la^{-1}}>> B'\aff@>{{\bf isom}_\mu}>>
\  {^{\Theta(\mu)}
B'\aff}.$$
This compatibility follows from the definitions and Proposition
\ref{properties of translation functors prop}. \epf

\srem The action of $B'\aff$ on  $\Db[\mod^c_\la(\tii \DD)]$,
and a closely related action on
$\Db[\Coh(\gt^*)]$, $\Db[\Coh(\Nt)]$ will be discussed in more detail in
a future publication.
There we will provide a geometric description of the action of
generators in
terms of convolution with some explicit sheaves on $(\gt^*)^2$,
analogous to the one for the characteristic zero setting given in
\cite{BeBe}.
Also, we will show that the induced action on the
K-groups of Springer fibers factors through the affine Weyl group and
extends to an action on the equivariant K-groups, where it
gives the action of the affine Hecke algebra on a standard module.

 Here we
state only a simple property of the action
on  $\Db[\mod^c_\la(\tii \DD)]$.

\pro {\em Fix a regular $\la \in \La$. For $\nu\in
\Lambda^+\subset \Wap$ on $\Db(\mod^c_\la(\tii \DD))$ the action
of $\bbI_{\tilde \nu}$ is given by $\bbI_{\tii \nu}:\F\mapsto
\F\otimes_{\OO_{\BB\tw}} \OO_\nu\tw$.
 }

\pf Pick $\la\in A_0$. Let $g_\nu$
be  the image of $\nu\in \Lambda^+\subset \Wap$
under the isomorphism $\isla$; thus $\la*g_\nu=\la+p\nu$ and $g_\nu$
 increases $\la$. Then for
$\MM\in \mod^c_\la(\tii\DD)$ we have
$\bbI_{\tii{\nu}}(\MM)=\LL(\bI_ {g_\nu}(R\Ga(\MM)))$
(cf. the characterizations of the actions in Corollary \ref{bbbI_acts}).
 On the other hand we have
$$\LL(\bI_ {g_\nu}(R\Ga(\MM)))\cong {\mathrm{Tens}}_{\la+p\nu,\la}\LL(R\Ga(\MM))
 = (\OO_{\BB\tw})_\nu\ot_{\OO_{\BB\tw}}  \MM
,$$
where the first isomorphism is a particular case of \eqref{Tens} from Theorem
 \ref{braid_act_mod},
and the second one follows from the definition of ${\mathrm{Tens}}_{\la,\mu}$.
\epf

\se{\bf Dualities}

In this section we translate the duality operation on $\fg$-modules
to $\DD$-modules and coherent sheaves. This section is directly
inspired by Lusztig \cite[part II, remark ]{Lu} (see remark \ref{onLu} below).

\setcounter{subsubsection}{0}
\sss{
Homological duality
}
For a $\fg$-module $M$ the dual $\k$ vector space $M^*$ carries a
$\fg$ action, we denote the resulting $\fg$-module by $M\check{\ }$.
Thus we get a contravariant equivalence between finite dimensional
$\U_{\hat\chi}^{\hat\lambda}$-modules and finite dimensional
$\U_{\hatt{-\chi}}^{\hatt{\lambda\check{\ }}}$-modules, which sends
$M$ to $M\check{\ }$; here $\lambda\check{\ }=-w_0(\la)$ is the
{\em dual
weight} (recall that $w_0\in W$ is the long element). We use the same notation $M\mapsto M\check{\ }$
for the anti-equivalence of the derived categories induced by the exact functor
on abelian categories. 

The aim of this
section is to describe the localization of $M\check{\ }$ in terms of
localization of $M$.
We start by extending the duality to $\Db[\mod^{fg}(\U)]$. To this end
recall the {\em homological duality} $\bbD_\U:M\mapsto \RHom_\U(M,
\U)[\dim \fg]$, where the latter is understood to be an object in the
derived category of right $\U$-modules (which can, of course, be
identified with left modules in a standard way).

In the following Lemma we assume that $\fg$ is any finite dimensional
unimodular Lie algebra (i.e. the adjoint action of $\fg$ on its top
exterior power is trivial) over any field $k$ (not of characteristic
2 or 3).

\slem
\label{MdD}
{\em
Let $M$ be a finite complex of finite
dimensional $\fg$-modules. We have a canonical isomorphism in
$\Db[\mod(\U(\fg))]$: $M\check{\ }\cong \bbD_\U(M)$. }

\pf   Let $C^\bu$ be the standard projective resolution for the
trivial $\U(\fg)$-module, thus $C^i=\Lambda^{-i}(\fg)\otimes \U(\fg)$
with Koszul-Chevalley differential. Choosing an isomorphism $k\cong
\La^{\dim \fg}(\fg)$ we get an isomorphism $\La^{\dim \fg - d}(\fg)\cong
\La^d(\fg)^*$, which yields an isomorphism of complexes of $\fg$
modules $\Hom_\U(C^\bu, \U)\cong C^\bu[-\dim \fg]$.
% {\bf is this right?}.

We have isomorphisms in $\Db[\mod(\U(\fg))]$:
 $$\bbD_\U(M)\cong  \Hom_\U(M\otimes C^\bu, \U[\dim \g])\cong M\check{\ } \otimes
 \Hom _\U(C^\bu, \U[\dim \g])\cong  M\check{\ }\otimes C^\bu\cong
 M\check{\ },$$
 which yields the result. \epf

\subsubsection{Duality for $\tii \DD$ modules}
We define the homological duality for $\tii \DD$ modules,
$\bbD_D:\Db[\mod^c(\tii \DD)]\to \Db[\mod^c(\tii \DD^{op})] $,
$\bbD_D(\MM)=\RHHom (M,\tii \DD)[\dim \fg ]$.

Recall the functors
%\fttt{Notice
%that these are {\em not} the functors used above to get the localization
% theorem: the latter are direct summands of the former.}
$\RGa:\Db[\mod^c(\tii \DD)]\to \Db[\mod^{fg}(\U)]$, $\LL:
\Db[\mod^{fg}(\U)] \to \Db[\mod^c(\tii \DD)]$, $\LL:M\to \tii \DD
\Lotimes _{\U}M$. We use the same notation for similarly defined
functors between derived categories of right $\DD$-modules and right
$\U$-modules.

\lem\label{LDU} {\em  We have canonical isomorphisms of functors
$\LL\circ \bbD_\U \cong \bbD_D \circ \LL $, $\RGa\circ \bbD_D\cong
\bbD_\U \circ \RGa$.
 }

\pf  It is clear that we have natural isomorphisms
$$\LL\circ \bbD_\U(\U)\cong \LL (\U[\dim \g])\cong \tii \DD [\dim \g]
 \cong \bbD_D(\tii
\DD)\cong \bbD_D(\LL(\U));$$ moreover, it is easy to see that this
isomorphism is compatible with the action of $\U=\End_\U(\U)$. Thus we
get an isomorphism of the two functors on the category of bounded
complexes of finitely generated projective $\U$-modules induced,
respectively, by $\LL \circ \bbD_\U$ and $\bbD_D \circ \LL$. The
ensuing isomorphism of functors on the derived category is the first
of the two desired isomorphisms.

The functor $\LL$ is left adjoint to $\RGa$. It has been shown
in the first chapter (see Corollary \ref{U is Calabi-Yau}.c)
 that both categories
$\Db[\mod^{fg}(\U)]$ and $\Db[\mod^c(\tii \DD)]$ are Calabi-Yau
categories over $\OO(\fg^*\tw)$. It follows that $\LL$ is also right
adjoint to $\RGa$. In view of the two adjunctions the second
isomorphism follows from the first one. \epf

%\rem This is a particular case of  NC Grothendieck-Serre.

\subsubsection{Right $\DD$-modules as twisted left $\DD$-modules}
Let $X$ be a smooth $d$-dimensional  variety over the field $\k$,
and let $\Omega_X$ be the sheaf of volume forms. We have the Lie
derivative action of vector fields on $\Omega_X$, $\xi: \omega \to
\Lie_\xi(\omega)= d(i_\xi \omega)$ (the term $i_\xi d$ vanishes on top forms).

\slem
{\em
a) Let $\MM$ be a sheaf of $\DD_X$ modules. The action
of vector fields
$\xi:\omega \ot \sigma \mapsto -\Lie_\xi(\omega)\ot \sigma -
\omega \ot \xi(\sigma)$ extends to an action of  $\DD_X^{op}$
on
$\Omega_X\otimes _{\OO_X} \MM$. The resulting functor
$\MM\mm\MM\ten\Om_X$ provides a
(covariant) equivalence between the categories of left and right
$\DD_X$ modules.

b)
The action of $\DD_X^{op}$ on $\Om_X$ (the case $\MM=\OO_X$ in (a)),
identifies  $\DD_X^{op}$ with the algebra
$^{\Om_X}\DD_X=\Om_X\ten\DD_X\ten\Om_X\inv$
of differential operators on
the line bundle $\Om_X$.

c) Let $\omega$ be a nowhere vanishing section of $\Omega_X$. Then
there exists an isomorphism $I_\omega:\DD_X \cong \DD_X^{op}$, such
that the actions of $\DD_X$ and $\DD_X^{op}$
on
$\DD_X\ten \Om_X=\Om_X\ten \DD_X^{op}$
are related by
$d(\omega\otimes 1)=\omega\otimes I_\omega(d)$ for any section
$d$ of $\DD_X$.}

\pf The proof is the same as in
the characteristic zero case, cf., e.g.,  \cite[Chapter VI.3]{Borel}.
\epf

\sss{ Duality and splitting bundles } \lab{Duality and splitting
bundles} We now apply these general considerations to the case of
the flag variety. In the remainder $\nu$ lies in $\La=X^*(T)$  and
$\chi\in\fg^*\tw$ is nilpotent. We will sometimes identify the formal neighborhood 
%  $\FN_{\ZZ(\U)}(\chi,d\nu)$ with $\FN_{\g^*\tw(\chi)$ and
of $\BB_\chi\tw\times \{d\nu\}$
in $\ZZ(\tii \DD_\BB)=\tii \g ^{*}\tw\times_{\h^{*}\tw} \h^*$
 with $\FN_{\tii \g^{*}\tw}(\BB_\chi\tw)$, we have such 
a canonical
identification, since the Artin-Schreier map $\fh^*\to \fh^*\tw$ is etale.

Recall a canonical choice of a
splitting bundle $\MM^{\nu}_\chi$ for the Azumaya algebra $\tii \DD$
on the formal neighborhood of a Springer fiber $\BB_\chi\tw\times
\{d\nu\}\subset \ZZ(\tii \DD)$, 
see \ref{Spring is there}. First,
$\MM^{-\rho}_\chi$ is a pull-back from a formal neighborhood of
$(\chi,-\rho)$ in $\ZZ(U)=\fg^*\tw\times _{\h^*\tw/W}\h/W$, so it is trivial as a vector bundle.
%Since any
%two splitting bundles differ by a line bundle, this is the unique
%splitting bundle for the data $(-\rho,\chi)$ which is trivial as a
%vector bundle.
% I've percented the above b/c it uses that $\Pic$ has no p-torsion,
%and is unnecessary
 The general case is obtained by twisting\ie
$\MM^{\nu}_\chi \dff \MM^{-\rho}_\chi\ten_{\OO_\BB}\OO_{\nu+\rho} $.

The above choice of the splitting bundles $\MM^{\nu}_\chi$
is used  to
construct the equivalence of coherent sheaves and $\fg$-modules:
$$
\ga_{\chi,\nu}(\FF)
\dff
\RGa(\MM^\nu_\chi\ten\FF)
,\ \ \
\FF\in \Db[\Coh_{\BB_\chi}(\gt^*)]
.$$

\lem\label{MMMM} {\em a) We have a canonical isomorphism
%of algebras
%$\DD_{-\rho+\nu}^{op}\cong \DD_{-\rho-\nu}.$
$\fI:\tii \DD\cong \tii \DD^{op}$ compatible with the standard
isomorphism $\U^{op}\cong \U$.
 Its restriction to
 the center coincides with the map induced by the involution
  $\iota$ of $\gt^*\tw \times _{\fh^*\tw} \fh^*$ given by
  $\iota(\fb,x; \nu)=(\fb,-x; -2\rho -\nu)$.

b) There exists an isomorphism of $\tii \DD$ modules
$$
\MM^{-\rho+\nu}_\chi
\cong
\fI^* \bbD_D(
\MM^{-\rho-\nu}_{-\chi})[-\dim \fg]
$$
compatible with the above
isomorphism of algebras. }

\pf (a) The canonical line bundle on $G/J$ is the pull back of the
canonical line bundle on $\BB$. We have a well-known isomorphism
$\Omega_\BB\cong \OO_{-2\rho}$, thus the canonical line bundle of
$G/J$   is $G$-equivariantly trivial. Thus a choice of a
$G$-invariant nonzero section $\om$ is unique up to scaling. By the
previous lemma $\om$ provides a (canonical) isomorphism
$\DD_{G/J}^{op}\cong\DD_{G/J}$, which leads to $\fI$ since $\tii\DD$
is obtained from $\DD_{G/J}$ by taking direct image to $\BB$ and
invariants with respect to $H$. The formula for the action on the
center comes from the fact that $\omega$ transforms by the character
$-2\rho$ under the action of the abstract Cartan group $H$.

It remains to  check (b). For $\nu=0$ the isomorphism follows from
the fact that the Azumaya algebra $\tii \DD|_{\FN(\BB_{\chi}\times
\{\-\rho\})}$ and its two splitting bundles $\MM^{-\rho}_{\chi}$ and
$\fI^*\bbD_D(\MM_{-\chi}^{-\rho})[-\dim \g]\cong \fI^*
(\MM^{-\rho}_{-\chi})^*$ (where the last $*$ denotes the dual vector
bundle) are pull-backs %under $\pi_\PP$ 
of the corresponding
structures on $\FN_{\g^*\tw\times _{\h^*\tw/W}\h^*/W}(\chi,-\rho)$,
while any two splitting bundles for an Azumaya algebra over the ring
of Taylor series are isomorphic.

 For a general $\nu$ the isomorphism follows from $\MM^{-\rho+
\nu}_\chi=\OO_{\nu} \otimes _{\OO_\BB}\MM^{-\rho}_\chi$ (see
\ref{Duality and splitting bundles}).
 \epf

\sss{ Grothendieck-Serre duality $\bbD_\OO$ } The canonical line
bundle for $\gt^*$ is trivial, and then the same is true  for
$\gt^*\tw\times_{\fh^*\tw} \fh^*$. Therefore, the Grothendieck-Serre
duality functor for $\Db[\Coh(\gt^*\tw\times _{\fh^*\tw} \fh^*)]$ is
given by $\bbD_\OO:\F\mapsto \RHHom (\F,\OO)[\dim \fg]$. We
let $\sigma$ denote the involution of $\gt^*$ sending $(\b, \chi)$
to $(\b, -\chi)$.

\cor
\label{DFM}
{\em
For $\F\in \Db[\Coh_{\BB_\chi}(\gt^*)]$ we have a
canonical isomorphism
$$
\bbD_D(\F\otimes \MM_\chi^{-\rho+\nu})
\cong
\sigma^*\bbD_\OO (\F)\otimes
\MM^{-\rho-\nu}_{-\chi}.$$
}

\pf
For a coherent sheaf $\F$ on
$\gt^*\tw\times_{\fh^*\tw}\fh^*$
and an
object $\MM\in \mod^c(\tii \DD)$ the tensor product sheaf
$\F\otimes_{\OO}\MM$ on $\gt^*\tw\times_{\fh^*\tw}\fh^*$
carries a $\tii \DD$
action. It is obvious that for a locally free sheaf $\F$ and $\MM\in
\Db[\mod^c(\tii\DD)]$ we have $\bbD_D(\F\otimes \MM)\cong \F^*\otimes
\bbD_D(\MM)$ canonically. Extending this isomorphism to bounded
complexes of locally free coherent sheaves we get a canonical
isomorphism $$\bbD_D(\F\otimes \MM)\cong \bbD_\OO(\F) \otimes
\bbD_D(\MM)[-\dim \fg].$$ We restrict attention to the formal
neighborhood of $\BB_\chi\times \{\nu\}$ and plug in
 $\MM=\MM^{-\rho+\nu}_\chi$. The desired isomorphism follows then from Lemma
 \ref{MMMM}.b. \epf

We are ready to prove
% the main result of this section.

\pro
{\em
Let $\nu\in\La$ be  regular
and $\chi\in\fg^*\tw$ be nilpotent.
Then for $\F$ in $
\Db[\Coh_{\BB_\chi}(\gt^*)]
$ we have a canonical
isomorphism
$$
\left[ \ga_{\chi, -\rho+\nu }(\F)\right] \check{\ } \cong
\ga_{-\chi,-\rho-\nu}(\sigma^* \bbD _\OO (\F))
.$$
}

\pf Compare Lemma \ref{LDU} with Corollary \ref{DFM} and Lemma
\ref{MdD}. \epf

\sss{ } In the assumptions of the Proposition, assume moreover that
$\nu\in\La$ lies in the fundamental alcove $A_0$. Recall that for
such $\nu$ there is an isomorphism $ {\bf isom}_\nu:B'\aff\to
\LB'\aff$, where $\Theta$ is the $\Wap$ orbit of $\nu$.
Denote  $\bw_0\dff {\bf isom}_\nu(w_0)$. 
Let $\nu\check{\ }=-w_0(\nu)$ be the dual weight.

\cor
{\em
For $\F$ in $\Db[\Coh_{\BB_\chi\tw}(\gt^*\tw)]$ we have a canonical
isomorphism
$$
\left[ \ga_{\chi, \nu }(\F)\right] \check{\ } \cong
\bI_{\tii{\bw_0}}^{-1}\circ
\ga_{-\chi,\nu\check{\ }}\circ \bbD_\OO\circ \sigma^*(\F)
.$$
}

\pf In view of the
last Proposition we have
$$
\left[ \ga_{\chi, \nu }(\F)\right] \check{\ } \cong
\ga _{-\chi,-\nu-2\rho}(\sigma^*\bbD_\OO(\F)).$$
We have $w_0\bu(-\nu-2\rho)=\nu\check{\ }=(-\nu-2\rho)*\bw_0$.
It is clear that $\bw_0$ increases the weight $-\nu-2\rho$.
Thus 
 Theorem \ref{braid_act_mod} and the definition of $\gamma_{\chi,\nu}$
show that $$\gamma_{-\chi,\nu\check{\ }}\cong \bI_{\bw_0}\circ \gamma_{-\chi,
-\nu-2\rho}.$$
The two displayed isomorphisms yield the result.
\epf

\rem \label{onLu} This description of duality on the category of
modules (especially the particular case $\nu=\nu\check{\ }=0$)
should be compared to the definition of a certain involution
on the $K$-group of a Springer fiber in \cite[part II]{Lu}.
%The latter is
%defined as a composition of three involutions
This remark will be elaborated in a future paper, where we plan to
use the last Corollary to show that Lusztig's involution preserves
the classes of irreducible modules (lifted in  a certain canonical
way to the equivariant $K$-group); cf. also \cite[2.13,2.17]{ICM}.

\end{document}